\documentclass[reqno, 11pt]{amsart}
\pdfoutput=1
\makeatletter % changes the catcode of @ to 11
\let\origsection=\section \def\section{\@ifstar{\origsection*}{\mysection}}
\def\mysection{\@startsection{section}{1}\z@{.7\linespacing\@plus\linespacing}{.5\linespacing}{\normalfont\scshape\centering\S}}
\makeatother % changes the catcode of @ back to 12

\usepackage{amsmath,amssymb,amsthm}
\usepackage{mathrsfs}
\usepackage{mathabx}\changenotsign
\usepackage{dsfont}
\usepackage{multicol}

\usepackage{xcolor}
\usepackage[backref]{hyperref}
\hypersetup{
    colorlinks,
    linkcolor={red!60!black},
    citecolor={green!60!black},
    urlcolor={blue!60!black}
}

\usepackage[open,openlevel=2,atend]{bookmark}

\usepackage[abbrev,msc-links,backrefs]{amsrefs}
\usepackage{doi}

\renewcommand{\PrintDOI}[1]{\doi{#1}}

\usepackage[T1]{fontenc}
\usepackage{lmodern}
\usepackage[babel]{microtype}
\usepackage[english]{babel}

\usepackage{geometry}
\numberwithin{equation}{section}
\numberwithin{figure}{section}

\usepackage{enumitem}
\def\rmlabel{\upshape({\itshape \roman*\,})}

\def\Alabel{\upshape({\itshape \Alph*\,})}

\let\polishlcross=\l
\def\l{\ifmmode\ell\else\polishlcross\fi}

\def\qand{\quad\text{and}\quad}
\def\qqand{\qquad\text{and}\qquad}

\let\emptyset=\varnothing
\let\setminus=\smallsetminus

\makeatletter
\def\moverlay{\mathpalette\mov@rlay}
\def\mov@rlay#1#2{\leavevmode\vtop{   \baselineskip\z@skip \lineskiplimit-\maxdimen
   \ialign{\hfil$\m@th#1##$\hfil\cr#2\crcr}}}
\newcommand{\charfusion}[3][\mathord]{
    #1{\ifx#1\mathop\vphantom{#2}\fi
        \mathpalette\mov@rlay{#2\cr#3}
      }
    \ifx#1\mathop\expandafter\displaylimits\fi}
\makeatother

\newcommand{\dcup}{\charfusion[\mathbin]{\cup}{\cdot}}
\newcommand{\bigdcup}{\charfusion[\mathop]{\bigcup}{\cdot}}

\DeclareFontFamily{U}  {MnSymbolC}{}
\DeclareSymbolFont{MnSyC}         {U}  {MnSymbolC}{m}{n}
\DeclareFontShape{U}{MnSymbolC}{m}{n}{
    <-6>  MnSymbolC5
   <6-7>  MnSymbolC6
   <7-8>  MnSymbolC7
   <8-9>  MnSymbolC8
   <9-10> MnSymbolC9
  <10-12> MnSymbolC10
  <12->   MnSymbolC12}{}
\DeclareMathSymbol{\powerset}{\mathord}{MnSyC}{180}

\usepackage{tikz}
\usetikzlibrary{calc,decorations.pathmorphing}
\pgfdeclarelayer{background}
\pgfdeclarelayer{foreground}
\pgfdeclarelayer{front}
\pgfsetlayers{background,main,foreground,front}
\usepackage{wrapfig}

\usepackage{subcaption}
\let\epsilon=\varepsilon

\let\rho=\varrho
\let\theta=\vartheta

\newcommand{\cF}{\mathcal{F}}

\newcommand{\cH}{H}

\newcommand{\cP}{\mathcal{P}}

\newcommand{\cc}{\mathcal{C}_k(n,s)}

\newcommand{\Ex}{\mathrm{Ex}}
\newcommand{\ex}{\mathrm{ex}}

\newcommand{\qedge}[7]{
	
	\ifx\relax#4\relax
	\def\qoffs{0pt}
	\else
	\def\qoffs{#4}
	\fi
	
	\def\qhedge{
		($#1+#3!\qoffs!-90:#2-#3$) --
		($#2+#1!\qoffs!-90:#3-#1$) --
		($#3+#2!\qoffs!-90:#1-#2$) -- cycle}

	\coordinate (12) at ($#1!\qoffs!90:#2$);
	\coordinate (13) at ($#1!\qoffs!-90:#3$);
	\coordinate (23) at ($#2!\qoffs!90:#3$);
	\coordinate (21) at ($#2!\qoffs!-90:#1$);
	\coordinate (31) at ($#3!\qoffs!90:#1$);
	\coordinate (32) at ($#3!\qoffs!-90:#2$);
	
	\def\nqhedge{
		(13) let \p1=($(13)-#1$), \p2=($(12)-#1$) in
		arc[start angle={atan2(\y1,\x1)}, delta angle={atan2(\y2,\x2)-atan2(\y1,\x1)-360*(atan2(\y2,\x2)-atan2(\y1,\x1)>0)}, x radius=\qoffs, y radius=\qoffs] --
		(21) let \p1=($(21)-#2$), \p2=($(23)-#2$) in
		arc[start angle={atan2(\y1,\x1)}, delta angle={atan2(\y2,\x2)-atan2(\y1,\x1)-360*(atan2(\y2,\x2)-atan2(\y1,\x1)>0)}, x radius=\qoffs, y radius=\qoffs] --
		(32) let \p1=($(32)-#3$), \p2=($(31)-#3$) in
		arc[start angle={atan2(\y1,\x1)}, delta angle={atan2(\y2,\x2)-atan2(\y1,\x1)-360*(atan2(\y2,\x2)-atan2(\y1,\x1)>0)}, x radius=\qoffs, y radius=\qoffs] --
		cycle}
	
	\ifx\relax#5\relax
	\def\qlwidth{1pt}
	\else
	\def\qlwidth{#5}
	\fi
	
	\ifx\relax#7\relax
	\fill \nqhedge;
	\else
	\fill[#7]\nqhedge;
	\fi
	
	\ifx\relax#6\relax
	\draw[line width=\qlwidth,rounded corners=\qoffs]\nqhedge;
	\else
	\draw[line width=\qlwidth,#6]\nqhedge;
	\fi
}

\newcommand{\qqedge}[8]{
	
	\ifx\relax#4\relax
	\def\qoffs{0pt}
	\else
	\def\qoffs{#4}
	\fi
	
	\def\qhedge{
		($#1+#3!\qoffs!-90:#2-#3$) --
		($#2+#1!\qoffs!-90:#3-#1$) --
		($#3+#2!\qoffs!-90:#1-#2$) --
		($#2+#1!\qoffs!-90:#3-#1$) --cycle}

	\coordinate (12) at ($#1!\qoffs!90:#2$);
	\coordinate (13) at ($#1!\qoffs!-90:#2$);
	\coordinate (23) at ($#2!\qoffs!90:#3$);
	\coordinate (21) at ($#2!\qoffs!-90:#1$);
	\coordinate (31) at ($#3!\qoffs!90:#2$);
	\coordinate (32) at ($#3!\qoffs!-90:#2$);
	\coordinate (41) at ($#2!\qoffs!-270:#1!#8!(13)$);
	\coordinate (43) at ($#2!\qoffs!270:#3!#8!(31)$);
	
	\def\nqhedge{
		(13) let \p1=($(13)-#1$), \p2=($(12)-#1$) in
		arc[end angle={atan2(\y2,\x2)}, delta angle={-180}, x radius=\qoffs, y radius=\qoffs] --
		(21) let \p1=($(21)-#2$), \p2=($(23)-#2$) in
		arc[start angle={atan2(\y1,\x1)}, delta angle={atan2(\y2,\x2)-atan2(\y1,\x1)-360*(atan2(\y2,\x2)-atan2(\y1,\x1)>0)}, x radius=\qoffs, y radius=\qoffs] --
		(32) let \p1=($(32)-#3$), \p2=($(31)-#3$) in
		arc[start angle={atan2(\y1,\x1)}, delta angle={-180}, x radius=\qoffs, y radius=\qoffs][rounded corners = 1pt] --
		(43)--(41)[sharp corners]--
		cycle}
	
	\ifx\relax#5\relax
	\def\qlwidth{1pt}
	\else
	\def\qlwidth{#5}
	\fi
	
	\ifx\relax#7\relax
	\fill \nqhedge;
	\else
	\fill[#7]\nqhedge;
	\fi
	
	\ifx\relax#6\relax
	\draw[line width=\qlwidth,rounded corners=\qoffs]\nqhedge;
	\else
	\draw[line width=\qlwidth,#6]\nqhedge;
	\fi
}

\newtheoremstyle{note}  {4pt}  {4pt}  {\sl}  {}  {\bfseries}  {.}  {.5em}          {}
\newtheoremstyle{introthms}  {3pt}  {3pt}  {\itshape}  {}  {\bfseries}  {.}  {.5em}          {\thmnote{#3}}
\newtheoremstyle{remark}  {2pt}  {2pt}  {\rm}  {}  {\bfseries}  {.}  {.3em}          {}

\theoremstyle{plain}
\newtheorem{theorem}{Theorem}[section]
\newtheorem{lemma}[theorem]{Lemma}

\newtheorem{prop}[theorem]{Proposition}

\newtheorem{cor}[theorem]{Corollary}
\newtheorem{fact}[theorem]{Fact}
\newtheorem{claim}[theorem]{Claim}

\theoremstyle{note}

\theoremstyle{remark}

\usepackage{accents}

\usepackage{lineno}
\newcommand*\patchAmsMathEnvironmentForLineno[1]{%
\expandafter\let\csname old#1\expandafter\endcsname\csname #1\endcsname
\expandafter\let\csname oldend#1\expandafter\endcsname\csname end#1\endcsname
\renewenvironment{#1}%
{\linenomath\csname old#1\endcsname}%
{\csname oldend#1\endcsname\endlinenomath}}%
\newcommand*\patchBothAmsMathEnvironmentsForLineno[1]{%
\patchAmsMathEnvironmentForLineno{#1}%
\patchAmsMathEnvironmentForLineno{#1*}}%
\AtBeginDocument{%
\patchBothAmsMathEnvironmentsForLineno{equation}%
\patchBothAmsMathEnvironmentsForLineno{align}%
\patchBothAmsMathEnvironmentsForLineno{flalign}%
\patchBothAmsMathEnvironmentsForLineno{alignat}%
\patchBothAmsMathEnvironmentsForLineno{gather}%
\patchBothAmsMathEnvironmentsForLineno{multline}%
}

\def\cc{{C}}
\def\sp{{SP}}
\def\sk{{SK}}

\def\C{\mathcal{C}}

\def\cP{\mathcal{P}}

\usepackage{todonotes}

	\def\c{C_4}
	\def\p{\cP_4}
	\def\m{M_3}

	\def\rr{\textcolor{red!60!black}{R}}
	\def\bb{\textcolor{blue!60!black}{B}}
	\def\gg{\textcolor{green!60!black}{G}}
	
	\def\mr{\textcolor{red!60!black}{M_2^R}}
	\def\mb{\textcolor{blue!60!black}{M_2^B}}
	\def\mg{\textcolor{green!60!black}{M_2^G}}

	\def\ff{\cF_{12}}
	\def\fff{\cF_{123}}
	\def\fr{\cF_{12}^\lhd}
	\def\fl{\cF_{12}^\rhd}
	\def\ffi{\cF_{12}^{\textrm{in}}}
	\def\ffo{\cF_{12}^{\textrm{out}}}
	
	\newcommand{\crr}[1]{\textcolor{red!60!black}{#1}}
	\newcommand{\cbb}[1]{\textcolor{blue!60!black}{#1}}
	\newcommand{\cgg}[1]{\textcolor{green!60!black}{#1}}

\begin{document}
	\title{Tur\'an and Ramsey numbers for $3$-uniform minimal paths of length $4$}
	
	\author[J. Han]{Jie Han}
	
	\address{Department of Mathematics, University of Rhode Island, Kingston, RI, 02881, USA}
	
	\email{\tt jie\_han@uri.edu}
	
	\author[J.~Polcyn]{Joanna Polcyn}
	
	\address{Adam Mickiewicz University,
		Faculty of Mathematics and Computer Science
		ul.~Uniwersytetu Poznañskiego 4,
		61-614 Pozna\'n, Poland}
	
	\email{\tt joaska@amu.edu.pl}
	
	\author[A.~Ruci\'nski]{Andrzej Ruci\'nski}
	
	\address{Adam Mickiewicz University,
		Faculty of Mathematics and Computer Science
		ul.~Uniwersytetu Pozna\'nskiego 4,
		61-614 Pozna\'n, Poland}
	
	\email{\tt rucinski@amu.edu.pl}

	\thanks{The third author supported by the Polish NSC grant 2018/29/B/ST1/00426. Part of the research was done when the first author visited Adam Mickiewicz University.}
	\date{\today}
	
	

	\begin{abstract}
		We determine  Tur\'an numbers for the family of 3-uniform minimal paths of length four \emph{for all $n$}. We also establish the second and third order Tur\'an numbers and use them to compute the corresponding Ramsey numbers for up to four colors.
	\end{abstract}
	
	\maketitle
	
	%%%%%%%%%%%%%%%%%%%%%%%%%%%%%%%%%%%%%%%%%%%%%%%%%%%%%%%%%%%%%%%%%%%%%%
	%																																		%
	%																																		%
	%											Introduction            																%
	%																																		%
	%																																		%
	%%%%%%%%%%%%%%%%%%%%%%%%%%%%%%%%%%%%%%%%%%%%%%%%%%%%%%%%%%%%%%%%%%%%%%
	
\section{Introduction}	
	Tur\'an-type problems concern the maximum number of edges in a (hyper)graph without certain forbidden substructures.
	They are central to extremal combinatorics and have a long and influential history initiated by Tur\'an in 1944 \cite{Tu} who solved the problem for all complete graphs. A few years later Erd\H os and Stone \cite{ES} determined asymptotically the Tur\'an numbers for all non-bipartite graphs.
 Such questions for hypergraphs are, however, notoriously difficult in general, and several natural problems are wide open, most notably Tur\'an's conjecture  for the tetrahedron. And, again, asymptotic results are perhaps  a little easier to obtain. For a comprehensive survey on Tur\'an numbers for hypergraphs see \cite{Keevash}.

Similar stature and even longer history are enjoyed by Ramsey Theory, started by Ramsey's paper \cite{R} and developed in the mid thirties of the twenties century by Erd\H os and Szekeres \cite{ESz}. Here the object of interest is the smallest order of a complete (hyper)graph which, when edge-partitioned into a given number of colors, possesses a desired substructure entirely in one color. When the  substructure is itself complete, an exact solution of this problem is still beyond our reach already for graphs and becomes hopeless for hypergraphs, except for some very small cases.

The two problems are immanently related by a (trivial) observation that if the number of edges of one color exceeds the Tur\'an number for the target substructure, then there is a monochromatic copy of it in that color. However, since one is typically interested in a small number of colors (as we are), the corresponding Tur\'an numbers should also be known  for small number of vertices.

In general, both, Tur\'an and Ramsey problems are more difficult for dense hypergraphs. Consequently, the area of research interest has broadened to include sparser structures like paths and cycles.
In this paper we focus on a particular family of 3-uniform hypergraphs,  minimal paths of length four, for which the Tur\'an numbers have been already determined \emph{for large $n$} in~\cite{FJS}. We compute them \emph{for all $n$} and, consequently, obtain the corresponding Ramsey numbers for up to four colors.

	\subsection{Basic definitions}
For $k\ge2$, a $k$-graph ($k$-uniform hypergraph) is an ordered pair $H=(V,E)$, where $V=V(H)$ is a finite set (of vertices) and $E=E(H)$ is a subset of the set $\binom Vk$ of  $k$-element subsets of $V$ (called edges). If $E=\binom Vk$, we call $H$ \emph{complete} and denote by $K_n^{(k)}$, where $n=|V(H)|$.

	For $k$-graphs $H'$ and $H$ we say that $H'$ is a sub-$k$-graph of $H$ and write $H'\subseteq H$ if $V(H')\subseteq V(H)$ and $E(H')\subseteq E(H)$.
Given a family of $k$-graphs $\mathcal F$, we call a $k$-graph $H$ \emph{$\mathcal F$-free}
	if for all $F \in \mathcal F$ we have  $F \nsubseteq H$, that is, no sub-$k$-graph of $H$ is isomorphic to $F$. Given a family of $k$-graphs $\mathcal F$ and an integer $n\ge1$, the \textit{Tur\'an number for $\mathcal F$ and $n$} is defined as
		\[
			\ex_k(n; \mathcal F)\coloneq \max\{|E(H)|: |V(H)|=n\;\mbox{ and $H$ is
			$\mathcal F$-free}\}.
		\]
	Every $n$-vertex $\mathcal F$-free $k$-graph with exactly $\ex_k(n;\mathcal F)$
	edges is called  \emph{extremal  for~$\mathcal F$}.
	We denote by $\mathrm{Ex}_k(n;\mathcal F)$
	the family of all $n$-vertex $k$-graphs which are extremal for $\mathcal F$.
	In the case when $\mathcal F=\{F\}$, we will often write $\mathrm{ex}_k(n;F)$ for
	$\mathrm{ex}_k(n;\{F\})$ and $\mathrm{Ex}_k(n;F)$ for $\mathrm{Ex}_k(n;\{F\})$.

Let $\mathcal F$ be a family of $k$-graphs and $r\ge2$ be an integer. The \emph{Ramsey number} $R(\mathcal F;r)$ is the smallest integer $n$ such that every $r$-edge-coloring of the complete $k$-graph $K_n^{(k)}$ yields a monochromatic copy of a member of~$\mathcal F$.
The relationship between Tur\'an and Ramsey numbers allured to above is best exemplified by the following implication:
\begin{equation}\label{RamTur}
\frac1r\binom nk>\mathrm{ex}_k(n;\mathcal F)\quad\Rightarrow\quad R(\mathcal F;r)\le n.
\end{equation}
		
	%The Tur\'an problems for $k$-graphs are notoriously difficult in general, and a series of natural problems are wide open (e.g., Tur\'an's conjecture on the Tur\'an number for tetrahedron).
	As mentioned earlier, we shall consider the Tur\'an problem for a special family of $3$-uniform paths. At this point the reader should be alerted that there are several other notions of paths and cycles in $k$-graphs (e.g., Berge, loose, linear, tight) and that authors take a great liberty in using those names (except for Berge).
	In this paper we restrict our attention to  minimal paths and cycles defined as follows.
	
	Given $k, \ell\ge 2$, a $k$-uniform \emph{minimal $\ell$-path} (a.k.a. \emph{loose}) is a $k$-graph with edge set $\{a_0, a_1, \dots, a_{\ell-1}\}$ such that $a_i\cap a_j\neq \emptyset$ if and only if $|i-j|\le 1$, while a \emph{$k$-uniform minimal $\ell$-cycle} is a $k$-graph with edge set $\{a_0, a_1, \dots, a_{\ell-1}\}$ such that $a_i\cap a_j\neq \emptyset$ if and only if $|i-j|\le 1 \pmod {\ell}$. So, minimal paths and cycles form  special subclasses of, resp., Berge paths and cycles (see, e.g., \cite{MV}), with no redundant edge intersections. Put another way, the minimality manifests itself by no vertex belonging to more than two edges.

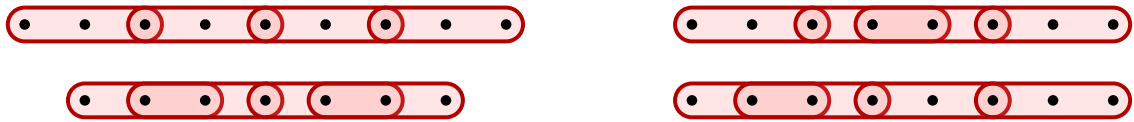
\begin{figure}[h!]
		\begin{multicols}{2}
			\begin{tikzpicture}[scale = .8]
			
			\foreach \i in {1,...,9} {\coordinate (v\i) at (\i,0);}
			
			\qedge{(v3)}{(v2)}{(v1)}{8pt}{1.5pt}{red!70!black}{red!50!white,opacity=0.2};
			\qedge{(v5)}{(v4)}{(v3)}{8pt}{1.5pt}{red!70!black}{red!50!white,opacity=0.2};
			\qedge{(v7)}{(v6)}{(v5)}{8pt}{1.5pt}{red!70!black}{red!50!white,opacity=0.2};
			\qedge{(v9)}{(v8)}{(v7)}{8pt}{1.5pt}{red!70!black}{red!50!white,opacity=0.2};
			
			\foreach \i in {1,...,9} {\fill (v\i) circle (2.5pt);}
			
			\end{tikzpicture}\\
					
			\bigskip
			
				\begin{tikzpicture}[scale = .8]
				
				\qedge{(v1)}{(v2)}{(v3)}{8pt}{1.5pt}{red!70!black}{red!50!white,opacity=0.2};
				\qedge{(v2)}{(v3)}{(v4)}{8pt}{1.5pt}{red!70!black}{red!50!white,opacity=0.2};
				\qedge{(v4)}{(v5)}{(v6)}{8pt}{1.5pt}{red!70!black}{red!50!white,opacity=0.2};
				\qedge{(v5)}{(v6)}{(v7)}{8pt}{1.5pt}{red!70!black}{red!50!white,opacity=0.2};
				
				\foreach \i in {1,...,7} {\fill (v\i) circle (2.5pt);}
				
				\end{tikzpicture}
				
			\begin{tikzpicture}[scale = .8]
			
				\qedge{(v1)}{(v2)}{(v3)}{8pt}{1.5pt}{red!70!black}{red!50!white,opacity=0.2};
				\qedge{(v3)}{(v4)}{(v5)}{8pt}{1.5pt}{red!70!black}{red!50!white,opacity=0.2};
				\qedge{(v4)}{(v5)}{(v6)}{8pt}{1.5pt}{red!70!black}{red!50!white,opacity=0.2};
				\qedge{(v6)}{(v7)}{(v8)}{8pt}{1.5pt}{red!70!black}{red!50!white,opacity=0.2};
				
			\foreach \i in {1,...,8} {\fill (v\i) circle (2.5pt);}
			\end{tikzpicture}\\
			
			\bigskip
			
				\begin{tikzpicture}[scale = .8]
				
				\qedge{(v1)}{(v2)}{(v3)}{8pt}{1.5pt}{red!70!black}{red!50!white,opacity=0.2};
				\qedge{(v2)}{(v3)}{(v4)}{8pt}{1.5pt}{red!70!black}{red!50!white,opacity=0.2};
				\qedge{(v4)}{(v5)}{(v6)}{8pt}{1.5pt}{red!70!black}{red!50!white,opacity=0.2};
				\qedge{(v6)}{(v7)}{(v8)}{8pt}{1.5pt}{red!70!black}{red!50!white,opacity=0.2};
				
				\foreach \i in {1,...,8} {\fill (v\i) circle (2.5pt);}
				\end{tikzpicture}\\
			
			\bigskip

		\end{multicols}
		\caption{All 3-uniform minimal 4-paths from $\p^{(3)}$.}
		\label{fig:cP}
	\end{figure}

	We write $\cP_{\ell}^{(k)}$ for the family of all $k$-uniform minimal $\ell$-paths and $\C_{\ell}^{(k)}$ for the family of all $k$-uniform minimal $\ell$-cycles (see Figure \ref{fig:cP} for all 3-uniform minimal 4-paths). Note that the longest path in $\cP_{\ell}^{(k)}$ has $\ell(k-1)+1$ vertices. It is called \emph{linear} (a.k.a. \emph{loose}), since edges intersect pairwise in at most one vertex, and denoted by $P_{\ell}^{(k)}$.
For convenience, in what follows we shall write $\p$ instead of $\p^{(3)}$.
For $k=2$ the families $\cP_{\ell}^{(2)}$ and $\C_{\ell}^{(2)}$ each consists of a single graph, the ordinary (graph) path and cycle, which will be denoted by, respectively, $P_{\ell}^{(2)}$ and $C_{\ell}^{(2)}$.

\subsection{Main results}
	
	Mubayi and Verstra\"ete \cite{MV} showed that $\ex_k(n;\cP_3^{(k)}) = \binom{n-1}{k-1}$ for all $n\ge2k$ and $\ex_3(n;\cP_\ell^{(3)}) \le\tfrac{5\ell-1}6\binom{n-1}{2}$ for all $n\ge3(\ell+1)/2$. F\"uredi, Jiang, and Seiver \cite{FJS} proved that, for $k\ge3$, $t\ge1$, and \emph{for sufficiently large $n$},
	\begin{equation}\label{43}
	\ex_k(n;\cP_{2t+1}^{(k)})=\binom nk-\binom{n-t}k\quad\mbox{ and }\quad \ex_k(n;\cP_{2t+2}^{(k)})=\binom nk-\binom{n-t}k+1,
	\end{equation}
	and that the unique extremal $k$-graph consists of all $k$-tuples intersecting a given set $T$ of $t$ vertices plus, for even $\ell$, one extra edge disjoint from $T$. Note that for $t=1$, the above expressions become, respectively, $\binom{n-1}{k-1}$ and $\binom{n-1}{k-1}+1$.

In fact, in \cite{FJS} the authors focused on linear paths and determined Tur\'an numbers $\ex_k(n;P_\ell^{(k)}) $ for large $n$ and $k\ge4$, while Kostochka, Mubayi, and Verstra\"ete \cite{KMV} did the same for large $n$ and $\ell\ge4$. The remaining case of $\ell=k=3$ was also  implicit in their proof, but again for large $n$. In \cite{JPR}, it was proved for all $n$ that $\ex_3(n;P_3^{(3)})=\binom{n-1}2$.

	In this paper we similarly extend \eqref{43} in the smallest open case, that is,   we determine the Tur\'an numbers $\ex_3(n;\p)$ for \emph{all~$n$}. All special 3-graphs appearing in  Theorem \ref{P4} below, as well as in Theorems \ref{th:hm}-\ref{e3} in Subsection \ref{ss:turan}, are defined, for clarity of exposition, only in Section \ref{useful}.

		\begin{theorem}\label{P4}
		For $n\ge 1$,
		\[
		\ex_3(n;\p)=\left\{ \begin{array}{ll}
		{n\choose 3} &\text{\rm and} \quad  \Ex_3(n;\p)=\{K_n\}\hskip 2.8cm \text{\rm for } n\le 6,\\
		20  & \text{\rm and}  \quad \Ex_3(n;\p)=\{K_6^{(3)} \cup K_1\}\hskip 1.6cm \text{\rm for $n=7$},\\
		22  & \text{\rm and}  \quad \Ex_3(n;\p)=\{S_8^{+1},SP_8,SK_8\}\hskip .9cm \text{\rm for $n=8$},\\
		\binom{n-1}{2}+1&\text{\rm and}\quad\Ex_3(n;\p)=\{S_n^{+1}\}\hskip 2.7cm \text{\rm for }n\ge 9.
		\end{array} \right.
		\]
	\end{theorem}
	\noindent (Note that for $n=8$, we have $\binom{n-1}{2}+1=22$.)
	%Let $\mathcal F$ be a family of $k$-uniform hypergraphs and. $r\ge2$ be an integer. The \emph{Ramsey number} $R(\mathcal F;r)$ is the smallest integer $n$ such that every $r$-edge-coloring of the complete $k$-uniform hypergraph $K_n^{(k)}$ yields a monochromatic copy of a member of~$\mathcal F$. Our second main result is the following.

As an immediate consequence of Theorem  \ref{P4} and the relation \eqref{RamTur}, plugging $n=3r+1$, we infer that, for $r\ge3$, $R(\p;r)\le 3r+1$. On the other hand, a simple construction originated in \cite{GR} (see Section \ref{RamNum} for more details) yields a lower bound $R(\p;r)\ge r+6$ for all $r\ge1$. Using Theorem \ref{P4} along with some more technical results from the next subsection, we confirm that, at least for up to four colors, the lower bound is, indeed, the correct value.

	\begin{theorem}\label{RN}
		For $r\le 4$, we have $R(\p;r)= r+6$.
	\end{theorem}

\subsection{Tur\'an numbers of higher orders}\label{ss:turan}
	
	To calculate Ramsey numbers based on Tur\'an numbers, it is sometimes necessary to consider Tur\'an numbers of higher orders (see, e.g., \cite{JPR2}), which can be defined iteratively as follows. The Tur\'an number of the first order is the ordinary Tur\'an number.
	For a family of $k$-graphs $\mathcal F$ and  integers $s,n\ge1$,
	\textit{the Tur\'an number of the $(s+1)$-st order} is defined as
	\begin{eqnarray*}
		\mathrm{ex}^{(s+1)}_k(n;\mathcal F)=\max\{|E(H)|:|V(H)|=n,\; \mbox{$H$ is
			$\mathcal F$-free, and }\\
		\forall H'\in \mathrm{Ex}^{(1)}_k(n;\mathcal F)\cup...\cup\mathrm{Ex}^{(s)}_k(n;\mathcal F),  H\nsubseteq H'\},
	\end{eqnarray*}
	if such a $k$-graph $H$ exists.
	An $n$-vertex $\mathcal F$-free $k$-graph $H$ is called \textit{(s+1)-extremal for} $\mathcal F$ if $|E(H)| = \mathrm{ex}^{(s+1)}_k(n;\mathcal F)$ and  $\forall H'\in \mathrm{Ex}^{(1)}_k(n;\mathcal F)\cup...\cup\mathrm{Ex}^{(s)}_k(n;\mathcal F),  H\nsubseteq H'$; we denote by $\mathrm{Ex}^{(s+1)}_k(n;\mathcal F)$ the family of $n$-vertex $k$-graphs which are $(s+1)$-extremal for $\mathcal F$.
	
	A historically first example of a Tur\'an number of the 2nd order is due to Hilton and Milner \cite{HM} who determined the maximum size of \emph{a nontrivial} intersecting $k$-graph, that is, one which is not a star (see the definition in Section \ref{useful}). Recall that a 3-graph is intersecting if and only if it is $M_2$-free and that, by Erd\H os-Ko-Rado theorem \cite{EKR}, $\ex(n,M_2)=\binom{n-1}2$ for $n\ge6$, while for $n\ge7$ the only extremal 3-graph is a full star. Hilton and Milner proved that
		\begin{theorem}{\cite{HM}}\label{th:hm}
			For $n\ge 7$ we have $\ex_3^{(2)}(n;M_2)=3n-8$.
		\end{theorem}
\noindent In \cite{HK} the authors determined $\ex_k^{(3)}(n;M^{(k)}_2)$ for all $k$; in \cite{PR2}  the complete hierarchy of 3-uniform Tur\'an numbers $\ex_3^{(s)}(n;M_2)$, $s=1,\dots,6$, has been found (for $s\ge7$ they do not exist).
		
	In this paper we determine for $\p$ the Tur\'an numbers  of the second and third order.
	
	\begin{theorem}\label{e2}For $n\ge9$,
		\[
		\ex_3^{(2)}(n;\p)=\left\{ \begin{array}{ll}5n-18&\text{\rm and} \quad  \Ex_3^{(2)}(n;\p)=\{SP_n\}\hskip 1.1cm\text{\rm for $n\le 11$},\\
		\binom{n-3}2+7&\text{\rm and} \quad  \Ex_3^{(2)}(n;\p)=\{CB_n\}\hskip 1cm\text{\rm for $n\ge 12$}.
		\end{array} \right.
		\]
	\end{theorem}
	
	\begin{theorem}\label{e3} For $n\ge9$,
		\[
		\ex_3^{(3)}(n;\p)=\left\{ \begin{array}{ll}4n-10&\text{\rm and} \quad  \Ex_3^{(3)}(n;\p)=\{SK_n\}\hskip 1.1cm\text{\rm for $n\le 10$},\\
		\binom{n-3}2+7=35&\text{\rm and} \quad  \Ex_3^{(3)}(n;\p)=\{CB_{n}\}\hskip 1.1cm\text{\rm for $n=11$},\\
		5n-18=42&\text{\rm and} \quad  \Ex_3^{(3)}(n;\p)=\{SP_{n}\}\hskip 1.2cm\text{\rm for $n=12$},\\
		47&\text{\rm and} \quad  \Ex_3^{(3)}(n;\p)=\{SP_{n}, B_n\}\hskip 0.5cm\text{\rm for $n=13$},\\
		\binom{n-4}2+11&\text{\rm and} \quad  \Ex_3^{(3)}(n;\p)=\{B_n\}\hskip 1.4cm\text{\rm for $n\ge 14$}.
		\end{array} \right.
		\]
	\end{theorem}
	\noindent  Note that for $n=13$ we have $5n-18=\binom{n-4}2+11=47$.
	
\subsection{Notation}
For a $k$-graph $\cH$ and a vertex $v\in V(\cH)$, the \emph{link graph}  of $v$ in $\cH$ is the $(k-1)$-graph on the vertex set $V(\cH)$ and the edge set
		\[
			L_{\cH}(v) = \{e\setminus \{v\}: v\in e\in \cH\}.
		\]
\emph{The degree of $v$} in $\cH$ is defined as $\deg_{\cH}(v)=|L_{\cH}(v)|$, while  maximum and  minimum degrees in $\cH$ are denoted by $\Delta_1(\cH)$ and $\delta_1(\cH)$, respectively. For $k=2$, we obtain the ordinary notions of degrees and maximum and minimum degree in a graph. Also, in the case $k=2$, the link graph is just a set of singletons and coincides with the standard notion of the neighborhood $N_G(v)$.
The subscript $_1$ in $\Delta_1(\cH)$ and $\delta_1(\cH)$ is often omitted.

For a 3-graph $\cH$ on $V$, the set of neighbors of a pair $x,y\in V$ in $\cH$ is defined as 
$$N_{\cH}(x,y)=\{z: \{x,y,z\}\in \cH\}.$$
The number  $\deg_{\cH}(x,y) = |N_{\cH}(x,y)|$ is called  \emph{degree of the pair of vertices $x,y$}  and we set $\Delta_2(\cH)=\max_{x,y\in V}\deg_{\cH}(x,y)$ for the \emph{maximum pair degree} in $\cH$.

	We identify a $k$-graph $\cH$ with its edge set $E(\cH)$.
	Throughout the paper we will use the name ``edge'' for both, the edges of a 3-graph (triples) and the edges of a 2-graph (pairs). It will always be clear from the context which one is meant.
For a $k$-graph $\cH$ with vertex set $V$ we write  
		
		\[
			V[\cH] \coloneq  \bigcup_{h\in \cH}h
		\]
for the set of all non-isolated vertices, i.e., vertices $v$ with $\deg_{\cH}(v)>0$.
Given $W\subseteq V$ we write
\[
			\cH[W] \coloneq  \{h\in \cH\colon h\subseteq W\}
		\]
for the sub-$k$-graph of $\cH$ induced by $W$.
		
	For simplicity, if there is no danger of confusion,
	we sometimes denote edges $\{x,y\}$ of graphs and edges $\{x,y,z\}$
	of $3$-graphs by $xy$ and $xyz$, respectively.
	Also, if $f=\{x,y\}$ is a pair of vertices and $v\in V$ is a single vertex, we may write $fv$ for the edge $\{x,y,v\}\in \cH$.

	  Notation $f_1f_2\cdots f_\ell$ will represent a minimal path with edges $f_1,f_2,\dots, f_\ell$ in this order and, likewise, notation $v_1v_2\cdots v_m$ will represent a minimal path with vertices $v_1,v_2,\dots,v_m$ in this order. The same shorthand notation may apply to cycles as well.

%On the other hand, we sometimes denote a minimal 4-path, a path, or a cycle by listing their vertices in the order.

	For two $k$-graphs $G, \cH$, let $G\cup \cH$ denote the disjoint union of them.
	If $\cH$ is a $k$-graph on $V$, $v\in V$, and $e\in \cH$ is an edge of $\cH$, then
	we denote by $\cH-v$ the $k$-graph obtained from $\cH$ by deleting vertex $v$ together with all edges containing it, whereas by $\cH-e$ we mean the $k$-graph obtained from $\cH$ by deleting the single edge $e$.
	For a  $k$-graph $\cH$, by $\cH^c$ we mean the complement of $\cH$, that is, $\cH^c=\binom Vk\setminus \cH$.
	
		%
	
%	For given $t\ge 3$ and the vertex set $V$, $|V|=t$, we write $K_t^{(k)}[V]$ for the complete $k$-graph on the vertex set $V$. Having given two sets of vertices $S$, $|S|=s\ge 2$, and $T$, $|T|=t\ge 2$, by $K^{(2)}_{s,t}[S\dcup T]$ we denote a complete bipartite 2-graph with the bipartition $S\dcup T$.
%Given $k$-graphs $H$ and $F$ and vertex sets $S\subseteq V(H)$, when the host graph $H$ is clear, denote by $F[S]$ the induced subgraph of $H[S]$ that is isomorphic to $F$.
%Moreover, Given a graph $H$ and a bipartite graph $F$ and disjoint vertex sets $S, T\subseteq V(H)$, denote by $F[S\dcup T]$ the bipartite induced subgraph of $H[S\dcup T]$ that is isomorphic to $F$.

	\subsection{Organization}
The rest of the paper is organized as follows. In the next section we construct  3-graphs which play a special role in the statements and proofs of our results.
 In Section 3 we introduce several lemmas and use them to deduce Theorems~\ref{P4}, \ref{e2}, and \ref{e3}.
The proofs of these lemmas are presented in Sections 4--6.  We prove Theorem~\ref{RN} in Section \ref{RamNum}. This proof  relies only on the statements of Theorems 1.1, 1.4, and 1.5, and thus can be understood without reading the earlier sections. Finally, the last section contains a couple of open problems.

	\section{Special $3$-graphs}\label{useful}
	
	In this section we define  $3$-graphs which play a special role in the paper, either as tools in the proofs or as extremal 3-graphs. By default, we drop the superscript $^{(3)}$.

The (unique) 6-vertex minimal 4-cycle $\c$ is a 3-graph with
	\[
	V(\c)=\{x_1,x_2,y_1,y_2,z_1,z_2\}\quad\mbox{ and }\quad E(\c)=\{x_1y_1y_2, y_1y_2x_2, x_2z_1z_2, z_1z_2x_1\}
	\]
	(see Figure \ref{fig:c4}).
Further, let $K\coloneq K_4$ stand for the complete 3-graph on four vertices and let $P\coloneq P_2$ denote the minimal 2-path with five vertices, that is, two edges sharing one vertex.

	%Recall that by $K^{(k)}_n$ we denote a complete $k$-graph, namely a $k$-graph with $n$ vertices consisting of all possible ${n\choose k}$ edges. For given $\ell$, by $P^{(2)}_\ell$ we denote a (graph) path with $\ell$ edges and by $C^{(2)}_\ell$ we denote a (graph) cycle with $\ell$ edges.
	For $s\ge2$, let $M_s$ stand for the matching of size $s$, that is, a $3$-graph consisting of $s$ disjoint edges. %For simplicity, we write $M_s\coloneq M^{(3)}_s$.
	
	\subsection{Stars}
A \emph{star} is a $3$-graph $S$ with a vertex $v$ (called sometimes the center) contained in all the edges of $S$. A star is \emph{full} if it consists of all sets in $\binom V3$ containing~$v$, that is, if $\deg_S(v)=\binom{|V|-1}{2}$. Normally, we write $S_n$ for the full star with $n$ vertices, but if we want to specify the vertex set and the star center, we may sporadically use symbol $S_V^v$ instead.
	By $S_n^{+1}$ we denote the unique (up to isomorphism) $n$-vertex 3-graph obtained from the full star $S_n$ by adding one extra edge. We call  $S_n^{+1}$ a \emph{starplus}.
	
\begin{figure}
     \centering
     \begin{subfigure}[b]{0.3\textwidth}
         \centering
       	\begin{tikzpicture}[scale = .9]
       	
       	\coordinate (x2) at (0,-1);
       	\coordinate (x1) at (0,1);
       	
       	\coordinate (y1) at (-2,0);
       	\coordinate (y2) at (-1,0);
       	
       	\coordinate (z1) at (1,0);
       	\coordinate (z2) at (2,0);

       	\node at (0.3,-1.4) {$x_2$};
       	\node at (0.3,1.3) {$x_1$};
       	
       	\node at (-2.3,.35) {$y_1$};
       	\node at (-0.8,.35) {$y_2$};
       	
       	\node at (0.8,.35) {$z_1$};
       	\node at (2.3,.35) {$z_2$};

       	\qedge{(z2)}{(z1)}{(x1)}{5pt}{1.5pt}{red!70!black}{red!50!white,opacity=0.2};
       	\qedge{(z1)}{(z2)}{(x2)}{5pt}{1.5pt}{red!70!black}{red!50!white,opacity=0.2};
       	\qedge{(x2)}{(y1)}{(y2)}{5pt}{1.5pt}{red!70!black}{red!50!white,opacity=0.2};
       	\qedge{(x1)}{(y2)}{(y1)}{5pt}{1.5pt}{red!70!black}{red!50!white,opacity=0.2};
       	
       	\foreach \i in {x1,x2,y1,y2,z1,z2}
       	\fill  (\i) circle (2pt);

       	\end{tikzpicture}
         \caption{$\c$}\label{fig:c4}
     \end{subfigure}
     \hfill
     \begin{subfigure}[b]{0.3\textwidth}
         \centering
       	\begin{tikzpicture}[scale = .8]

       	    \foreach \i in {1,...,5} {\coordinate (v\i) at (\i,0);}
       	    \coordinate (v0) at (3,2.5);
       	    \coordinate (v6) at (4.2,1.1);
       	    \coordinate (v7) at (5.1, 1.5);
       	    \coordinate (v8) at (.9, 1.5);
       	    
       	        \node [above right = .1] at (v0) {$v$};
       	      \node [left = .3] at (v1) {\Large $P_2$};

       	    \qedge{(v3)}{(v2)}{(v1)}{9pt}{1.5pt}{blue!70!black}{blue!50!white,opacity=0.2};
       	    \qedge{(v5)}{(v4)}{(v3)}{9pt}{1.5pt}{blue!70!black}{blue!50!white,opacity=0.2};
       	    \qqedge{(v0)}{(v6)}{(v5)}{6pt}{1.5pt}{red!70!black}{red!50!white,opacity=0.2}{8pt};
       	    \qqedge{(v0)}{(v7)}{(v5)}{6pt}{1.5pt}{red!70!black}{red!50!white,opacity=0.2}{10pt};
       	    \qqedge{(v1)}{(v8)}{(v0)}{6pt}{1.5pt}{red!70!black}{red!50!white,opacity=0.2}{10pt};
       	    \qedge{(v0)}{(v2)}{(v1)}{6pt}{1.5pt}{red!70!black}{red!50!white,opacity=0.2};
       	    \qedge{(v0)}{(v4)}{(v3)}{6pt}{1.5pt}{red!70!black}{red!50!white,opacity=0.2};
       	    
       	    \foreach \i in {0,...,8} {\fill (v\i) circle (2.5pt);}
       	    
       	\end{tikzpicture}
         \caption{$\sp_n$}\label{fig:sp}
     \end{subfigure}
     \hfill
     \begin{subfigure}[b]{0.3\textwidth}
         \centering
        	\begin{tikzpicture}[scale = .8]

        \coordinate (v1) at (-2,0);
        \coordinate (v2) at (0,-.5);
        \coordinate (v3) at (2,0);
        \coordinate (v4) at (0,.5);
        \coordinate (v0) at (0,2);
        \coordinate (v5) at (2.7,1.7);
        \coordinate (v6) at (-2.7,1.7);
        \coordinate (v7) at (2, 1.3);
        \coordinate (v8) at (-2, 1.3);
        
        \node [above right = .1] at (v0) {$v$};
        \node [left = .25] at (v1) {\Large $K_4$};

        \qedge{(v3)}{(v2)}{(v1)}{9pt}{1.5pt}{blue!70!black}{blue!50!white,opacity=0.2};
        \qedge{(v4)}{(v3)}{(v2)}{9pt}{1.5pt}{blue!70!black}{blue!50!white,opacity=0.2};
        \qedge{(v1)}{(v4)}{(v3)}{9pt}{1.5pt}{blue!70!black}{blue!50!white,opacity=0.2};
        \qedge{(v2)}{(v1)}{(v4)}{9pt}{1.5pt}{blue!70!black}{blue!50!white,opacity=0.2};

        \qedge{(v0)}{(v4)}{(v1)}{7pt}{1.5pt}{red!70!black}{red!50!white,opacity=0.2};
        \qedge{(v0)}{(v2)}{(v4)}{6pt}{1.5pt}{red!70!black}{red!50!white,opacity=0.2};
        \qqedge{(v0)}{(v3)}{(v2)}{6pt}{1.5pt}{red!70!black}{red!50!white,opacity=0.2}{20pt};
        %	\qqedge{(v1)}{(v6)}{(v0)}{6pt}{1.5pt}{red!70!black}{red!50!white,opacity=0.2}{20pt};
        \qqedge{(v1)}{(v8)}{(v0)}{6pt}{1.5pt}{red!70!black}{red!50!white,opacity=0.2}{10pt};
        %	\qqedge{(v0)}{(v5)}{(v3)}{6pt}{1.5pt}{red!70!black}{red!50!white,opacity=0.2}{20pt};
        \qqedge{(v0)}{(v7)}{(v3)}{6pt}{1.5pt}{red!70!black}{red!50!white,opacity=0.2}{10pt};
        
        \foreach \i in {0,...,4,7,8} {\fill (v\i) circle (2.5pt);}

        	\end{tikzpicture}\\
         \caption{$\sk_n$}\label{fig:sk}
     \end{subfigure}
        \caption{4-cycle $\c$, $\sp_n$ and $\sk_n$.}
        \label{fig:cpks}
\end{figure}
	
\subsection{$F$-stars} For  a set $V$ of $n\ge 6$ vertices, a subset $A\subset V$, and a vertex $v\in V\setminus A$, let $S(v,A)=S_V^v\setminus S_{V\setminus A}^v$ be the star obtained from the full star $S_V^v$ by deleting all edges disjoint from~$A$. In other words, $S(v,A)$ consists of all triples containing $v$ and at least one vertex of~$A$.

	Given a 3-graph $F$, we define the $F$-star by $SF_n\coloneq F\cup S(v,V(F))$, where $V\supset V(F)$, $|V|=n$, and $v\in V\setminus V(F)$.
	We will focus on two instances of $F$-stars: with $F=P$ and $F=K$ (see Figure \ref{fig:cpks}(b-c)).
	It is easy to check that both, $\sk_n$ and $\sp_n$, are $\{\p,\m\}$-free and contain a copy of $\c$. Moreover, $|\sk_n|=4n-10$ and $|\sp_n|=5n-18$. Notice that for $n=8$ these two expressions are equal to each other.

	\subsection{Balloons} 
\begin{figure}[h!]
	\centering
	\begin{subfigure}[b]{0.45\textwidth}
		\centering
		\begin{tikzpicture}[scale = .9]

	\coordinate (x) at (0,0);
	\coordinate (y1) at (1.5,-1.5);
	\coordinate (y2) at (1.5,0);
	\coordinate (y3) at (1.5,1.5);
	\coordinate (z1) at (3.5, -1.5);
	\coordinate (z2) at (3.5,0);
	\coordinate (z3) at (3.5,1.5);
	
	\foreach \i in {1,...,5}{
		\coordinate (v\i) at (105+\i*25:2cm);
	}

	\qedge{(y2)}{(y1)}{(y3)}{7pt}{1.5pt}{red!70!black}{red!50!white,opacity=0.2};
	\qedge{(z2)}{(z1)}{(z3)}{7pt}{1.5pt}{red!70!black}{red!50!white,opacity=0.2};
	
	\foreach \i/\j/\k in {
		y3/x/v1,
		y3/x/v2,
		y2/x/v3,
		y2/x/v1,
		v4/x/y2,
		v3/x/y1,
		v4/x/y1,
		v5/x/y1}
	\qedge{(\i)}{(\j)}{(\k)}{8pt}{1.5pt}{red!70!black}{red!50!white,opacity=0.2};%{30pt};
	
	\foreach \i/\j/\k in {
		v2/v1/x,
		v3/v2/x, 
		v4/v3/x, 
		x/v5/v4}
	\qedge{(\i)}{(\j)}{(\k)}{6pt}{1.5pt}{red!70!black}{red!50!white,opacity=0.2};
	
	\foreach \i in {1,2,3}
	\foreach \j in {1,2,3}{
		\path[green!80!black, line width = 2pt]
		(z\i) edge (y\j)
		;}
	
	\foreach \i in {x, y1, y2,y3, z1, z2,z3} {
		\fill (\i) circle (2pt);
	}
	
	\foreach \i in {1,...,5}{
		\fill (v\i) circle (2pt);
	}

	\node at ($(z3)+(.5,0)$) {$z_1$};			
	\node  at ($(x)+(.2,.5)$) {$x$};
	\node at ($(y1)+(.3,-.45)$) {$y_3$};
	\node  at ($(y2)+(0,-.5)$) {$y_2$};
	\node  at ($(y3)+(.3,.45)$) {$y_1$};
	\node at ($(z1)+(.5,0)$) {$z_3$};
	\node  at ($(z2)+(.5,0)$) {$z_2$};
	
	\end{tikzpicture}
					
				\caption{Balloon $B_n$. }\label{fig:b}
			\end{subfigure}
			\hfill    
			\begin{subfigure}[b]{0.45\textwidth}
				\centering
				\begin{tikzpicture}[scale = .9]
			
			\phantom {\fill (0,-2.2) circle (1pt);}
			
					\coordinate (x) at (0,0);
					\coordinate (y1) at (1.5,1);
					\coordinate (y2) at (1.5,-1);
					\coordinate (z1) at (3.5, 1);
					\coordinate (z2) at (3.5,-1);
					
					\foreach \i in {1,...,5}{
						\coordinate (v\i) at (105+\i*25:2cm);
					}
					
					\qedge{(y2)}{(y1)}{(z1)}{6pt}{1.5pt}{red!70!black}{red!50!white,opacity=0.2};
					\qedge{(y2)}{(y1)}{(z2)}{6pt}{1.5pt}{red!70!black}{red!50!white,opacity=0.2};
					
					\foreach \i/\j/\k in {
						y1/x/v1,
						y1/x/v2,
						v3/x/y2,
						v4/x/y2,
						v5/x/y2}
					\qedge{(\i)}{(\j)}{(\k)}{8pt}{1.5pt}{red!70!black}{red!50!white,opacity=0.2};%{30pt};
					
					\foreach \i/\j/\k in {
						v2/v1/x,
						v3/v2/x, 
						v4/v3/x, 
						x/v5/v4}
					\qedge{(\i)}{(\j)}{(\k)}{6pt}{1.5pt}{red!70!black}{red!50!white,opacity=0.2};
					
					\path[green!80!black, line width = 2pt]
					(z1) edge (z2)
					(z1) edge (y1)
					(z1) edge (y2)
					(z2) edge (y1)
					(z2) edge (y2)
					(y1) edge (y2)
					;
					
					\foreach \i in {x, y1, y2, z1, z2} {
						\fill (\i) circle (2pt);
					}
					
					\foreach \i in {1,...,5}{
						\fill (v\i) circle (2pt);
					}
					
					\node  at ($(x)+(.2,.45)$) {$x$};
					\node at ($(y1)+(.1,.5)$) {$y_1$};
					\node  at ($(y2)+(.1,-.5)$) {$y_2$};
					\node at ($(z1)+(.1,.45)$) {$z_1$};
					\node  at ($(z2)+(.1,-.45)$) {$z_2$};

		\end{tikzpicture}
	\caption{Compact balloon $CB_n$. }\label{fig:cb}
	\end{subfigure}  
	\hfill
	\hfill  
	\caption{Balloons. The green pairs form 3-edges with the vertex $x$.}	\label{fig:bcb}
	\vspace{-1em}
\end{figure}
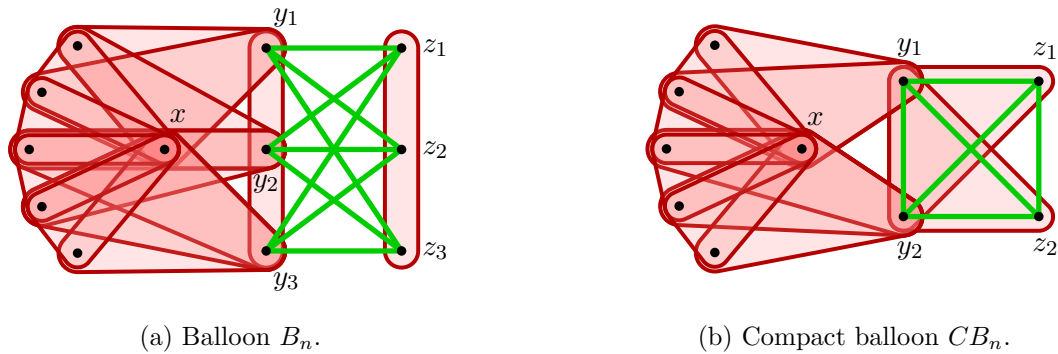

Finally, we define two more deformations of stars. For $n\ge9$, let $B_n$ be a 3-graph on $n$ vertices, called the \emph{balloon}, obtained from the full star $S_{n-3}$ with center $x$ by selecting three vertices $y_1,y_2,y_3\in V(S_{n-3})\setminus\{x\}$, adding three new vertices $z_1,z_2,z_3$, and adding eleven new edges: $\{y_1,y_2,y_3\}$, $\{z_1,z_2,z_3\}$, and all nine edges of the form $\{x,y_i,z_j\}$, $i,j=1,2,3$ (see Figure \ref{fig:b}).
	Note that the balloon $B_n$ is $\p$-free, contains $M_3$, and has $\binom{n-4}2+11$ edges.
	
	For $n\ge 8$, let $CB_n$ be a 3-graph on $n$ vertices, called the \emph{compact balloon}, obtained from the full star $S_{n-2}$ with center $x$ by selecting two vertices $y_1,y_2\in V(S_{n-2})\setminus\{x\}$, adding two new vertices $z_1,z_2$, and adding seven new edges: $\{y_1,y_2,z_1\}$, $\{y_1,y_2,z_2\}$,  all four edges of the form $\{x,y_i,z_j\}$, $i,j=1,2$, and the edge $\{x,z_1,z_2\}$ (see Figure \ref{fig:cb}).
	Note that the compact balloon $CB_n$ is $\p$-free, is not a sub-3-graph of the starplus $S_n^{+1}$, and  has $\binom{n-3}2+7$ edges.

%%%%%%%%%%%%%%%%%%%%%%%%%%%%%%%%%%%%%%%%%%%%%%%%%%%%%%%%%%%%%%%%%%%%%%
%																																		%
%																																		%
%											Turan numbers																%
%																																		%
%																																		%
%%%%%%%%%%%%%%%%%%%%%%%%%%%%%%%%%%%%%%%%%%%%%%%%%%%%%%%%%%%%%%%%%%%%%%

\section{Tur\'an numbers}

\begin{figure}[h!]
	\centering
	\begin{tikzpicture}
	\coordinate (a) at (0,2.8);
	\coordinate (b) at (-4,1.8);
	\coordinate (c) at (4,1.8);
	\coordinate (x) at (-6,.3);
	\coordinate (y) at (-2,.4);
	\coordinate (z) at (2.5,0);
	\coordinate (v) at (6,0);
	\coordinate (p) at (-4,-1);
	\coordinate (q) at (.3,-2.5);
	\coordinate (r) at (-2.8,-4.2);
	\coordinate (s) at (0,-4.2);
	\coordinate (t) at (2.8,-4.2);
	\coordinate (tt) at (6.5,-4.2);
	\coordinate (u) at (4,-2);
	\coordinate (w) at (8,-2);
	\coordinate (p1) at ($(p)-(1.5,1.5)$);
	\coordinate (p2) at ($(p)+(1.5,-1.5)$);

	\node at ($(b)!.3!(x)$) [anchor=east] {\tiny Lemma \ref{l:M3}};
	\node at ($(p)!.4!(p1)$) [anchor=east] {\tiny Lemma \ref{l:c41}};
	\node at ($(p)!.4!(p2)$) [anchor=west] {\tiny Lemma \ref{l:c42}};
	\node at ($($(q)!.33!(tt)$)+(.3,0)$) [anchor=west] {\tiny Lemma \ref{l:nu2}};
	\node at ($(c)!.26!(z)$) [anchor=east] {\tiny Lemma \ref{l:disc}};
	\node at ($(c)!.4!(v)$) [anchor=west] {\tiny by induction};
	
	\draw[ line width=1.5pt, line cap=round, ->, shorten >= 27] (a) -- (b);
	\draw[ line width=1.5pt, line cap=round, ->, shorten >= 33] (a) -- (c);
	\draw[ line width=1.5pt, line cap=round, ->, shorten >= 33] (b) -- (x);
	\draw[ line width=1.5pt, line cap=round, ->, shorten >= 21] (b) -- (y);
	\draw[ line width=1.5pt, line cap=round, ->, shorten >= 36] (c) -- (z);
	\draw[ line width=1.5pt, line cap=round, ->, shorten >= 15] (c) -- (v);
	\draw[ line width=1.5pt, line cap=round, ->, shorten >= 22] (y) -- (p);
	\draw[ line width=1.5pt, line cap=round, ->, shorten >= 15] (y) -- (q);
	\draw[ line width=1.5pt, line cap=round, ->, shorten >= 30] (q) -- (t);
	\draw[ line width=1.5pt, line cap=round, ->, shorten >= 60] (q) -- (tt);
	\draw[ line width=1.5pt, line cap=round, ->, shorten >= 37] (q) -- (r);
	\draw[ line width=1.5pt, line cap=round, ->, shorten >= 18] (q) -- (s);
	\draw[ line width=1.5pt, line cap=round, ->, shorten >= 26] (v) -- (u);
	\draw[ line width=1.5pt, line cap=round, ->, shorten >= 42] (v) -- (w);
	\draw[ line width=1.5pt, line cap=round, ->, shorten >= 25] (p) -- (p1);
	\draw[ line width=1.5pt, line cap=round, ->, shorten >= 25] (p) -- (p2);
	
	\node at (a) [draw, fill=white]  {$\p\nsubseteq \cH$};
	\node at (b) [draw, fill=white]  {connected};
	\node at (c) [draw,fill=white, text width = 2cm, text centered]  {disconnected};
	\node at (x) [draw,fill=white!90!blue, text width = 2cm, text centered]
			{$\m\subseteq \cH$ \tiny{$|\cH|\le {n-4\choose 2}+11$}};
	\node at (y) [draw,fill=white, text width = 2cm, text centered]
			{$\m\nsubseteq \cH$ };
	\node at (z) [draw,fill=white!90!black,text width=3.5cm, text centered]
			{$\delta_1(\cH)\ge 1$ \tiny{ ${|\cH|\le \begin{cases}
					{n-3\choose 3} + 1, n\in [6,9]\cr
					{n-6\choose 3} + 20, n \in [10,12] \cr
					{n-4\choose 2} +3, n \ge 13
					\end{cases}}$}};
	\node at (v) [draw,fill=white, text width = 2cm, text centered]
			{$\delta_1(\cH)=0$};
	\node at (p) [draw,fill=white, text width = 2cm, text centered]
			{$\c\subseteq \cH$};
	\node at (p1) [draw,fill=white, text width = 2cm, text centered]
			{$\cH\subseteq SP_n$ \tiny{$|\cH|\le 5n-18$} };
	\node at (p2) [draw,fill=white, text width = 2cm, text centered]
			{$\cH\nsubseteq SP_n$ \tiny{$|\cH|\le 4n-10$}};
	\node at (q) [draw,fill=white, text width = 2cm, text centered] {$\c\nsubseteq \cH$};
	\node at (r) [draw,fill=white!90!red, text width = 2cm, text centered]
			{$\cH\subseteq S_n^{+1}$ \tiny $|\cH|\le {n-1\choose 2} +1$};
	\node at (s) [draw,fill=white!90!green, text width = 2cm, text centered]
			{$\cH\subseteq CB_n$ \tiny $|\cH|\le {n-3\choose 2} +7$};
	\node at (t) [draw,fill=white!90!black, text width = 2cm, text centered]
			{$\cH\subseteq SP_n$ \tiny $|\cH|\le 5n-19$};
	\node at (tt) [draw,fill=white!90!black, text width = 3.8cm, text centered]
			{$\cH\nsubseteq S_n^{+1}, CB_n, SP_n$ \tiny $|\cH|\le \max\left\{4n-11,{n-4\choose 2}+10 \right\}$};
	\node at (u) [draw,fill=white!90!black, text width = 3.9cm, text centered]
				{$\cH\subseteq S^{+1}_n$ \tiny $|\cH|\le \ex(n-1;\p)={n-2\choose 2} +1$};
	\node at (w) [draw,fill=white!90!black, text width = 3cm, text centered]
				{$\cH\nsubseteq S^{+1}_n$ \tiny ${|\cH|\le \ex^{(2)}(n-1;\p)}$ $=\begin{cases} 5n-23,&n\le 12\cr {n-4\choose 2}+7,&n\ge 13\end{cases}$};
	\end{tikzpicture}
	\caption{Division of the family of $\p$-free 3-graphs. The gray blocks contain 3-graphs not appearing in extremal families  of the first three orders. For $n\ge 14$ the \crr{red}, \cgg{green}, and \cbb{blue} block represent, respectively, the \crr{first}, \cgg{second}, and  \cbb{third} order Tur\'an number for $\p$.}
\label{fig:split}
\end{figure}
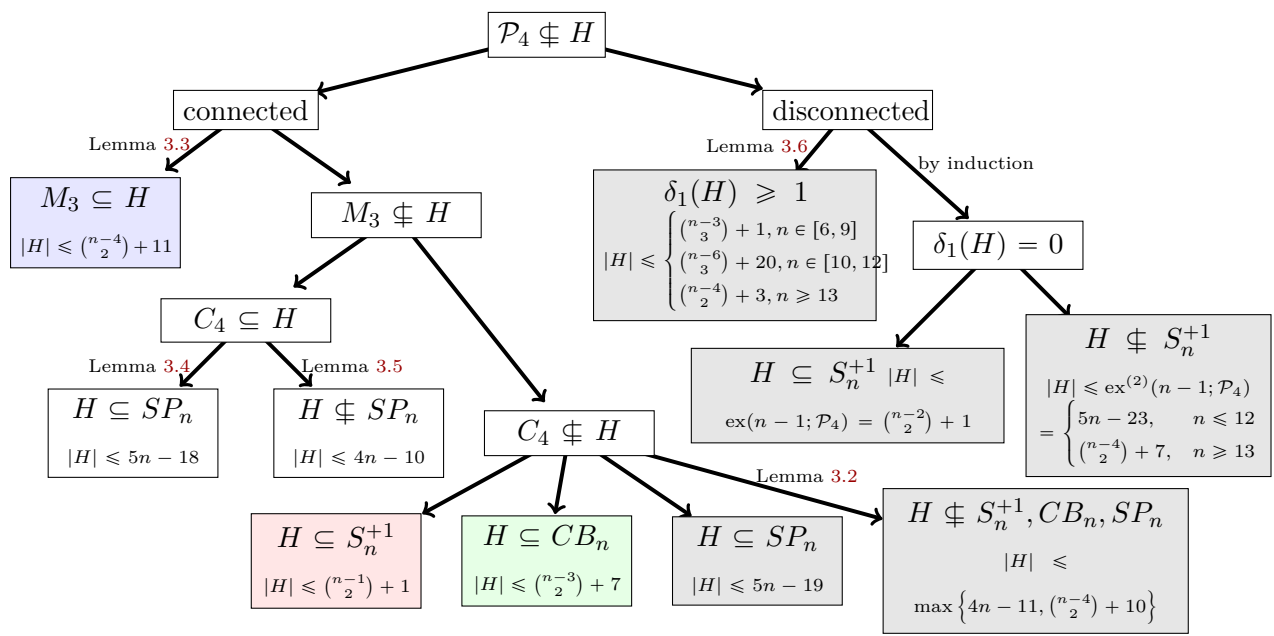

The goal of this section is to prove Theorems \ref{P4}, \ref{e2}, and \ref{e3}. In order to do this we divide the family of all $\p$-free 3-graphs into some special subfamilies and then count the maximum number of edges within them separately (see Figure \ref{fig:split}).

Next, we compare to each other bounds obtained in Lemmas \ref{l:nu2}-\ref{l:disc}. For $n\ge 14$ we have,
	\begin{equation}\label{eq:turan}
		4n-10< 5n-18 < {n-4\choose 2} + 11 < {n-3\choose 2} +7 < {n-1\choose 2}+1,
	\end{equation}
whereas for $n\in [8,14]$ we gather these bounds in Table \ref{tab}.

\begin{table}[h!]
	\begin{center}
		\begin{tabular}{ |c ||c| c |c|c|c|c|c|}
			\hline
			&$S^{+1}_n$ & $SP_n$ & $SK_n$ & $CB_n$ & $B_n$&$\max\{4n-11,$&disconn\\
			$n$& ${n-1\choose 2}+1$ & $5n-18$ & $4n-10$ & ${n-3\choose 2} +7$ & ${n-4\choose 2} + 11$&
			${n-4\choose 2}+10\}$&
			Lemma \ref{l:disc}\\
			\hline\hline
			8 & \crr{22} & \crr{22} & \crr{22}& 17 &17&21&11\\ \hline
			9 & \crr{29} & \cgg{27} & \cbb{26}& 22&21&25&21 \\\hline
			10 & \crr{37} & \cgg{32} & \cbb{30}& 28&26&29&24 \\\hline
			11 & \crr{46} & \cgg{37} & 34& \cbb{35}&32&33&30 \\\hline
			12 & \crr{56} & \cbb{42} & 38& \cgg{43}&39&38&40 \\\hline
			13 & \crr{67} & \cbb{47} & 42& \cgg{52}&\cbb{47}&46&39 \\\hline
			14 & \crr{79} & 52 & 46& \cgg{62}&\cbb{56}&55&48 	\\\hline
			
		\end{tabular}
		\caption{The Tur\'an numbers for $\p$ and $n\in [9,14]$ of the \crr{first}, \cgg{second} and \cbb{third} order.}	\label{tab}
	\end{center}
\end{table}

However, before we do this precisely, we need one more piece of notation. A $3$-graph $\cH$ is said to be \emph{connected} if for every partition of the vertex set $V(\cH)=U\dcup W$, there is an edge in $\cH$ with non-empty intersection with both subsets, $U$ and $V$.

A forced presence of a sub-$k$-graph can be expressed in terms of conditional Tur\'an numbers, introduced in \cite{JPR}.
For a  $k$-graph $ F$, an  $ F$-free $k$-graph $ G$, and an integer $n\ge | G|$, the  \textit{conditional Tur\'an number} is defined as
\begin{eqnarray*}
	\ex_k(n; F| G)=\max\{|E(\cH)|:|V(\cH)|=n,\;
	\mbox{$\cH$ is
		$ F$-free, and } \cH\supseteq G \}.
\end{eqnarray*}
Every $n$-vertex $F$-free $k$-graph $\cH$ with $\ex_k(n; F| G)$ edges and such
that $\cH\supseteq G$ is called \emph{$ G$-extremal for $ F$}.
We denote by $\Ex_k(n; F| G)$ the family of all  $n$-vertex $k$-graphs which are $ G$-extremal for $ F$. 
%If $\mathcal F=\{F\}$ or $\mathcal G=\{G\}$, we will simply write
%$\ex_k(n;F| \mathcal G)$, $\ex_k(n;\mathcal F|G)$, or $\ex_k(n;F|G)$,  respectively, similar for $\Ex_k(n;\cF| \cG)$.
For $k=3$ we drop the subscript $_3$.
The \emph{conditional Tur\'an number of the $s^{th}$ order} is defined in a similar way as the ordinary Tur\'an number of the $s^{th}$ order (see Subsection \ref{ss:turan}). Finally, for $k=3$, if in the above definition one restricts oneself to connected $3$-graphs, we add the subscript $_{conn}$ and denote the corresponding extremal numbers and families, respectively, by $\ex_{conn}(n; F|G)$, $\Ex_{conn}(n; F|G)$, $\ex^{(s)}_{conn}(n; F|G)$, and $\Ex^{(s)}_{conn}(n; F|G)$.

Now let us state a few lemmas from which Theorems \ref{P4}, \ref{e2}, and \ref{e3} follow. The case $n=7$ is treated separately.

\begin{lemma}\label{l:n7}
	$\ex(7;\p) = 20$, $\Ex(7;\p) = \{K^{(3)}_6\cup K_1\}$.
\end{lemma}

\begin{lemma}\label{l:nu2}
	Let $\cH$ be a $\{\p,\c, \m\}$-free connected 3-graph on $n\ge 8$ vertices. If $\cH\nsubseteq S_n^{+1}$, ${\cH\nsubseteq SP_n}$, and $\cH\nsubseteq CB_n$ then 	
		\[
			|\cH| \le \max\left\{4n-11,{n-4\choose 2}+10\right\}.
		\]
\end{lemma}

\begin{lemma}\label{l:M3}
	For $n\ge 9$, $\ex_{conn}(n; \p|\m) = \binom{n-4}{2}+11$ and the balloon $B_n$ is the only extremal $3$-graph.
\end{lemma}

\begin{lemma}\label{l:c41} For $n\ge 8$,
	\begin{align*}
	&\ex_{conn}(n;\p\cup\{\m\}|\c)=5n-18, \\& \Ex_{conn}(n;\p\cup\{\m\}|\c)=\left\{ \begin{array}{ll}
	\{\sp_8,\sk_8\} &\textrm{for } n=8,\\
	\{\sp_n\} & \textrm{for } n\ge 9.
	\end{array} \right.
	\end{align*}
\end{lemma}

\begin{lemma}\label{l:c42} For $n\ge 9$,
	\begin{align*}
	&\ex^{(2)}_{conn}(n;\p\cup\{\m\}|\c)=4n-10, \\& \Ex^{(2)}_{conn}(n;\p\cup\{\m\}|\c)=\{\sk_n\}.
	\end{align*}
\end{lemma}

\begin{lemma}\label{l:disc}
	If $\cH$ is a disconnected $\p$-free 3-graph on $n$ vertices, with $\delta_1(\cH)\ge 1$, then
		\[
			|\cH|\le \begin{cases}
				{n-3\choose 3} + 1, &\textrm{ for } 6\le n\le 9,\cr
				{n-6\choose 3} + 20, &\textrm{ for } 10 \le n \le 12, \cr
				{n-4\choose 2} +3, &\textrm{ for } n \ge 13.
			\end{cases}.
		\]
\end{lemma}
\begin{proof}
	Let $\cH_1$ be a connected component of $\cH$ with the smallest number of vertices. Set $\cH_2=\cH\setminus \cH_1$, $n_i=|V[\cH_i]|$, $i=1,2$. Clearly $3\le n_1\le n_2=n-n_1\le n-3$.
	
	We argue by induction on $n$. For the base case $6\le n\le 9$, we use the fact that $|\cH_i|\le |K^{(3)}_{n_i}|= {n_i\choose 3}$, $i=1,2$. Therefore, a simple optimization shows
	\[
	|\cH|= |\cH_1|+|\cH_2|\le {n_1\choose 3} + {n_2\choose 3} \le 1+{n-3\choose 3},
	\]
	as required.
	
	For the induction step assume $n\ge 10$ and that Lemma \ref{l:disc} is true for all disconnected $\p$-free 3-graphs with less than $n$ vertices and $\delta_1(\cH)\ge1$. Then, as $ n_2\le n-3$ we are in position to apply the induction hypothesis to $\cH_2$ in case it is disconnected.  For $\cH_1$, as well as, for connected $\cH_2$ we apply Lemmas \eqref{l:n7}-\eqref{l:c41}. Altogether, we claim that, for $i=1,2$,
	\begin{equation}\label{eq:H2}
		|\cH_i|\le \begin{cases}
			{n_i\choose 3}, & \textrm{ for }n_i\le 6,\cr
			19, & \textrm{ for }n_i=7,\cr
			{n_i-1\choose 2}+1, &\textrm{ for }n_i\ge 8.
		\end{cases}
	\end{equation}
	Indeed, for $n_i\le 6$ clearly $|\cH_i|\le |K^{(3)}_{n_i}|\le {n_i\choose 3}$, whereas for $n_i=7$ $|\cH_i|\le 19$ follows from Lemma \ref{l:n7} combined with $\delta_1(\cH)\ge 1$. Finally, to show that $|\cH_i|\le {n_i-1\choose 2}+1$ for $n_i\ge 8$, in view of Lemmas \ref{l:nu2}, \ref{l:M3}, \ref{l:c41}, and \ref{l:disc} (see also Figure \ref{fig:split}) it is enough to observe that
		\begin{equation}\label{eq:win}
				{n-1\choose 2}+1 > \begin{cases}
				{n-3\choose 3} + 1, &\textrm{ for } n\le 10,\cr
				{n-6\choose 3} + 20, &\textrm{ for } 8 \le n \le 14, \cr
				\max\left\{{n-3\choose 2} + 7, {n-4\choose 2} + 11, 5n-19\right\}, &\textrm{ for } n \ge 8.
				\end{cases}
		\end{equation}
	 In particular, for $n\ge8$, ${n-1\choose 2}+1 >4n-11$, as well as, ${n-1\choose 2}+1 \ge5n-18$.
		
	Now we use \eqref{eq:H2} to bound the number of edges in $\cH$. Considering separately cases $n_1=3,4,\dots, \lfloor n/2 \rfloor$ one gets,
		\[
			|\cH|\le \begin{cases}
				\max\{1+19, 4 + 20, 10 + 10\} = {n-6\choose 3} + 20, &\textrm{ for } n=10,\cr
				\max\{1 + 22, 4+19,10+20\} = {n-6\choose 3} + 20, &\textrm{ for } n=11,\cr
				\max\{1+29, 4 + 22,10+19, 20+20  \} ={n-6\choose 3} + 20, &\textrm{ for } n=12,\cr
				\max\{1+37, 4+29,10+22,20+19\}={n-4\choose 2} + 3, &\textrm{ for } n=13.
			\end{cases}
		\]

	Therefore it remains to take care of $n\ge 14$. If $n_1\le 6$, then $n_2\ge 8$ and thus $|\cH_1|\le {n_1\choose 3}$, $|\cH_2|\le {n_2-1\choose 2}+1$, yielding
		\[
			|\cH| \le
			\max\left\{{n-4\choose 2} + 2,{n-5\choose 2}+5, {n-6\choose 2} +11, {n-7\choose 2} + 21\right\}={n-4\choose 2} + 2.
		\]
	For $n_1 = 7$, $|\cH_1|\le 19$ and hence
		\[
			|\cH|\le \begin{cases}
				19+19< {n-4\choose 2} + 3, &\textrm{ for } n=14,\cr
				{n-8\choose 2} +20 < {n-4\choose 2} + 3, &\textrm{ for } n\ge 15.
			\end{cases}
		\]
		Finally, if $n_1\ge 8$, then also $n_2\ge 8$ and thus $|\cH_i|\le {n_i\choose 2}+1$, $i=1,2$. But then, clearly
		\[
			|\cH|= |\cH_1|+|\cH_2|\le {n_1-1\choose 2} + {n_2-1\choose 2} +2< {n-4\choose 2} +3. \qedhere
		\]
\end{proof}

	Now we are ready to prove Theorems \ref{P4}, \ref{e2}, and \ref{e3}.

\begin{proof}[Proof of Theorem \ref{P4}]
	We argue by induction on $n$. For the base case $n\le 6$ the assumption easily follows from the fact that every minimal 4-path has at least 7 vertices, whereas for $n= 7$ we use Lemma \ref{l:n7}.
	
	Next, we let $n\ge 8$ and observe that as $S_n^{+1}$ is a $\p$-free 3-graph with ${n-1\choose 2}+1$ edges, we get
		\[
			\ex(n;\p)\ge {n-1\choose 2} + 1.
		\]

	To obtain the reverse bound on $\ex(n;\p)$ we let $\cH$ to be a $\p$-free 3-graph on $n\ge 8$ vertices and with at least ${n-1\choose 2}+1$ edges. We argue that $\cH=S^{+1}_n$ for $n\ge 9$, whereas for $n=8$, $\cH = S_n^{+1}$ or $\cH\in \{\sp_8,\sk_8\}$, which will end the proof.
	To this end we consider separately connected and disconnected $\p$-free 3-graphs.
	In the former case Lemma \ref{l:M3} together with ${n-4\choose 2} + 11 < {n-1\choose 2} +1$ tells us that $\m\nsubseteq \cH$.
	Further, as for $n\ge 8$ we have $5n-18\le {n-1\choose 2} +  1$ with the equality only for $n=8$, in view of Lemma \ref{l:c41} we learn that for $n\ge 9$, $\cH$ is $\c$-free, whereas for $n=8$ the only possibility to have $\c\subseteq\cH$ is $\cH\in \{\sp_8,\sk_8\}$.
	Finally we use Lemma~\ref{l:nu2} to deduce that the only $\{\p,\c,\m\}$-free 3-graph with at least ${n-1\choose 2}+1$ edges is $S_n^{+1}$, as required (see Figure \ref{fig:split}, Table \ref{tab} and \eqref{eq:turan}).

	Now, to exclude the disconnected case we first assume that $\delta_1(\cH)\ge 1$ and use Lemma \ref{l:disc} combined with~\eqref{eq:win}. Finally, if $\cH$ contains an isolated vertex $v$, then we can apply the induction hypothesis to $\cH-v$, obtaining
		\[
			|\cH| = |\cH-v|\le \ex(n-1; \p) < {n-1\choose 2}+1,
		\]
	which ends the proof.
\end{proof}

\begin{proof}[Proof of Theorem \ref{e2}]
	The proof is similar to the proof of Theorem \ref{P4}. Let $\cH$ be a $\p$-free 3-graph on the set of vertices $V$, $|V|=n\ge 9$ with $|\cH|=\ex^{(2)}(n;\p)$. Moreover, as we are computing the second order Tur\'an number and $\Ex(n;\p)=\{S_n^{+1}\}$ for $n\ge 9$, we may assume that $\cH\nsubseteq S_n^{+1}$. Because both 3-graphs $SP_n$ and $CB_n$ are $\p$-free and are not contained in $S_n^{+1}$, we have the lower bound
		\begin{equation}\label{eq:turan2}
			|\cH|=\ex^{(2)}(n;\p)\ge\max\left\{5n-18, {n-3\choose 2}+7\right\}= \begin{cases}
				5n-18, &\textrm{ for } n\le 11,\cr
				{n-3\choose 2}+7, &\textrm{ for } n\ge 12.
			\end{cases}
		\end{equation}
	 We argue that $\cH=SP_n$ for $n\le 11$ and $\cH=CB_n$ for $n\ge 12$. The proof is by induction on $n$.
	
	 First assume that $\cH$ is connected and notice that since ${n-4\choose 2}+11 < {n-3\choose 2}+7$ for $n\ge 9$, Lemma~\ref{l:M3} yields $\m\nsubseteq \cH$. Therefore, since $4n-11 < 5n-18$, in view of Lemmas \ref{l:nu2} and \ref{l:c41} combined with $\cH\nsubseteq S^{+1}_n$, either $\cH=SP_n$ or $\cH=CB_n$, as required (see Figure \ref{fig:split}, Table \ref{tab} and \eqref{eq:turan}).
	
	 In the disconnected case Lemma \ref{l:disc} tells us that $\delta_1(\cH) = 0$, because clearly ${n-4\choose 2}+3 < {n-3\choose 2} + 7$ and for $n\le 12$ the bound obtained in this lemma is smaller than $5n-18$ (see Table \ref{tab}). Thus we let $v$ be an isolated vertex of $\cH$. For the base case, $n=9$ we use Theorem \ref{P4}, getting
		 \[
			 |\cH| = |\cH-v|\le \ex(8,\p) = 22 < 27 = 5n-18.
		 \]
	For the induction step assume $n\ge 10$ and that Theorem \ref{e2} is true for $n-1$ in place of $n$. Now observe, that because $\cH\nsubseteq S_n^{+1}$, we also have $\cH-v\nsubseteq S_{n-1}^{+1}$, and consequently,
		\[
			|\cH|=|\cH-v|\le \ex^{(2)}(n-1,\p) = \begin{cases}
				5n-23, &\textrm{ for } n\le 12,\cr
				{n-4\choose 2} + 7, &\textrm{ for } n\ge 13,
			\end{cases}
		\]
	contradicting \eqref{eq:turan2}.
\end{proof}

The proof of Theorem \ref{e3} is very similar to the one of Theorem \ref{e2}, and therefore we left it to the Reader (see Figure \ref{fig:split}, Table \ref{tab}, and \eqref{eq:turan}).

%%%%%%%%%%%%%%%%%%%%%%%%%%%%%%%%%%%%%%%%%%%%%%%%%%%%%%%%%%%%%%%%%%%%%%
%																																		%
%																																		%
%								 Seven vertices - the proof of Lemma 4.1											%
%																																		%
%																																		%
%%%%%%%%%%%%%%%%%%%%%%%%%%%%%%%%%%%%%%%%%%%%%%%%%%%%%%%%%%%%%%%%%%%%%%

\section{Seven vertices - proof of Lemma \ref{l:n7}}

%%%%%%%%%%%%%%%%%%%%%%%%%%%%%%%%%%%%%%%%%%%%%%%%%%%%%%%%%%%%%%%%%%%%%%
%																																		%
%																																		%
%											Two-colored graphs																%
%																																		%
%																																		%
%%%%%%%%%%%%%%%%%%%%%%%%%%%%%%%%%%%%%%%%%%%%%%%%%%%%%%%%%%%%%%%%%%%%%%

\subsection{2-colored graphs without a forbidden pattern.}\label{2color}

\def\rb{$\textcolor{red!60!black}{rr}$-$\textcolor{blue!60!black}{bb}$-path\ }
\def\RB{P_4^{\textcolor{red!60!black}{R}\textcolor{blue!60!black}{B}}}

In the whole subsection we consider only ordinary 2-graphs, therefore for simplicity of notation we omit the superscript $^{(2)}$ here. We prove two lemmas needed in the proof of Lemma~\ref{l:n7}, where link graphs, $\rr$ and $\bb$, of two given vertices are considered. However, before we state them, one more piece of notation is needed. Let two graphs, $\rr$ and $\bb$, on the same vertex set be given. We define an
	\emph{\rb $\RB =\tikz{\draw [blue!80!black](0,-.12) -- (.2,-.06)--(0,0);\draw [red!80!black] (0,0) -- (.2,.06)--(0,.12); \foreach \i/\j in {0/-.12,.2/-.06,0/0,.2/.06,0/.12}  \fill (\i,\j) circle (.6pt);}\ $}
to be a subgraph of $\rr\cup\bb$ consisting of 4 edges, $\crr{r_1,r_2}\in \rr$ and $\cbb{b_1,b_2}\in \bb$, such that $\crr{r_1r_2}\cbb{b_1b_2}$ is the 4-edge path $P_4$.
By $T\cup \{e\}$ we denote a graph on five vertices consisting of a complete graph on three vertices $T=K_3$ and a single edge $e$, disjoint from $V[T]$. We start with two technical facts used in further proofs.

\begin{figure}[h!]
	\centering			
	\begin{multicols}{7}
		
		\begin{tikzpicture}		
		\node at (-.8,.7){{\tiny $(a)$}};
			\foreach \i in {0,...,4}{
				\coordinate (v\i) at (90+\i*72:.7cm);
			}				
		
		\qedge{(v1)}{(v1)}{(v2)}{4pt}{.4pt}{black}{black!30!white,opacity=0.2};
		\qedge{(v3)}{(v0)}{(v4)}{4pt}{.4pt}{black}{black!30!white,opacity=0.2};
		
		\foreach \i/\j in {0/1,1/4,1/3}{
			\draw (v\i) edge [color = red!90!black, thick] (v\j);
		}			
		\draw ($(v2)+(0,0.023)$) edge [color = blue!90!black,thick] ($(v3)+(0,.023)$);
		\draw ($(v2)-(0,0.023)$) edge [color = red!90!black, thick] ($(v3)-(0,.023)$);
		
		\draw ($(v2)+(0,0.03)$) edge [color = blue!90!black,thick, dashed] ($(v4)+(0,.03)$);
		\draw ($(v2)-(0,0.03)$) edge [color = red!90!black, thick] ($(v4)-(0,.03)$);
		
		\draw ($(v2)-(0.025,0)$) edge [color = blue!90!black,thick, dashed] ($(v0)-(.025,0)$);
		\draw ($(v2)+(0.025,0)$) edge [color = red!90!black, thick] ($(v0)+(.025,0)$);
		
		\foreach \i/\j in {0/3,3/4}{
			\draw (v\i) edge [color = blue!90!black, thick, dashed] (v\j);
		}			
		\foreach \i in {0,...,4}{
			\fill (v\i) circle (1.5pt);
		}		
		\end{tikzpicture}\\

		\begin{tikzpicture}		
		\node at (-.8,.7){{\tiny $(b)$}};
	%	\foreach \i in {0,...,4}{
	%		\coordinate (v\i) at (90+\i*72:.7cm);
	%	}		
		\qedge{(v1)}{(v1)}{(v2)}{4pt}{.4pt}{black}{black!30!white,opacity=0.2};
		\qedge{(v3)}{(v0)}{(v4)}{4pt}{.4pt}{black}{black!30!white,opacity=0.2};

		\draw ($(v2)+(0,0.023)$) edge [color = blue!90!black,thick] ($(v3)+(0,.023)$);
		\draw ($(v2)-(0,0.023)$) edge [color = red!90!black, thick] ($(v3)-(0,.023)$);
		
		\draw ($(v1)+(0,0.02)$) edge [color = blue!90!black,thick, dashed] ($(v4)+(0,.02)$);
		\draw ($(v1)-(0,0.02)$) edge [color = red!90!black, thick] ($(v4)-(0,.02)$);
		
		\draw ($(v2)+(0,0.03)$) edge [color = blue!90!black,thick, dashed] ($(v4)+(0,.03)$);
		\draw ($(v2)-(0,0.03)$) edge [color = red!90!black, thick] ($(v4)-(0,.03)$);
		
		\draw ($(v1)+(0,0.03)$) edge [color = blue!90!black,thick, dashed] ($(v3)+(0,.03)$);
		\draw ($(v1)-(0,0.03)$) edge [color = red!90!black, thick] ($(v3)-(0,.03)$);
		
		\draw ($(v0)+(0,0.03)$) edge [color = blue!90!black,thick] ($(v1)+(0,.03)$);
		\draw ($(v0)-(0,0.03)$) edge [color = red!90!black, thick] ($(v1)-(0,.03)$);
		
		\draw ($(v2)-(0.025,0)$) edge [color = blue!90!black,thick, dashed] ($(v0)-(.025,0)$);
		\draw ($(v2)+(0.025,0)$) edge [color = red!90!black, thick] ($(v0)+(.025,0)$);
		
		\foreach \i/\j in {0/3,3/4,0/4}{
			\draw (v\i) edge [color = blue!90!black, thick, dashed] (v\j);
		}							
		\foreach \i in {0,...,4}{
			\fill (v\i) circle (1.5pt);
		}		
		\end{tikzpicture}\\

			\begin{tikzpicture}		
			\node at (-.8,.7){{\tiny $(c)$}};
			%	\foreach \i in {0,...,4}{
			%		\coordinate (v\i) at (90+\i*72:.7cm);
			%	}		
			\qedge{(v1)}{(v1)}{(v2)}{4pt}{.4pt}{black}{black!30!white,opacity=0.2};
			\qedge{(v3)}{(v0)}{(v4)}{4pt}{.4pt}{black}{black!30!white,opacity=0.2};

			\draw ($(v2)+(0,0.023)$) edge [color = blue!90!black,thick] ($(v3)+(0,.023)$);
			\draw ($(v2)-(0,0.023)$) edge [color = red!90!black, thick] ($(v3)-(0,.023)$);
			
			\draw ($(v1)+(0,0.03)$) edge [color = blue!90!black,thick] ($(v0)+(0,.03)$);
			\draw ($(v1)-(0,0.03)$) edge [color = red!90!black, thick] ($(v0)-(0,.03)$);
			
			\draw (v2)edge [color = blue!90!black, thick] (v1);		
			
			\draw (v2)edge [color = red!90!black, thick] (v4);			
			\draw (v1) edge [color = red!90!black, thick] (v3);						
			\draw (v4) edge [color = red!90!black, thick] (v1);		
			\draw (v2)edge [color = red!90!black, thick] (v0);
			
			\foreach \i/\j in {3/4,0/4}{
				\draw (v\i) edge [color = red!90!black, thick, dashed] (v\j);
			}							
			\foreach \i in {0,...,4}{
				\fill (v\i) circle (1.5pt);
			}		
			\end{tikzpicture}\\

	\begin{tikzpicture}		
	\node at (-.8,.7){{\tiny $(d)$}};
	%	\foreach \i in {0,...,4}{
	%		\coordinate (v\i) at (90+\i*72:.7cm);
	%	}		
	\qedge{(v1)}{(v1)}{(v2)}{4pt}{.4pt}{black}{black!30!white,opacity=0.2};
	\qedge{(v3)}{(v0)}{(v4)}{4pt}{.4pt}{black}{black!30!white,opacity=0.2};

	\draw ($(v2)+(0,0.023)$) edge [color = blue!90!black,thick] ($(v3)+(0,.023)$);
	\draw ($(v2)-(0,0.023)$) edge [color = red!90!black, thick] ($(v3)-(0,.023)$);

		\draw ($(v1)+(0,0.02)$) edge [color = blue!90!black,thick, dashed] ($(v4)+(0,.02)$);
		\draw ($(v1)-(0,0.02)$) edge [color = red!90!black, thick] ($(v4)-(0,.02)$);

		\draw ($(v0)+(0,0.03)$) edge [color = blue!90!black,thick, dashed] ($(v1)+(0,.03)$);
		\draw ($(v0)-(0,0.03)$) edge [color = red!90!black, thick] ($(v1)-(0,.03)$);

	\draw ($(v2)+(0,0.03)$) edge [color = blue!90!black,thick, dashed] ($(v4)+(0,.03)$);
	\draw ($(v2)-(0,0.03)$) edge [color = red!90!black, thick] ($(v4)-(0,.03)$);
	
	\draw ($(v2)-(0.025,0)$) edge [color = blue!90!black,thick, dashed] ($(v0)-(.025,0)$);
	\draw ($(v2)+(0.025,0)$) edge [color = red!90!black, thick] ($(v0)+(.025,0)$);
	
	\draw ($(v1)+(0,0.03)$) edge [color = blue!90!black,thick] ($(v3)+(0,.03)$);
	\draw ($(v1)-(0,0.03)$) edge [color = red!90!black, thick] ($(v3)-(0,.03)$);

	\draw ($(v2)+(0.025,0)$) edge [color = red!90!black, thick] ($(v0)+(.025,0)$);
	
		\foreach \i/\j in {0/3,3/4}{
			\draw (v\i) edge [color = blue!90!black, thick, dashed] (v\j);
		}	
							
	\foreach \i in {0,...,4}{
		\fill (v\i) circle (1.5pt);
	}		
	\end{tikzpicture}\\
		
		\begin{tikzpicture}	
		\node at (-.8,.7){{\tiny $(e)$}};	
		%	\foreach \i in {0,...,4}{
		%		\coordinate (v\i) at (90+\i*72:.7cm);
		%	}		
		\foreach \i/\j in {0/1,0/4}{
			\draw (v\i) edge [color = red!90!black,thick] (v\j);
		}	
		\draw ($(v2)+(0,0.023)$) edge [color = blue!90!black,thick] ($(v3)+(0,.023)$);
		\draw ($(v2)-(0,0.023)$) edge [color = red!90!black, thick] ($(v3)-(0,.023)$);
		
		\draw ($(v1)-(0.025,0)$) edge [color = red!90!black,thick] ($(v2)-(.025,0)$);
		\draw ($(v1)+(0.025,0)$) edge [color = blue!90!black, thick, dashed] ($(v2)+(.025,0)$);
		
		\draw ($(v3)+(0.025,0)$) edge [color = red!90!black,thick] ($(v4)+(.025,0)$);
		\draw ($(v3)-(0.025,0)$) edge [color = blue!90!black, thick, dashed] ($(v4)-(.025,0)$);
		
		\draw (v2) edge [color = blue!90!black, thick, dashed] (v4);
		\draw (v3) edge [color = blue!90!black, thick, dashed] (v1);					
		\foreach \i in {0,...,4}{
			\fill (v\i) circle (1.5pt);
		}		
		\end{tikzpicture}\\

		\begin{tikzpicture}
		\node at (-.8,.7){{\tiny $(f)$}};		
%		\foreach \i in {0,...,4}{
%			\coordinate (v\i) at (90+\i*72:.7cm);
%		}		
		\foreach \i/\j in {0/3,1/4,1/3,2/4}{
			\draw (v\i) edge [color = blue!90!black, thick, dashed] (v\j);
		}			
		\draw ($(v2)+(0,0.023)$) edge [color = blue!90!black,thick] ($(v3)+(0,.023)$);
		\draw ($(v2)-(0,0.023)$) edge [color = red!90!black, thick] ($(v3)-(0,.023)$);
		
		\draw ($(v0)-(0,0.03)$) edge [color = blue!90!black,thick] ($(v4)-(0,.03)$);
		\draw ($(v0)+(0,0.03)$) edge [color = red!90!black, thick] ($(v4)+(0,.03)$);
	
		\draw ($(v1)-(0.025,0)$) edge [color = red!90!black,thick] ($(v2)-(.025,0)$);
		\draw ($(v1)+(0.025,0)$) edge [color = blue!90!black, thick, dashed] ($(v2)+(.025,0)$);
		
		\draw ($(v0)-(0.035,0)$) edge [color = red!90!black,thick] ($(v1)-(.035,0)$);
		\draw ($(v0)+(0.035,0)$) edge [color = blue!90!black, thick, dashed] ($(v1)+(.035,0)$);
	
		\draw ($(v3)+(0.025,0)$) edge [color = red!90!black,thick] ($(v4)+(.025,0)$);
		\draw ($(v3)-(0.025,0)$) edge [color = blue!90!black, thick, dashed] ($(v4)-(.025,0)$);
								
		\foreach \i in {0,...,4}{
			\fill (v\i) circle (1.5pt);
		}		
		\end{tikzpicture}\\

		\begin{tikzpicture}	
		\node at (-.8,.7){{\tiny $(g)$}};	
	%	\foreach \i in {0,...,4}{
	%		\coordinate (v\i) at (90+\i*72:.7cm);
	%	}		
		\foreach \i/\j in {0/1,1/2,2/3,3/4,0/4}{
			\draw (v\i) edge [color = red!90!black, thick] (v\j);
		}	
		\foreach \i/\j in {1/3,1/4,0/2,0/3,2/4}{
			\draw (v\i) edge [color = blue!90!black, thick] (v\j);
		}			
		\foreach \i in {0,...,4}{
			\fill (v\i) circle (1.5pt);
		}		
		\end{tikzpicture}\\

	\end{multicols}
	
	\caption{The illustration to the proofs of Facts \ref{f:K23} and \ref{f:C5}.}
	\label{fig:ff3}
	\vspace{-1em}
\end{figure}
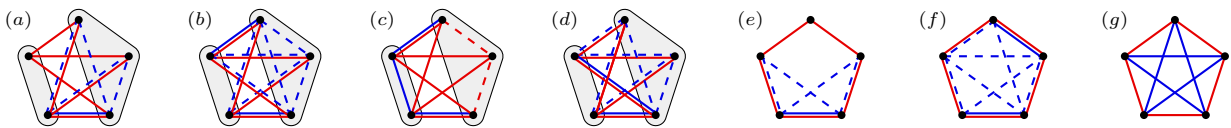

\begin{fact}\label{f:K23}
	Let $\rr$ and $\bb$ be two graphs on the same 5-vertex set, such that $\RB\nsubseteq\rr\cup \bb$. If $K_{2,3}\subseteq \rr$, then $|\bb|\le 4$ and either $\bb\subseteq T\cup \{e\}$ or $|\rr|+|\bb|\le 11$.
\end{fact}
\begin{proof}
	 We let $\crr{K_{2,3}}\subseteq \rr$ and $\bb\nsubseteq T\cup \{e\}$, since otherwise $|\bb|\le 4$, and the assertion follows. Note that due to $\RB\nsubseteq \rr\cup \bb$, whenever $|\crr{K_{2,3}}\cap \bb|=1$, then four pairs of ${V\choose 2}$, shown in Figure~\ref{fig:ff3}$(a)$ with dashed lines, are forbidden for $\bb$. In particular $|\crr{K_{2,3}}\cap \bb|\le 1$ causes $\bb\subseteq T\cup \{e\}$, and thus we may assume $|\crr{K_{2,3}}\cap \bb|\ge 2$. Further, $M_2\subseteq\crr{K_{2,3}}\cap \bb$ entails $|\bb|\le 3$ (see Figure \ref{fig:ff3}$(b)$) and $|\bb|=3$ yields $|\rr|\le 8$ (see Figure \ref{fig:ff3}$(c)$). Therefore in this case, either $\bb\subseteq T\cup \{e\}$ or $|\rr| + |\bb|\le 11$, as required. Finally, if $P_2\subseteq\crr{K_{2,3}}\cap \bb$, then $\bb\subseteq T\cup \{e\}$ (see Figure \ref{fig:ff3}$(d)$), and the assertion follows again.
\end{proof}

%\vspace{-3em}

\begin{fact}\label{f:C5}
	Let $\rr$ and $\bb$ be two graphs on the same 5-vertex set, such that $\RB\nsubseteq\rr\cup \bb$. If $C_5\subseteq \rr$ and $|\rr|\ge 6$, then $|\bb|\le 4$.
\end{fact}
\begin{proof}
		We let $\crr{C_5}\subseteq \rr$. 		
		Now, if $|\crr{C_5}\cap \bb|= 1$ then $|\crr{C}^c_{\crr{5}}\cap \bb|\le 3$ and thus $|\bb|\le 4$, as required (see Figure \ref{fig:ff3}$(e)$).
		Further, for $|\crr{C_5}\cap \bb|\ge 2$ we have $|\bb|\le 3$ (see Figure \ref{fig:ff3}$(f)$) and we are done again. Finally, let $\crr{C_5}\cap \bb= \emptyset$, that is $\bb\subseteq \crr{C}^c_{\crr{5}}$. Then $|\bb|\ge 5$ entails $\bb=\crr{C}^c_{\crr{5}}$ (see Figure \ref{fig:ff3}$(g)$). This, in turn, due to the symmetry, yields $\rr = \crr{C_5}$, contradicting $|\rr|\ge 6$.
\end{proof}

 It turns out that if two graphs, $\rr$ and $\bb$, on the same 5-vertex set do not contain $\RB$, then $|\rr|+|\bb|\le 13$.

\begin{lemma}\label{l:five}
	Let $\rr$ and $\bb$ be two graphs on the same vertex set $V=\{v,a,b,x,y\}$, such that $\RB\nsubseteq\rr\cup \bb$. Then $|\rr|+|\bb|\le13$ and, if $|\rr|+|\bb|\ge 12$, $|\rr|\ge|\bb|$, then up to the isomorphism one of the following holds, (see Figure \ref{fig:edge12}),
		\begin{enumerate}[label=\Alabel]
			\item\label{it:typa} $\rr\subseteq K_5[V]- \{ab\}$, $\bb\subseteq T\cup \{ab\}$, where $T=K_3[\{v,x,y\}]=\{vx,vy,xy\}$;
			\item\label{it:typb} $\rr=K_5[V]$, $\bb = \{ab,xy\}$;
			\item\label{it:typc} $\rr=\bb=K_4[\{a,b,x,y\}]$, where $K_4$ is a complete graph on the vertex set $\{a,b,x,y\}$;
			\item\label{it:typd} $\rr=S_5\cup \{ab,xy\}$, $\bb=S_5\cup \{ax,by\}$, where $S_5=\{va,vb,vx,vy\}$.
		\end{enumerate}
\end{lemma}
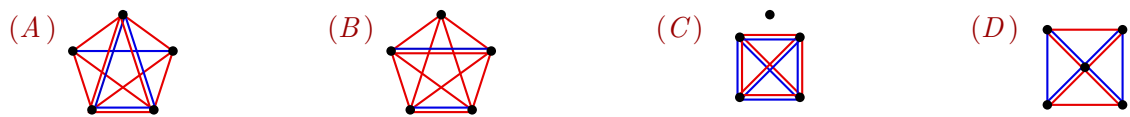
\begin{figure}[h!]
	\centering			
	\begin{multicols}{4}
		
		\begin{tikzpicture}
		
		\node at (-1.2,.5) {\ref{it:typa}};
		
					\foreach \i in {0,...,4}{
						\coordinate (v\i) at (90+\i*72:.7cm);
					}
					
					\foreach \i/\j in {0/1,0/4,1/2,1/3,2/4,3/4}{
						\draw ($(v\i)+(0,0)$) edge [color = red!90!black,  thick] ($(v\j)+(0,0)$);
					}
					
					\draw ($(v1)+(0,0)$) edge [color = blue!90!black, thick] ($(v4)+(0,0)$);
					
					\draw ($(v0)+(-.045,-.045)$) edge [color = red!90!black, thick] ($(v2)+(-.045,-.045)$);
					\draw ($(v0)+(.045,.045)$) edge [color = blue!80!black,  thick] ($(v2)+(.045,.045)$);
					
					\draw ($(v0)+(-.03,-.03)$) edge [color = red!90!black,  thick] ($(v3)+(-.03,-.03)$);
					\draw ($(v0)+(.02,.02)$) edge [color = blue!90!black,  thick] ($(v3)+(.02,.02)$);
					
					\draw ($(v2)+(0,-.03)$) edge [color = red!90!black,  thick] ($(v3)+(0,-.03)$);
					\draw ($(v2)+(0,.03)$) edge [color = blue!90!black,  thick] ($(v3)+(0,.03)$);	
					
	%				\node [anchor = south west] at (v0) {$a$};
	%				\node [anchor = north east] at (v2) {$b$};
	%				\node [anchor = north west] at (v3) {$c$};
	%				\node [anchor = south west] at (v4) {$y$};		
	%				\node [anchor = south east] at (v1) {$x$};
					
					\foreach \i in {0,...,4}{
						\fill (v\i) circle (1.8pt);
					}
			
		\end{tikzpicture}\\

		\begin{tikzpicture}
			\node at (-1.2,.5) {\ref{it:typb}};

		\foreach \i in {0,...,4}{
			\coordinate (v\i) at (90+\i*72:.7cm);
		}
		
		\foreach \i/\j in {0/1,0/4,1/2,1/3,2/4,3/4,0/2,0/3}{
			\draw ($(v\i)+(0,0)$) edge [color = red!90!black,  thick] ($(v\j)+(0,0)$);
		}
		
	%	\draw ($(v1)+(0,0)$) edge [color = red!90!black, line width=2pt] ($(v4)+(0,0)$);
		
		\draw ($(v1)+(0,-.03)$) edge [color = red!90!black, thick] ($(v4)+(0,-.03)$);
		\draw ($(v1)+(0,.03)$) edge [color = blue!80!black,  thick] ($(v4)+(0,.03)$);
		
		\draw ($(v2)+(0,-.03)$) edge [color = red!90!black,  thick] ($(v3)+(0,-.03)$);
		\draw ($(v2)+(0,.03)$) edge [color = blue!90!black,  thick] ($(v3)+(0,.03)$);
		
%		\draw ($(v2)+(0,-.04)$) edge [color = red!90!black, line width=2pt] ($(v3)+(0,-.04)$);
%		\draw ($(v2)+(0,.04)$) edge [color = blue!90!black, line width=2.3pt] ($(v3)+(0,.04)$);	

		\foreach \i in {0,...,4}{
			\fill (v\i) circle (1.8pt);
		}

		%		\node [anchor = south west] at (v0) {$a$};
		%		\node [anchor = north east] at (v2) {$b$};
		%		\node [anchor = north west] at (v3) {$c$};
	%	\node [anchor = south west] at (v4) {$y$};		
	%	\node [anchor = south east] at (v1) {$x$};
		
		\end{tikzpicture}\\

	\begin{tikzpicture}
	
	\node at (-1.2,.5) {\ref{it:typc}};
	
%	\foreach \i in {0,...,4}{
%		\coordinate (v\i) at (90+\i*72:1.5cm);
%	}

		\coordinate (v0) at (0,.7);
		\coordinate (v1) at (-.4,.4);
		\coordinate (v2) at (-.4,-.4);
		\coordinate (v3) at (.4,-.4);
		\coordinate (v4) at (.4,.4);

	\foreach \i/\j in {1/3,2/4}{
		\draw ($(v\i)+(0,.04)$) edge [color = red!90!black,  thick] ($(v\j)+(0,.04)$);
		\draw ($(v\i)-(0,.04)$) edge [color = blue!90!black,  thick] ($(v\j)-(0,.04)$);
	}

	\foreach \i/\j in {1/4,3/2}{
		\draw ($(v\i)+(0,.03)$) edge [color = red!90!black,  thick] ($(v\j)+(0,.03)$);
		\draw ($(v\i)-(0,.03)$) edge [color = blue!90!black,  thick] ($(v\j)-(0,.03)$);
	}
	
\draw ($(v1)+(.03,0)$) edge [color = red!90!black,  thick] ($(v2)+(.03,0)$);
\draw ($(v3)+(.03,0)$) edge [color = red!90!black, thick] ($(v4)+(.03,0)$);

\draw ($(v1)-(.03,0)$) edge [color = blue!90!black,  thick] ($(v2)-(.03,0)$);
\draw ($(v3)-(.03,0)$) edge [color = blue!90!black,  thick] ($(v4)-(.03,0)$);

	%			\node [anchor = south west] at (v0) {$a$};
	%			\node [anchor = north east] at (v2) {$b$};
	%			\node [anchor = north west] at (v3) {$c$};
	%			\node [anchor = south west] at (v4) {$y$};		
	%			\node [anchor = south east] at (v1) {$x$};
	
	\foreach \i in {0,...,4}{
		\fill (v\i) circle (1.8pt);
	}
	
	\end{tikzpicture}

			\begin{tikzpicture}
				\node at (-1.2,.5) {\ref{it:typd}};
				
			\coordinate (v0) at (0,0);
			\coordinate (v1) at (-.5,.5);
			\coordinate (v2) at (-.5,-.5);
			\coordinate (v3) at (.5,-.5);
			\coordinate (v4) at (.5,.5);	
			
			\foreach \i in {1,...,4}{
				\draw ($(v\i)+(0,.04)$) edge [color = blue!90!black,  thick] ($(v0)+(0,.04)$);
				\draw ($(v\i)-(0,.04)$) edge [color = red!90!black,  thick] ($(v0)-(0,.04)$);
			}
			
			\draw (v1) edge [color = blue!90!black,  thick] (v2);
			\draw (v3) edge [color = blue!90!black,  thick] (v4);
			
			\draw (v1) edge [color = red!90!black, thick] (v4);
			\draw (v2) edge [color = red!90!black, thick] (v3);

			\foreach \i in {0,...,4}{
				\fill (v\i) circle (1.8pt);
			}
			\end{tikzpicture}
				
	\end{multicols}
	\caption{All $\rr\cup\bb$ on 5 vertices and with $|\rr|+|\bb|\ge 12$, such that $\RB\nsubseteq \rr\cup \bb$.}	\label{fig:edge12}
	\vspace{-1em}
\end{figure}
\begin{proof}
	Let two graphs, $\rr$ and $\bb$, on the same vertex set $V=\{v,a,b,x,y\}$, with $\RB\nsubseteq \rr\cup \bb$ be given. Moreover, let $|\rr|+|\bb|\ge 12$, $|\rr|\ge|\bb|$, and thereby $2\le |\bb|\le |\rr|\le 10$ and $|\rr|\ge 6$. We will show that one of \ref{it:typa}-\ref{it:typd} occurs. In what follows we assume that $\bb\nsubseteq M_2$, because otherwise \ref{it:typb} holds.
	
	First observe that $|\rr|\ge 8$ entails $K_{2,3}\subseteq \rr$. Then Fact \ref{f:K23} combined with $|\rr|+|\bb| \ge 12$ tells us that $\bb\subseteq T\cup \{e\}$. Moreover, $\bb\nsubseteq M_2$ yields $|\bb\cap T|\ge 2$. But $|\rr|\ge 8$ and thus there are at least 4 edges of $\rr$ between $V[T]$ and $e$. Therefore, to avoid $\RB\subseteq \rr\cup \bb$, we have $e\notin \rr$, and hence \ref{it:typa} follows.
		
	Further, for $|\rr|\le 7$ we have $|\bb|\ge 5$ and thus Facts \ref{f:K23} and \ref{f:C5} yield that $\rr$ contains neither $K_{2,3}$ nor $C_5$. If $\crr{K_4}\subseteq \rr$, then to avoid $\RB$ in $\rr\cup \bb$, every edge $e\in \bb$ with $|e\cap V[\crr{K_4}]|= 1$ is an isolated edge in $\bb$ (see Figure \ref{fig:f3}$(a)$), entailing $|\bb|\le 4$. Therefore $\bb \subseteq \crr{K_4}$ and hence, using again $\RB\nsubseteq \rr\cup \bb$, also $\rr\subseteq \crr{K_4}$, yielding \ref{it:typc}.		
	
\begin{figure}[h!]
	\centering			
	\begin{multicols}{7}

		\begin{tikzpicture}		
		\node at (-.8,.7){{\tiny $(a)$}};
		\foreach \i in {0,...,4}{
			\coordinate (v\i) at (90+\i*72:.7cm);
		}		
		\foreach \i/\j in {2/3,2/4,3/4}{
			\draw (v\i) edge [color = red!90!black, thick] (v\j);
		}
		\draw [blue!90!black, thick, dashed] ($(v1)+(.02,0)$) edge ($(v2)+(.02,0)$);
		\draw [red!90!black, thick] ($(v1)-(.02,0)$) edge ($(v2)-(.02,0)$);
		
		\draw [blue!90!black, thick, dashed] ($(v1)+(.02,.01)$) edge ($(v3)+(.02,.01)$);
		\draw [red!90!black, thick] ($(v1)-(.02,.01)$) edge ($(v3)-(.02,.01)$);
		
		\draw [blue!90!black, thick, dashed] ($(v1)+(0,.02)$) edge ($(v4)+(0,.02)$);
		\draw [red!90!black, thick] ($(v1)-(0,.02)$) edge ($(v4)-(0,.02)$);
		
		\foreach \i/\j in {0/2,0/3,0/4}{
			\draw (v\i) edge [color = blue!90!black, thick, dashed] (v\j);
		}
		
		\foreach \i/\j in {2/3,2/4,3/4}{
			\draw (v\i) edge [color = red!90!black, ultra thick] (v\j);
		}			
			\draw (v0) edge [color = blue!90!black, thick] (v1);
		%	\draw ($(v2)+(0,0)$) edge [color = blue!90!black, line width=1.5pt, dotted] ($(v1)+(0,0)$);					
		\foreach \i in {0,...,4}{
			\fill (v\i) circle (1.5pt);
		}		
		\end{tikzpicture}\\

		\begin{tikzpicture}		
		\node at (-.8,.7){{\tiny $(b)$}};
	%	\foreach \i in {0,...,4}{
	%		\coordinate (v\i) at (90+\i*72:.7cm);
	%	}		
		\foreach \i/\j in {0/1,0/2,0/3,0/4}{
			\draw (v\i) edge [color = red!90!black, thick] (v\j);
		}			
		\draw (v2) edge [color = blue!90!black, thick] (v3);
		
		\foreach \i/\j in {2/1,2/4,3/1,3/4}{
			\draw (v\i) edge [color = blue!90!black, thick, dashed] (v\j);
		}			
		
		\node at ($(v0)+(.15,.05)$) {\tiny $v$};
		%	\draw ($(v2)+(0,0)$) edge [color = blue!90!black, line width=1.5pt, dotted] ($(v1)+(0,0)$);					
		\foreach \i in {0,...,4}{
			\fill (v\i) circle (1.5pt);
		}		
		\end{tikzpicture}\\

		\begin{tikzpicture}		
		\node at (-.8,.7){{\tiny $(c)$}};
%		\foreach \i in {0,...,4}{
%			\coordinate (v\i) at (90+\i*72:.7cm);
%		}		
		\foreach \i/\j in {1/4,2/3}{
			\draw (v\i) edge [color = blue!90!black, thick] (v\j);
		}			
		\draw ($(v0)-(0,0.03)$) edge [color = blue!90!black,thick] ($(v1)-(0,.03)$);
		\draw ($(v0)+(0,0.03)$) edge [color = red!90!black, thick] ($(v1)+(0,.03)$);
		
		\draw ($(v0)+(0.025,0)$) edge [color = blue!90!black,thick] ($(v2)+(.025,0)$);
		\draw ($(v0)-(0.025,0)$) edge [color = red!90!black, thick] ($(v2)-(.025,0)$);
		
		\draw ($(v0)+(0.025,0)$) edge [color = blue!90!black,thick] ($(v3)+(.025,0)$);
		\draw ($(v0)-(0.035,0)$) edge [color = red!90!black, thick] ($(v3)-(.025,0)$);
		
		\draw ($(v0)-(0,0.03)$) edge [color = blue!90!black,thick] ($(v4)-(0,.03)$);
		\draw ($(v0)+(0,0.03)$) edge [color = red!90!black, thick] ($(v4)+(0,.03)$);
							
		\foreach \i in {0,...,4}{
			\fill (v\i) circle (1.5pt);
		}		
		\end{tikzpicture}\\

	\begin{tikzpicture}		
	\node at (-.5,1.2){{\tiny $(d)$}};
	
	\coordinate (a) at (0,0);
	\coordinate (b) at (.7,0);
	\coordinate (x) at (0,.7);
	\coordinate (y) at (.7, .7);
	\coordinate (v) at (1,1);
	
	\node at ($(a)+(-.1,-.1)$) {\tiny $a$};
	\node at ($(b)+(.1,-.1)$) {\tiny $b$};
	\node at ($(x)+(-.1,.1)$) {\tiny $x$};
	\node at ($(y)+(.12,-.1)$) {\tiny $y$};
	\node at ($(v)+(.15,.05)$) {\tiny $v$};
	
	\foreach \i/\j in {a/b,b/y,x/y,a/x,b/x,v/y}{
		\draw (\i) edge [color = red!90!black, thick] (\j);
	}			
	
	\foreach \i in {a,b,x,y,v}{
		\fill (\i) circle (1.5pt);
	}		
	\end{tikzpicture}
	
	\begin{tikzpicture}		
	\node at (-.5,1.2){{\tiny $(e)$}};

	\foreach \i/\j in {a/b,b/y,x/y,a/x,b/x,v/y}{
		\draw (\i) edge [color = red!90!black, thick] (\j);
	}	
	
	\foreach \i/\j in {a/b,a/x,b/x}{
		\draw (\i) edge [color = blue!90!black, thick] (\j);
	}								
	\foreach \i in {a,b,x,y,v}{
		\fill (\i) circle (1.5pt);
	}		
	\end{tikzpicture}\\
	
	\begin{tikzpicture}		
	\node at (-.5,1.2){{\tiny $(f)$}};

	\foreach \i/\j in {a/b,b/y,x/y,a/x,b/x,v/y}{
		\draw (\i) edge [color = red!90!black, thick] (\j);
	}			
	\foreach \i/\j in {x/y,x/v,y/v}{
		\draw (\i) edge [color = blue!90!black, thick] (\j);
	}	
	\foreach \i in {a,b,x,y,v}{
		\fill (\i) circle (1.5pt);
	}		
	\phantom{\node[anchor = north west] at (b) {\tiny $b$};}
	\end{tikzpicture}
	
	\begin{tikzpicture}		
	\node at (-.5,1.2){{\tiny $(g)$}};

	\foreach \i/\j in {a/b,b/y,x/y,a/x,b/x,v/y}{
		\draw (\i) edge [color = red!90!black, thick] (\j);
	}			
	\foreach \i/\j in {b/y,b/v,y/v}{
		\draw (\i) edge [color = blue!90!black, thick] (\j);
	}	
	\foreach \i in {a,b,x,y,v}{
		\fill (\i) circle (1.5pt);
	}		
	\end{tikzpicture}

	\end{multicols}
	
	\caption{The illustration to the proof of Lemma \ref{l:five}.}
	\label{fig:f3}
	\vspace{-1em}
\end{figure}

%\vspace{-3em}		
	
	Now, as every 5-vertex graph with at least 7 edges contains at least one of the graphs, $K_{2,3}$, $C_5$, or $K_4$, as a subgraph, we may assume that $|\rr|\le 6$ and thereby $|\bb|=|\rr|= 6$. First consider $\Delta(\rr)=4$ and let $\deg_{\rr}(v)=4$. Note that $P_2\nsubseteq \bb[V\setminus \{v\}]$ (see Figure \ref{fig:f3}$(b)$). Thus $|\bb[V\setminus\{v\}]|\le 2$, and $|\bb|=6$ entails $S_5\subseteq\bb$, where $S_5=\{va,vb,vx,vy\}$, and $\bb[V\setminus\{v\}] = M_2$ (see Figure \ref{fig:f3}$(c)$). By the symmetry, $\rr[V\setminus \{v\}]\subseteq M_2$ and, to avoid a copy of $\RB$, $\rr \cap \bb[V\setminus\{v\}] =\emptyset$, yielding  \ref{it:typd}.
	
	Finally we let $\Delta(\rr)\le 3$ and $|\rr|=|\bb|=6$. The only (up to the isomorphism) $\{K_{2,3},C_5,K_4, S_5\}$-free graph $G=\{ab,by,xy,ax,bx,yv\}$ with 6 edges on the vertex set $V$ is given in Figure \ref{fig:f3}$(d)$. Observe that any two edges of one of the triangles $abx$, $xyv$ or $byv$ given in Figure \ref{fig:f3}$(e)$-$(g)$ in blue, create, together with $\rr$, a copy of $\RB$. Therefore $|\bb|\le 5$, a contradiction.
\end{proof}

\begin{lemma}\label{11}
	Let $\rr$ and $\bb$ be two graphs on the same 5-vertex set, such that $\RB\nsubseteq\rr\cup \bb$. If  $\Delta(\rr),\Delta (\bb) \le 3$, and at least three vertices of both $\rr$ and $\bb$ have degree at most $2$, then $|\rr|+|\bb|\le 10$.
\end{lemma}

\begin{proof}
	For the sake of contradiction assume that  $|\rr|+|\bb|\ge 11$ and let $|\rr|\ge 6$. Owing to the degree restriction we also have $\max\{|\rr|,|\bb|\}\le6$, so $5\le |\bb|\le |\rr|=6$. There are exactly two 5-vertex graphs with the degree sequence $(2,2,2,3,3)$: a pentagon $C_5$ with one diagonal, and $K_{2,3}$. But then, in view of Facts \ref{f:K23} and \ref{f:C5}, $|\bb|\le 4$, a contradiction.
\end{proof}

%%%%%%%%%%%%%%%%%%%%%%%%%%%%%%%%%%%%%%%%%%%%%%%%%%%%%%%%%%%%%%%%%%%%%%
%																																		%
%																												   						%
%   											Seven vertices         							    								%
%																																		%
%																																		%
%%%%%%%%%%%%%%%%%%%%%%%%%%%%%%%%%%%%%%%%%%%%%%%%%%%%%%%%%%%%%%%%%%%%%%

\subsection{Proof of Lemma \ref{l:n7}}\label{7}
	Let $\cH$ be a $\p$-free 3-graph on a 7-vertex set $V$ and with at least 20 edges. We will show that $\cH=K_6^{(3)}\cup K_1$, which will end the proof of Lemma \ref{l:n7}. To this end pick two vertices, $x,y\in V$, with the largest pair degree $\deg_\cH(x,y)=\Delta_2(\cH)$ and set $Z=V\setminus\{x,y\}$. We let
		\[
				\rr=L_{\cH}(x)[Z]\qand \bb=L_{\cH}(y)[Z]
		\]
 be the link graphs of $x$ and $y$, respectively, induced on $Z$. Then,
		\begin{equation}\label{eq:20}
			|\cH| = \deg_\cH(x,y) + |\rr| + |\bb| + |\cH[Z]|\ge 20.
		\end{equation}
	Moreover we have $3\le \deg_\cH(x,y) \le 5$. Indeed, the upper bound is a trivial consequence of $|Z|=5$, while the lower bound follows from $\sum_{x,y\in V}\deg_\cH(x,y)=3|\cH|\ge 60$. 	

	We start with estimating the number of edges in the 3-graph $\cH[Z]$ induced on $Z$.
		\begin{claim}\label{cl:z}
			\begin{enumerate}[label=\rmlabel]
				\item\label{it:z1} If $\deg_\cH(x,y)=5$, then $|\cH[Z]|\le 2$.
				\item\label{it:z2} If $\deg_\cH(x,y)=4$, then $|\cH[Z]|\le 4$. Moreover, if additionally $|\cH[Z]| = 4$, then $\cH[Z] = K^{(3)}_4[N_\cH(x,y)]$ is a complete 3-graph on the vertex set $N_\cH(x,y)$.
				\item\label{it:z3} If $\deg_\cH(x,y)=3$, then $|\cH[Z]|\le 6$.
			\end{enumerate}
		\end{claim}
		\begin{proof}
			Clearly, if $|\cH[Z]|\ge3$, then there
			are in $\cH[Z]$ two edges sharing two vertices,
			say, $abc$ and $bcd$. Set $z$ for the unique element of $Z\setminus\{a,b,c,d\}$. Observe that if both $z$ and $a$ are common neighbors of $x,y$, then the sequence $zxyabcd$ is a minimal 4-path in $\cH$ (see Figure~\ref{fig:f5}$(a)$). As for $\deg_\cH(x,y)=5$ each vertex of $Z$ is a common neighbor of $x,y$, the above observation establishes \ref{it:z1}.			
			
\begin{figure}[h!]
	\centering			
	\begin{multicols}{3}

			\begin{tikzpicture}		
			\node at (-1,.8){ $(a)$};
			
			\foreach  \i  in {0,...,4}{
				\coordinate (v\i) at (.6*\i,0);}
				
			\coordinate (x) at (.1,.8);
			\coordinate (y) at (.5,.8);
			
			\qedge{(x)}{(y)}{(v0)}{5pt}{1pt}{red!0!black}{black!20!white,opacity=0.2};
			\qedge{(x)}{(y)}{(v1)}{5pt}{1pt}{red!0!black}{black!20!white,opacity=0.2};
			\qedge{(v3)}{(v2)}{(v1)}{5pt}{1pt}{red!0!black}{black!20!white,opacity=0.2};
			\qedge{(v4)}{(v3)}{(v2)}{5pt}{1pt}{red!0!black}{black!20!white,opacity=0.2};

			\fill [red!0!black] (x) circle (2pt);
			\fill [blue!0!black] (y) circle (2pt);
			
			\foreach \i in {0,...,4}{
				\fill (v\i) circle (2pt);
			}

			\node[anchor=mid] at ($(v0)+(0,-.3)$){\tiny $z$};
			\node[anchor=mid] at ($(v1)+(0,-.3)$){\tiny $a$};
			\node[anchor=mid] at ($(v2)+(0,-.3)$){\tiny $b$};
			\node[anchor=mid] at ($(v3)+(0,-.3)$){\tiny $c$};
			\node[anchor=mid] at ($(v4)+(0,-.3)$){\tiny $d$};
			\node at ($(x)+(-.25,.2)$){\tiny $x$};
			\node  at ($(y)+(.25,.2)$){\tiny $y$};

			\end{tikzpicture}
		
		\begin{tikzpicture}		
			\node at (-1,.8){$ (b)$};
				
				\qedge{(x)}{(y)}{(v0)}{5pt}{1pt}{red!0!black}{black!20!white,opacity=0.2};
				\qedge{(x)}{(y)}{(v1)}{5pt}{1pt}{red!0!black}{black!20!white,opacity=0.2};
				
				\draw [ green!60!black, very thick]
				(v0) edge [bend right = 30] (v4)
				(v0) edge (v1);

			\foreach \i in {0,...,4}{
				\fill (v\i) circle (2pt);
			}		
			
				\fill [red!0!black] (x) circle (2pt);
				\fill [blue!0!black] (y) circle (2pt);
			
		\end{tikzpicture}\\

		\begin{tikzpicture}		
		\node at (-1,.8){$(c)$};
		
		\coordinate (x) at (.7,.8);
		\coordinate (y) at (1.1,.8);
		
		\qedge{(x)}{(y)}{(v0)}{5pt}{1pt}{red!0!black}{black!20!white,opacity=0.2};
		\qedge{(x)}{(y)}{(v1)}{5pt}{1pt}{red!0!black}{black!20!white,opacity=0.2};
			\qedge{(x)}{(y)}{(v2)}{5pt}{1pt}{red!0!black}{black!20!white,opacity=0.2};
			\qedge{(x)}{(y)}{(v3)}{5pt}{1pt}{red!0!black}{black!20!white,opacity=0.2};

		\draw [ green!60!black, very thick]
		(v0) edge [bend right = 30] (v4)
		(v1) edge [bend right = 30] (v4)
		(v2) edge [bend right = 30] (v4)
		(v3) edge [bend right = 30] (v4);

		\foreach \i in {0,...,4}{
			\fill (v\i) circle (2pt);
		}		
		
		\fill [red!0!black] (x) circle (2pt);
		\fill [blue!0!black] (y) circle (2pt);
		
		\end{tikzpicture}

	\end{multicols}
	\vspace{-2em}

	\caption{The illustration to the proof of Claim \ref{cl:z}.}
	\label{fig:f5}
	\vspace{-1em}
\end{figure}
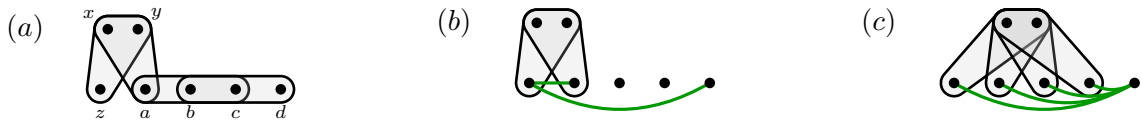		

			For the proof of \ref{it:z2}, instead of looking at edges $e\in \cH[Z]$, we will look at their \emph{complement edges} $e^c = Z\setminus e$ in $Z$. (For example the green 2-edges in Figure \ref{fig:f5}$(b)$ are complement edges of the 3-edges $abc,bcd\in \cH[Z]$ in Figure \ref{fig:f5}$(a)$.) In view of this definition, the above observation reads as follows. If there are two adjacent complement edges of $\cH[Z]$ such that at least one of them is contained in $N_\cH(x,y)$, then $\cH$ contains a minimal 4-path (see Figure \ref{fig:f5}$(a),(b)$). Therefore if $|\cH[Z]|\ge 4$, then all complement edges contain the unique vertex of $Z\setminus N_\cH(x,y)$ (see Figure \ref{fig:f5}$(c)$) and thereby $\cH[Z]=K^{(3)}_4[N_\cH(x,y)]$ is a complete 3-graph on the vertex set $N_\cH(x,y)$.
			
			Finally, to prove \ref{it:z3} note that \eqref{eq:20} together with $\deg_\cH(x,y)=\Delta_2(\cH)=3$ entails
				\[
					17 + 2|\cH[Z]|\le |\rr| + |\bb|+3|\cH[Z]| = \sum_{a,b\in Z} \deg_\cH(a,b)\le {5\choose 2}\cdot\Delta_2(\cH)=30. \qedhere
				\]
		\end{proof}	
		
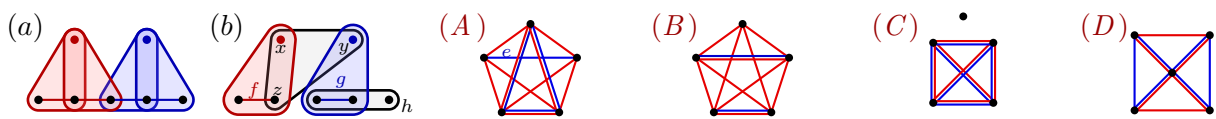
\begin{figure}[h!]
	\centering			
	\begin{multicols}{6}

			\begin{tikzpicture}	[scale = .8]
			\node at (-.3,1.2){{ $(a)$}};
			
			\coordinate (x) at (.6,1);
			\coordinate (y) at (1.8,1);
			
			\foreach  \i  in {0,...,4}{
				\coordinate (v\i) at (.6*\i,0);}

			\qedge{(y)}{(v3)}{(v2)}{5pt}{1pt}{blue!70!black}{blue!50!white,opacity=0.2};
			\qedge{(x)}{(v1)}{(v0)}{5pt}{1pt}{red!70!black}{red!50!white,opacity=0.2};
			\qedge{(y)}{(v4)}{(v3)}{5pt}{1pt}{blue!70!black}{blue!50!white,opacity=0.2};
			\qedge{(x)}{(v2)}{(v1)}{5pt}{1pt}{red!70!black}{red!50!white,opacity=0.2};
			
			\draw [thick]
			(v0) edge [red!70!black] (v1)
			(v1) edge [red!70!black] (v2)
			(v2) edge [blue!70!black] (v3)
			(v3) edge [blue!70!black](v4);

			\foreach \i in {0,...,4}{
				\fill (v\i) circle (2pt);
			}	
			\fill [red!70!black] (x) circle (2pt);
			\fill [blue!70!black] (y) circle (2pt);	
			\end{tikzpicture}\\
			
			\begin{tikzpicture}	[scale = .8]
			\node at (-.3,1.2){{$(b)$}};
			
			%		\coordinate (x) at (.5,1);
			%		\coordinate (y) at (1.3,1);
			
			\coordinate (v0) at (-.1,0);
			\coordinate (v1) at (.5,0);
			
			\qedge{(x)}{(y)}{(v1)}{4.5pt}{1pt}{blue!0!black}{black!20!white,opacity=0.2};
			\qedge{(v2)}{(v3)}{(v4)}{4.5pt}{1pt}{blue!0!black}{black!20!white,opacity=0.2};
			\qedge{(y)}{(v3)}{(v2)}{7pt}{1pt}{blue!70!black}{blue!50!white,opacity=0.2};
			\qedge{(x)}{(v1)}{(v0)}{7pt}{1pt}{red!70!black}{red!50!white,opacity=0.2};
			
			\node [red!70!black] at ($(.15,.15)$) {\tiny $f$};
			\node [blue!70!black] at ($(1.6,.3)$) {\tiny $g$};
			\node  at (2.7,-.1) {\tiny $h$};
			%\node  at (1.05,.7) {\tiny $h'$};
			\node at ($(x)+(0,-.15)$)  {\tiny $x$};
			\node at ($(y)+(-.15,-.15)$) {\tiny $y$};
			\node at ($(v1)+(.05,.14)$)  {\tiny $z$};

			\draw [thick]
			(v0) edge [red!70!black] (v1)
			(v2) edge [blue!70!black](v3);
			
			\foreach \i in {0,...,4}{
				\fill (v\i) circle (2pt);
			}	
			\fill [red!70!black] (x) circle (2pt);
			\fill [blue!70!black] (y) circle (2pt);
			\end{tikzpicture}\\

		\begin{tikzpicture}[scale = .5]
		
		\coordinate (n) at (-1.8,1.2);
		\node at (n) {\ref{it:typa}};
		
					\foreach \i in {0,...,4}{
						\coordinate (v\i) at (90+\i*72:1.3cm);
					}
					\node [blue!70!black] at ($(v1)+(.6,.13)$) {\tiny $e$};
					
					\foreach \i/\j in {0/1,0/4,1/2,1/3,2/4,3/4}{
						\draw ($(v\i)+(0,0)$) edge [color = red!90!black, thick] ($(v\j)+(0,0)$);
					}
					
					\draw ($(v1)+(0,0)$) edge [color = blue!90!black, thick] ($(v4)+(0,0)$);
					
					\draw ($(v0)+(-.065,-.065)$) edge [color = red!90!black, thick] ($(v2)+(-.065,-.065)$);
					\draw ($(v0)+(.065,.065)$) edge [color = blue!80!black, thick] ($(v2)+(.065,.065)$);
					
					\draw ($(v0)+(-.045,-.045)$) edge [color = red!90!black, thick] ($(v3)+(-.045,-.045)$);
					\draw ($(v0)+(.035,.035)$) edge [color = blue!90!black, thick] ($(v3)+(.035,.035)$);
					
					\draw ($(v2)+(0,-.045)$) edge [color = red!90!black, thick] ($(v3)+(0,-.045)$);
					\draw ($(v2)+(0,.045)$) edge [color = blue!90!black, thick] ($(v3)+(0,.045)$);	
					
	%				\node [anchor = south west] at (v0) {$a$};
	%				\node [anchor = north east] at (v2) {$b$};
	%				\node [anchor = north west] at (v3) {$c$};
	%				\node [anchor = south west] at (v4) {$y$};		
	%				\node [anchor = south east] at (v1) {$x$};
					
					\foreach \i in {0,...,4}{
						\fill (v\i) circle (3pt);
					}
			
		\end{tikzpicture}\\

		\begin{tikzpicture}[scale = .5]
			\node at (n) {\ref{it:typb}};

		\foreach \i in {0,...,4}{
			\coordinate (v\i) at (90+\i*72:1.3cm);
		}
		
		\foreach \i/\j in {0/1,0/4,1/2,1/3,2/4,3/4,0/2,0/3}{
			\draw ($(v\i)+(0,0)$) edge [color = red!90!black, thick] ($(v\j)+(0,0)$);
		}
		
	%	\draw ($(v1)+(0,0)$) edge [color = red!90!black, line width=2pt] ($(v4)+(0,0)$);
		
		\draw ($(v1)+(0,-.045)$) edge [color = red!90!black, thick] ($(v4)+(0,-.045)$);
		\draw ($(v1)+(0,.045)$) edge [color = blue!80!black, thick] ($(v4)+(0,.045)$);
		
		\draw ($(v2)+(0,-.045)$) edge [color = red!90!black, thick] ($(v3)+(0,-.045)$);
		\draw ($(v2)+(0,.045)$) edge [color = blue!90!black, thick] ($(v3)+(0,.045)$);
		
%		\draw ($(v2)+(0,-.04)$) edge [color = red!90!black, line width=2pt] ($(v3)+(0,-.04)$);
%		\draw ($(v2)+(0,.04)$) edge [color = blue!90!black, line width=2.3pt] ($(v3)+(0,.04)$);	

		\foreach \i in {0,...,4}{
			\fill (v\i) circle (3pt);
		}

		%		\node [anchor = south west] at (v0) {$a$};
		%		\node [anchor = north east] at (v2) {$b$};
		%		\node [anchor = north west] at (v3) {$c$};
	%	\node [anchor = south west] at (v4) {$y$};		
	%	\node [anchor = south east] at (v1) {$x$};
		
		\end{tikzpicture}\\

	\begin{tikzpicture}[scale = .5]
	
	\node at (n) {\ref{it:typc}};
	
%	\foreach \i in {0,...,4}{
%		\coordinate (v\i) at (90+\i*72:1.5cm);
%	}

		\coordinate (v0) at (0,1.5);
		\coordinate (v1) at (-.8,.8);
		\coordinate (v2) at (-.8,-.8);
		\coordinate (v3) at (.8,-.8);
		\coordinate (v4) at (.8,.8);

	\foreach \i/\j in {1/3,2/4}{
		\draw ($(v\i)+(0,.065)$) edge [color = red!90!black, thick] ($(v\j)+(0,.065)$);
		\draw ($(v\i)-(0,.065)$) edge [color = blue!90!black, thick] ($(v\j)-(0,.065)$);
	}

	\foreach \i/\j in {1/4,3/2}{
		\draw ($(v\i)+(0,.045)$) edge [color = red!90!black, thick] ($(v\j)+(0,.045)$);
		\draw ($(v\i)-(0,.045)$) edge [color = blue!90!black, thick] ($(v\j)-(0,.045)$);
	}
	
\draw ($(v1)+(.045,0)$) edge [color = red!90!black, thick] ($(v2)+(.045,0)$);
\draw ($(v3)+(.045,0)$) edge [color = red!90!black, thick] ($(v4)+(.045,0)$);

\draw ($(v1)-(.045,0)$) edge [color = blue!90!black, thick] ($(v2)-(.045,0)$);
\draw ($(v3)-(.045,0)$) edge [color = blue!90!black, thick] ($(v4)-(.045,0)$);

	%			\node [anchor = south west] at (v0) {$a$};
	%			\node [anchor = north east] at (v2) {$b$};
	%			\node [anchor = north west] at (v3) {$c$};
	%			\node [anchor = south west] at (v4) {$y$};		
	%			\node [anchor = south east] at (v1) {$x$};
	
	\foreach \i in {0,...,4}{
		\fill (v\i) circle (3pt);
	}
	
	\end{tikzpicture}

			\begin{tikzpicture}[scale = .5]
				\node at (n) {\ref{it:typd}};
				
			\coordinate (v0) at (0,0);
			\coordinate (v1) at (-1,1);
			\coordinate (v2) at (-1,-1);
			\coordinate (v3) at (1,-1);
			\coordinate (v4) at (1,1);	
			
			\foreach \i in {1,...,4}{
				\draw ($(v\i)+(0,.06)$) edge [color = blue!90!black, thick] ($(v0)+(0,.06)$);
				\draw ($(v\i)-(0,.06)$) edge [color = red!90!black, thick] ($(v0)-(0,.06)$);
			}
			
			\draw (v1) edge [color = blue!90!black, thick] (v2);
			\draw (v3) edge [color = blue!90!black, thick] (v4);
			
			\draw (v1) edge [color = red!90!black, thick] (v4);
			\draw (v2) edge [color = red!90!black, thick] (v3);

			\foreach \i in {0,...,4}{
				\fill (v\i) circle (3pt);
			}
			\end{tikzpicture}
				
	\end{multicols}
	
	\vspace{-1em}

	\caption{The illustration to the proof of Lemma \ref{l:n7}.}	\label{fig:l41}
	\vspace{-1em}
\end{figure}

		Having established Claim \ref{cl:z} we proceed with the proof of Lemma \ref{l:n7}. To this end look at the link graphs $\rr$ and $\bb$, and observe that the $\p$-freeness of $\cH$ entails $\RB\nsubseteq \rr\cup\bb$ (see Figure \ref{fig:l41}$(a)$).
		
		First assume $\deg_\cH(x,y)=\Delta_2(\cH)=3$. This implies that in each graph, $\rr$ and $\bb$, the vertices $z_1,z_2,z_3\in N_\cH(x,y)$ have degree at most 2, while the remaining two vertices of $Z$ have degree at most 3. Hence, by Lemma \ref{11}, $|\rr|+|\bb|\le 10$. On the other hand, Claim \ref{cl:z}\ref{it:z3} together with \eqref{eq:20} tell us that $|\rr|+|\bb| \ge 11$, a contradiction.
		
		Preparing for the remaining two cases, we make the following observation due to the $\p$-freeness of $\cH$.
		Suppose there is a 3-edge $h\in \cH[Z]$ and two 2-edges, $\crr{f}\in \rr$ and $\cbb{g}\in \bb$ such that $\crr{f}\cap h=\emptyset$, $\crr{f}\cap N_\cH(x,y)\neq\emptyset$ and $\cbb{g}\subset h$. Then, for any vertex $z\in \crr{f}\cap N_\cH(x,y)$, 3-edges $\crr{f}x,zxy,y\cbb{g},h$ form a minimal 4-path in $\cH$, a contradiction (see Figure \ref{fig:l41}$(b)$). Note further that in the above argument one can exchange the graphs $\rr$ and $\bb$.
			
		Next, let $\deg_{\cH}(x,y) = 5$. Then~\eqref{eq:20} and Claim~\ref{cl:z}\ref{it:z1} entails $|\rr|+|\bb|\ge 13$. Consequently, in view of Lemma \ref{l:five}, $|\cH[Z]|=2$ and $|\rr|+|\bb|=13$, and thereby there is a 2-edge $e\in {Z\choose 2}$ such that, $\rr= K^{(2)}_5[Z]-e$, $\bb=K^{(2)}_3[Z\setminus e]\cup e$, because all the other graphs described in \ref{it:typa}-\ref{it:typd} satisfy $|\rr|+|\bb|\le 12$ (see Figure \ref{fig:l41}$(A)$).
		Now writing $Z=\{a,b,c\}\cup e$, we let $h=ea$, $\crr{f}=bc$ and $\cbb{g}=e$, which satisfy the assumptions in the previous paragraph and thus yield a contradiction.
		
		Finally, let $\deg_{\cH}(x,y)= 4$, and write $N\coloneq N_\cH(x,y)$.
		In view of \eqref{eq:20} combined with Claim~\ref{cl:z}\ref{it:z2}, $|\cH[Z]|\le 4$ and $|\rr|+|\bb|\ge 12$. Then again, Lemma \ref{l:five} tells us that one of \ref{it:typa}-\ref{it:typd} holds. Moreover the condition $\Delta_2(\cH)=\deg_{\cH}(x,y)=4$ entails that only the unique vertex of $Z\setminus N$ can have degree 4 in $\rr$, and thus the cases $\rr=K^{(2)}_5[Z]$ and $\rr=K^{(2)}_5[Z]-e$ are excluded. 	
		Note that all the remaining 2-graphs $\rr\cup \bb$ with $|\rr|+|\bb|\ge 12$, described in Lemma \ref{l:five}, namely $\rr = K^{(2)}_5[V]-\{e,e'\}$, $\bb=T\cup e$, \ref{it:typc}, and \ref{it:typd}, satisfy $|\rr|+|\bb|=12$, implying that $|\cH[Z]|=4$ and thus $\cH[Z] = K^{(3)}_4[N]$. Moreover, they have the property that every 3-vertex set $h\subset Z$ contains an edge of both 2-graphs $\rr$ and $\bb$ (see Figure \ref{fig:l41}$(A), (C),(D)$), and, as $\deg_{\cH}(x,y)= 4$, every edge of $\rr\cup \bb$ intersects $N$. Therefore, if $\rr\cup \bb\nsubseteq K_4[N]$, one can take $(f,g)\in (\rr,\bb)\cup (\bb,\rr)$ with $f\nsubseteq N$ and $h=Z\setminus f\in \cH[Z]$, $g\subseteq h$, yielding a contradiction with the $\p$-freeness of $\cH$ (see Figure \ref{fig:l41}$(b)$). Thus, we conclude that $\rr\cup \bb\subseteq K^{(2)}_4[N]$ and $|\rr|+|\bb|=12$ implies that $\rr=\bb=K^{(2)}_4[N]$. Altogether $\cH=K^{(3)}_6\cup K_1$, as required.

%%%%%%%%%%%%%%%%%%%%%%%%%%%%%%%%%%%%%%%%%%%%%%%%%%%%%%%%%%%%%%%%%%%%%%
%																																		%
%																												   						%
%  									Structure of $\cP$-free 3-graphs					   								%
%																																		%
%																																		%
%%%%%%%%%%%%%%%%%%%%%%%%%%%%%%%%%%%%%%%%%%%%%%%%%%%%%%%%%%%%%%%%%%%%%%

\section{Proofs of Lemmas \ref{l:nu2} and \ref{l:M3}}

\subsection{Structure of $\p$-free 3-graphs}\label{str}

In this subsection we gather some basic information about the structure of connected $\p$-free 3-graphs. We begin by showing that such 3-graphs may contain at most three disjoint edges.  To this end, let us make the following observations.

\begin{fact}\label{clm:1} For every connected $\p$-free 3-graph $\cH$ the following holds.
	\begin{enumerate}[label=\rmlabel]
		\item\label{it:51} If $e_1,e_2\in \cH$ are disjoint, then there exists an edge $f\in \cH$ intersecting both $e_1$ and $e_2$.
		\item\label{it:52} If $e_1,e_2\in \cH$ are disjoint and $f,h\in \cH$ are such that $f\cap e_1\neq\emptyset$, $f\cap e_2\neq\emptyset$,  $h\cap e_1=\emptyset$, and
		$h\cap e_2\neq\emptyset$, then $f\cap h\neq\emptyset$.
		\item\label{it:53} If $e_1,e_2,e_3,f,h \in \cH$ are such that $e_1,e_2,e_3$ are pairwise disjoint, $f\cap e_1\neq\emptyset$, $f\cap e_2\neq\emptyset$, $f\cap e_3=\emptyset$,
		$h\cap e_2\neq\emptyset$, and $h\cap e_3\neq\emptyset$,
		then $h\cap e_1\neq\emptyset$.
		\item\label{it:54}  If $e_1,e_2,e_3 \in \cH$ are pairwise disjoint, then there exists an edge intersecting all the three edges $e_1$, $e_2$, and $e_3$.
	\end{enumerate}
\end{fact}
\begin{proof}
	To prove \ref{it:51} observe that in a connected 3-graph every pair of disjoint edges, $e_1$ and $e_2$, is connected by a minimal path $P$. If additionally there is no edge in $\cH$ intersecting both $e_1$ and $e_2$, then $P$ consists of at least four edges.
	
	For the proof of \ref{it:52} note that otherwise $e_1fe_2h$ would form a minimal 4-path in $\cH$. Next, to show \ref{it:53} observe that $h\cap e_1=\emptyset$ together with \ref{it:52} entails $f\cap h\neq \emptyset$ and, since $f\cap e_3= \emptyset$, $e_1fhe_3$ is a minimal 4-path in $\cH$, a contradiction.

	Finally, to deduce \ref{it:54} we apply \ref{it:51} twice getting two (not necessary different) edges $f,h\in \cH$, such that $f$ intersects $e_1$ and $e_2$, while $h$ intersects $e_2$ and $e_3$. If, additionally, $f\cap e_3\neq\emptyset$, we  are done. Otherwise \ref{it:53} yields $h\cap e_1\neq\emptyset$, which conclude the proof.
\end{proof}

Now we are ready to prove the promised, crucial fact.

\begin{lemma}\label{lem:1}
	If $\cH$ is a connected $\p$-free 3-graph, then $\nu(\cH)\le3$.
\end{lemma}
\begin{proof}
	Suppose that $\nu(\cH)\ge 4$ and fix four disjoint edges $e_1,e_2,e_3,e_4\in \cH$. Double application of Fact~\ref{clm:1}\ref{it:54} entails the existence of two edges, $f,h\in \cH$, such that $f$ intersects $e_1,e_2,e_3$, while $h$ intersects $e_2,e_3,e_4$. Clearly $h\cap e_1=\emptyset$ and thus, due to Fact \ref{clm:1}\ref{it:52}, $f\cap h \neq \emptyset$. But then $e_1fhe_4$ is a minimal $4$-path in $\cH$, a contradiction.
\end{proof}

As a preparation towards the proofs of Lemmas \ref{l:nu2} and \ref{l:M3}, we now make an attempt to characterize all connected $\p$-free 3-graphs with at least two disjoint edges. As an exception, in this section, in order to distinguish between ordinary graphs (2-graphs) and $3$-graphs, we will  use notation $\mathcal F$, with subscripts, for single $3$-graphs rather than families of $3$-graphs. (But we keep $H$ unchanged, as it clearly associates itself with hypergraphs.)

To this end, recall that a hypergraph $\cF$ is \emph{intersecting} if $f\cap f'\neq\emptyset$  for every $f,f'\in \cF$. Similarly, a pair $(\cF, \cF')$ of hypergraphs is called \emph{cross-intersecting}, if for all $f\in \cF$, $f'\in \cF'$ we have $f\cap f' \neq \emptyset$. It turns out that every connected $\p$-free 3-graph $\cH$ with $\nu(\cH)\in\{2,3\}$, can be described as follows.

\begin{lemma}\label{l:split}
	Every $\p$-free connected 3-graph $\cH$ with $\nu(\cH)=2$ on the set of vertices $V$, can be partitioned into three edge-disjoint 3-graphs $\cH = \cF_1\dcup \cF_2\dcup \ff$, such that
	\begin{enumerate}[label=\rmlabel]
		\item\label{it:st1} $V[\cF_1]\cap V[\cF_2]=\emptyset$,
		\item\label{it:st2} the 3-graphs $\cF_1$ and $\cF_2$ are non-empty intersecting families,
		\item\label{it:st3} $\ff\neq\emptyset$,
		\item\label{it:st4} the pair $(\cF_1\dcup \cF_2, \ff)$ is cross-intersecting.
	\end{enumerate}
\end{lemma}
\begin{lemma}\label{l:split3}
	Every $\p$-free connected 3-graph $\cH$ with $\nu(\cH) =3$ on the set of vertices $V$, can be partitioned into five edge-disjoint 3-graphs $\cH = \cF_1\dcup \cF_2\dcup \cF_3\dcup \ff\dcup \fff$, such that
	\begin{enumerate}[label=\rmlabel]
		\item the sets $V[\cF_1]$, $V[\cF_2]$, and $V[\cF_3]$ are pairwise disjoint, and $V[\ff]\cap V[\cF_3]=\emptyset$,
		\item the 3-graphs $\cF_1$, $\cF_2$, and $\cF_3$ are non-empty intersecting families,
		\item $\fff\neq\emptyset$,
		\item the pairs $(\cF_1\dcup \cF_2\dcup \cF_3\dcup \ff,  \fff)$ and $(\cF_1\dcup \cF_2,\ff)$ are cross-intersecting.
	\end{enumerate}
\end{lemma}
\begin{proof}[Proof of Lemmas \ref{l:split} and \ref{l:split3}]
	Let $\cH$ be a given $\p$-free connected 3-graph on $V$, and let $k=\nu(\cH)$, $k=2,3$. Fix a largest matching $M_k=\{e_1,\dots,e_k\}\subset \cH$. Now, for each $I\subseteq [k]$, $\cF_I$ is defined to be the set of all edges of $\cH$ that intersect every $e_i$, $i\in I$, and none of $e_j$, $j\in [k]\setminus I$. Clearly $e_i \in \cF_{\{i\}}$, $\cF_\emptyset = \emptyset$ and
		\begin{equation*}\label{eq:FI}
			\cH = \bigdcup_{I\subseteq [k]}\cF_I.
		\end{equation*}
	For simplicity of notation, we write $\fff$ instead of $\cF_{\{1,2,3\}}$, $\ff$ instead of $\cF_{\{1,2\}}$, etc.
	
	First note that in view of Fact \ref{clm:1}\ref{it:53}, for $k=3$ at most one of $\cF_{12}, \cF_{13}, \cF_{23}$, say $\ff$, is nonempty.	Now, if for some vertex $v\in V\setminus (\bigcup_{i\in [k]}e_i)$ there are two edges $f,h\in \cH$ such that $v\in f\cap h$, $f\in \cF_i$ and $h\in \cF_j\cup  \cF_{jk}$, $\{i,j,k\}=\{1,2,3\}$, then $e_ifhe_j$ is a minimal 4-path in $\cH$. But $\cH$ is $\p$-free and thus the sets $V[\cF_1]$, $V[\cF_2]$, and $V[\cF_3]$ are pairwise disjoint, and $V[\ff]\cap V[\cF_3]=\emptyset$, establishing \ref{it:st1}. Consequently, as $\nu(\cH)=k$ and $e_i\in \cF_i$ for each $i\in [k]$, every $\cF_i$ is a non-empty intersecting family, and thus \ref{it:st2} follows.
	
	Further, $\ff\neq\emptyset$ and $\fff\neq\emptyset$ result from Fact \ref{clm:1}\ref{it:51} and \ref{it:54}, respectively. Finally, Fact \ref{clm:1}\ref{it:52} tells us that the pairs $(\cF_1\dcup \cF_2,\ff)$ and $(\cF_1\dcup \cF_2\dcup \cF_3\dcup \ff,  \fff)$ (for $k=3$), are cross-intersecting.
\end{proof}

%%%%%%%%%%%%%%%%%%%%%%%%%%%%%%%%%%%%%%%%%%%%%%%%%%%%%%%%%%%%%%%%%%%%%%
%																																		%
%																												   						%
%  									nu   =   2																   						%
%																																		%
%																																		%
%%%%%%%%%%%%%%%%%%%%%%%%%%%%%%%%%%%%%%%%%%%%%%%%%%%%%%%%%%%%%%%%%%%%%%

\subsection{Proof of Lemma \ref{l:nu2}}

Let $\cH$ be a $\{\p,\c, \m\}$-free connected 3-graph on the set of vertices $V$, $|V|=n\ge 8$, and let $\cH\nsubseteq S_n^{+1}$, $\cH\nsubseteq SP_n$, and $\cH\nsubseteq CB_n$. We are to show that
	\begin{equation}\label{eq:lnu2}
		|\cH| \le \max\left\{4n-11,{n-4\choose 2}+10\right\}.
	\end{equation}

To prove this observe that because $\cH\nsubseteq S_n$, if $\nu(\cH)=1$, then in view of Theorem \ref{th:hm}, $|\cH|\le 3n-8 < 4n-11$, and we are done. Therefore, as $\cH$ is $\m$-free, we may assume $\nu(\cH)=2$ and take a partition
	\[
		\cH= \cF_1\dcup \cF_2\dcup \cF_{12}
	\]
guaranteed by Lemma \ref{l:split}. Recall that both $\cF_1$ and $\cF_2$ are non-empty intersecting families, $\ff\neq\emptyset$, and the pair $(\cF_1\dcup \cF_2, \cF_{12})$ is cross-intersecting.
For $i=1,2$, let $S_i\subseteq V(\cF_i)$ be the set of vertices $s$ that lie in all edges of $\cF_i$.
Clearly, $s_i\coloneq |S_i|$ satisfies $0\le s_i\le 3$ and without loss of generality we may assume $0\le s_2\le s_1\le 3$.

Set $V_1=V[\cF_1]$, $V_2=V\setminus V_1$, and note that $V[\cF_2]\subseteq V_2$.
A pair $p$ of vertices in $V_i$ is called a 2-cover of $\cF_i$ if it intersects every edge of $\cF_i$, i.e., $p\cap f\neq\emptyset$ holds for all $f\in \cF_i$.
Denote by $P_i\subseteq \binom{V_i}2$ the collection of all 2-covers of $\cF_i$.
Now \eqref{eq:lnu2}, and thereby Lemma \ref{l:nu2}, is a straightforward consequence of the following claim.

\begin{claim}\label{cl:ss}
	\begin{enumerate}[label=\rmlabel]
		\item\label{it:s1} If $s_1\ge s_2\ge 2$ then $|\cH|\le 4n-11$.	
		\item\label{it:s2} 	If $s_1=3$, $s_2 = 1$, $\cH\nsubseteq S_n^{+1}$, and $\cH\nsubseteq SP_n$, then $|\cH|\le 4n-11$.	
		\item\label{it:s3} 	If $s_1=2$, $s_2 = 1$, and $\cH\nsubseteq CB_n$, then $|\cH|\le \max\left\{4n-11,{n-4\choose 2}+10\right\}$.		
		\item\label{it:s4} If $s_1 =  s_2 = 1$ then $|\cH|\le \max\left\{4n-11,{n-4\choose 2}+10\right\}$.					
		\item\label{it:s5} If $s_2 = 0$ then $|\cH|\le \max\left\{4n-11,{n-4\choose 2}+10\right\}$.		
	\end{enumerate}	
\end{claim}
\begin{proof}
	Let us start with the proof of \ref{it:s1}, that is $s_1\ge s_2\ge 2$. To this end, for each $i=1,2$ pick an edge $e_i\in \cF_i$, and set $W=V\setminus (e_1\cup e_2)$, $|W|=n-6$. Then for every $z\in W$,
		\begin{equation}\label{eq:f1f2}
			\deg_{\cF_1\cup \cF_2}(z)\le 1
		\end{equation}
	follows from $s_1,s_2\ge 2$ and $V[\cF_1]\cap V[\cF_2]=\emptyset$.
	
	Now, let $u,w\in W$ in the case $n=8$, and $u,w,v\in W$ otherwise, be vertices with the largest degrees in $\ff$, such that
		\[
			\deg_{\cF_{12}}(u)\ge \deg_{\cF_{12}}(w)\ge \deg_{\cF_{12}}(v).
		\]
	We may assume that $\deg_{\cF_{12}}(w) \ge 2$. Otherwise, as $\hat \cH= \cH[e_1\cup e_2\cup \{u\}]$ has no isolated vertices, Lemma \ref{l:n7} tells us $|\hat \cH|\le 19$, and by \eqref{eq:f1f2} for $n\ge 8$ we have
		\[
			|\cH|= |\hat \cH| + \sum_{z\in W\setminus \{u\}}(\deg_{\cF_1\cup \cF_2}(z) +\deg_{\cF_{12}}(z))
				\le
			19 + 2(n-7) \le 4n-11.
		\]

	We contend
		\begin{equation}\label{eq:310}
			\deg_\cH(u) + \deg_\cH(w)\le 10 \qand \deg_{\cF_{12}}(v)\le 3,
		\end{equation}
	which ends the proof. Indeed, observe that the absence of $\c$ in $\cH$ entails $|\cH[e_1\cup e_2]|\le 11$, because, the set of edges of $K^{(3)}_6$ can be partitioned into 10 pairs of disjoint edges, and any two of these pairs form $\c$. Therefore \eqref{eq:310} combined with \eqref{eq:f1f2} tells us
		\[
			|\cH|= |\cH[e_1\cup e_2]|+\deg_\cH(u) + \deg_\cH(w)+
			\sum_{z\in W\setminus \{u,w\}}(\deg_{\cF_1\cup \cF_2}(z)+ \deg_{\cF_{12}}(z))\le 4n-11.
		\]

	To show \eqref{eq:310}, instead of looking at the degrees of $u$, $w$, and $v$ it is more convenient for us to look at their link graphs in $\cF_{12}$,
		\[
			\rr=L_{\cF_{12}}(u),\quad \bb=L_{\cF_{12}}(w), \qand \gg=L_{\cF_{12}}(v).
		\]
	Because every edge of $\ff$ intersects both $e_1$ and $e_2$, actually $\rr, \bb, \gg\subseteq K^{(2)}_{3,3}[e_1\dcup e_2]$.
	
	We first note that the $\{\c, \p\}$-freeness of $\cH$ entails some forbidden configurations of edges of $\rr$, $\bb$ and $\gg$. In particular, there are no two distinct vertices $x,y\in e_i$, $i=1,2$, such that $\deg_{\rr}(x),\deg_{\bb}(y)\ge 2$ (see Figure \ref{fig:f7}$(a),(b)$, similar with $\gg$ in place of $\rr$ or $\bb$). This immediately entails that if for some $i=1,2$, there are two distinct vertices $x,y\in e_i$ with
		\begin{equation}\label{eq:2deg2}
			\deg_{\rr}(x)\ge 2 \textrm{ and } \deg_{\rr}(y)\ge 2
			\quad\textrm{ then for all $z\in e_i$}\quad
			\deg_{\bb}(z)\le 1.
		\end{equation}
	In particular, whenever $|\rr|\ge 6$, then $\bb$ is a matching and thus $|\bb|\le 3$. Moreover, as there are no three disjoint edges in $K_{3,3}^{(2)}[e_1\dcup e_2]$, such that at least two of them are in $\rr$ and at least two of them are in $\bb$ (see Figure \ref{fig:f7}$(c),(d)$), because $|\bb| = \deg_{\cF_{12}}(w)\ge 2$, we have $|\rr|\le 7$ and if $|\rr|=7$ then $|\bb|=2$. Indeed, otherwise $|\rr|=7$ and $|\bb|=3$ yields either two disjoint edges in $\rr\cap \bb$ (see Figure \ref{fig:f7}$(c)$), or three disjoint edges, two in $\bb$ and one in $\rr$ and one edge in $\rr$ connecting the $\rr$-edge with the $\bb$-edge, entailing the existence of a minimal 4-path in $\cH$ (see Figure \ref{fig:f7}$(e)$).
	
\begin{figure}[h!]
		\begin{multicols}{10}

			\begin{tikzpicture}[scale = 1]	
				\coordinate (w1) at (0,.5);
				\coordinate (w2) at (0,0);
				\coordinate (v1) at (.7,.5);
				\coordinate (v2) at (.7,0);

				\node at (0,.8) {\tiny $(a)$};
				
				\draw (w1) edge [color = red!70!black,very thick] (v2);
				\draw (w2) edge [color = blue!70!black,very thick] (v2);
				\draw (w1) edge [color = red!70!black,very thick] (v1);
				\draw (w2) edge [color = blue!70!black,very thick] (v1);
				
				\foreach \i in {w1, w2, v1, v2}
				\fill (\i) circle (1.3pt);				
			\end{tikzpicture}

			\begin{tikzpicture}[scale = 1]				
				\coordinate (w) at (-.6, 0);
				\coordinate (v) at (.6, 0);
				
				\coordinate (x1) at (0, .45);
				\coordinate (x2) at (0, -.45);
				\coordinate (x3) at (.4, 0);
				\coordinate (x4) at (-.4, 0);
								
				\qedge{(x1)}{(v)}{(x3)}{4pt}{1pt}{blue!70!black}{blue!50!white,opacity=0.2};
				\qedge{(x3)}{(v)}{(x2)}{4pt}{1pt}{blue!70!black}{blue!50!white,opacity=0.2};
				\qedge{(x1)}{(x4)}{(w)}{4pt}{1pt}{red!70!black}{red!50!white,opacity=0.2};
				\qedge{(w)}{(x4)}{(x2)}{4pt}{1pt}{red!70!black}{red!50!white,opacity=0.2};
								
				\draw (x1) edge [color = blue!70!black] (x3) ;
				\draw (x1) edge [color = red!70!black] (x4);
				\draw (x2) edge [color = blue!70!black] (x3) ;
				\draw (x2) edge [color = red!70!black] (x4);
								
				\foreach \i in {x1, x2, x3, x4}
				\fill (\i) circle (1.3pt);
				
				\fill (w) [color = red!80!black] circle (1.3pt);
				\fill (v) [color = blue!80!black] circle (1.3pt);				
			\end{tikzpicture}

			\begin{tikzpicture}[scale = 1]		
				\node at (-.3,.4) {\tiny $(b)$};
				
				\coordinate (w1) at (0,.5);
				\coordinate (w2) at (0,0);
				\coordinate (w3) at (0,-.5);
				\coordinate (v1) at (.7,.5);
				\coordinate (v2) at (.7,0);
				\coordinate (v3) at (.7,-.5);
			
				\draw (w1) edge [color = red!70!black, very thick] (v1);			
				\draw (w2) edge [color = blue!70!black, very thick] (v2);
				\draw (w1) edge [color = red!70!black, very thick] (v2);
				\draw (w2) edge [color = blue!70!black, very thick] (v3);
			
				\foreach \i in {w1, w2, v1, v2, v3}
				\fill (\i) circle (1.3pt);			
			\end{tikzpicture}
					
			\begin{tikzpicture}[scale = 1]			
				\coordinate (w) at (.35, .4);
				\coordinate (v) at (.35, -.1);
			
				\coordinate (v1) at (.8,.5);
				\coordinate (v2) at (.8,0);
				\coordinate (v3) at (.8,-.5);
			
				\qedge{(w1)}{(v1)}{(w)}{3pt}{1pt}{red!70!black}{red!50!white,opacity=0.2};
				\qedge{(w1)}{(w)}{(v2)}{3pt}{1pt}{red!70!black}{red!50!white,opacity=0.2};
				\qedge{(w2)}{(v2)}{(v)}{3pt}{1pt}{blue!70!black}{blue!50!white,opacity=0.2};
				\qedge{(w2)}{(v)}{(v3)}{3pt}{1pt}{blue!70!black}{blue!50!white,opacity=0.2};
			
				\draw (w1) edge [color = red!70!black] (v1) ;
				\draw (w1) edge [color = red!70!black] (v2);
				\draw (w2) edge [color = blue!70!black] (v2) ;
				\draw (w2) edge [color = blue!70!black] (v3);
							
				\foreach \i in {w1, w2, v1, v2, v3}
				\fill (\i) circle (1.3pt);
			
				\fill (w) [color = red!80!black] circle (1.3pt);
				\fill (v) [color = blue!80!black] circle (1.3pt);		
			\end{tikzpicture}
			
				\begin{tikzpicture}[scale = 1]				
				\node at (0,.8) {\tiny $(c)$};
				
				\coordinate (w1) at (0,.5);
				\coordinate (w2) at (0,0);
				\coordinate (v1) at (.8,.5);
				\coordinate (v2) at (.8,0);

				\draw ($(w1)+(0,.027)$) edge [color = red!70!black, very thick] ($(v1)+(0,.027)$) ;
				\draw ($(w1)+(0,-.027)$) edge [color = blue!70!black, very thick] ($(v1)+(0,-.027)$);
				\draw ($(w2)+(0,.027)$) edge [color = red!70!black, very thick] ($(v2)+(0,.027)$) ;
				\draw ($(w2)+(0,-.027)$) edge [color = blue!70!black, very thick] ($(v2)+(0,-.027)$);
				
				\foreach \i in {w1, w2, v1, v2}
				\fill (\i) circle (1.3pt);
				\end{tikzpicture}	
				
				\begin{tikzpicture}[scale = 1]				
				\coordinate (w) at (.6, 0);
				\coordinate (v) at (-.6, 0);
				
				\coordinate (x1) at (-.17, .45);
				\coordinate (x2) at (.17, .45);
				\coordinate (x3) at (-.17, -.45);
				\coordinate (x4) at (.17, -.45);
				
				\qedge{(x1)}{(x2)}{(v)}{4pt}{1pt}{red!70!black}{red!50!white,opacity=0.2};
				\qedge{(v)}{(x4)}{(x3)}{4pt}{1pt}{red!70!black}{red!50!white,opacity=0.2};
				\qedge{(x1)}{(x2)}{(w)}{4pt}{1pt}{blue!70!black}{blue!50!white,opacity=0.2};
				\qedge{(w)}{(x4)}{(x3)}{4pt}{1pt}{blue!70!black}{blue!50!white,opacity=0.2};
				
				\draw ($(x1)+(0,-.015)$) edge [color = blue!70!black] ($(x2)+(0,-.015)$) ;
				\draw ($(x1)+(0,.015)$) edge [color = red!70!black] ($(x2)+(0,.015)$);
				\draw ($(x3)+(0,-.015)$) edge [color = blue!70!black] ($(x4)+(0,-.015)$) ;
				\draw ($(x3)+(0,.015)$) edge [color = red!70!black] ($(x4)+(0,.015)$);
				
				\foreach \i in {x1, x2, x3, x4}
				\fill (\i) circle (1.3pt);
				
				\fill (w) [color = blue!80!black] circle (1.3pt);
				\fill (v) [color = red!80!black] circle (1.3pt);				
				\end{tikzpicture}

		\begin{tikzpicture}[scale = 1]				
		\node at (-.3,.4) {\tiny $(d)$};
		
		\coordinate (w1) at (0,.5);
		\coordinate (w2) at (0,0);
		\coordinate (w3) at (0,-.5);
		\coordinate (v1) at (.7,.5);
		\coordinate (v2) at (.7,0);
		\coordinate (v3) at (.7,-.5);

	%	\draw ($(w1)+(0,.03)$) edge [color = blue!70!black, very thick] ($(v1)+(0,.03)$) ;
		\draw (w1) edge [color = red!70!black, very thick] (v1);
		\draw ($(w2)+(0,-.027)$) edge [color = blue!70!black, very thick] ($(v2)+(0,-.027)$) ;
		\draw ($(w2)+(0,.027)$) edge [color = red!70!black, very thick] ($(v2)+(0,.027)$);
		\draw (w3) edge [color = blue!70!black, very thick] (v3);
%		\draw ($(w3)+(0,-.03)$) edge [color = red!70!black, very thick] ($(v3)+(0,-.03)$);
		
		\foreach \i in {w1, w2, w3, v1, v2, v3}
		\fill (\i) circle (1.3pt);
		\end{tikzpicture}	
		
		\begin{tikzpicture}[scale = 1]				
		\coordinate (w) at (.35, .3);
		\coordinate (v) at (.35, -.3);

		\qedge{(w1)}{(v1)}{(w)}{3pt}{1pt}{red!70!black}{red!50!white,opacity=0.2};
		\qedge{(v2)}{(w2)}{(w)}{3pt}{1pt}{red!70!black}{red!50!white,opacity=0.2};
		\qedge{(w2)}{(v2)}{(v)}{3pt}{1pt}{blue!70!black}{blue!50!white,opacity=0.2};
		\qedge{(v3)}{(w3)}{(v)}{3pt}{1pt}{blue!70!black}{blue!50!white,opacity=0.2};

		% \draw ($(w1)+(0,.03)$) edge [color = blue!70!black] ($(v1)+(0,.03)$) ;
			\draw (w1) edge [color = red!70!black] (v1);
		\draw ($(w2)+(0,-.015)$) edge [color = blue!70!black] ($(v2)+(0,-.015)$) ;
		\draw ($(w2)+(0,.015)$) edge [color = red!70!black] ($(v2)+(0,.015)$);
				\draw (w3) edge [color = blue!70!black] (v3) ;
		%\draw ($(w3)+(0,-.03)$) edge [color = red!70!black, very thick] ($(v3)+(0,-.03)$);

		\foreach \i in {w1,w2,w3,v1,v2,v3}
		\fill (\i) circle (1.3pt);
		
		\fill (w) [color = red!80!black] circle (1.3pt);
		\fill (v) [color = blue!80!black] circle (1.3pt);				
		\end{tikzpicture}

				\begin{tikzpicture}[scale = 1]				
				\node at (-.3,.4) {\tiny $(e)$};
				
				\coordinate (w1) at (0,.5);
				\coordinate (w2) at (0,0);
				\coordinate (w3) at (0,-.5);
				\coordinate (v1) at (.7,.5);
				\coordinate (v2) at (.7,0);
				\coordinate (v3) at (.7,-.5);

				\draw (w1) edge [color = red!70!black, very thick] (v1);
				\draw (w2) edge [color = blue!70!black, very thick] (v2);
				\draw (w1) edge [color = red!70!black, very thick] (v2);
				\draw (w3) edge [color = blue!70!black, very thick] (v3);
			
				\foreach \i in {w1, w2, w3, v1, v2, v3}
				\fill (\i) circle (1.3pt);
				\end{tikzpicture}	
				
				\begin{tikzpicture}[scale = 1]				
				\coordinate (w) at (.36, .39);
				\coordinate (v) at (.35, -.25);

				\qedge{(w1)}{(v1)}{(w)}{2.5pt}{1pt}{red!70!black}{red!50!white,opacity=0.2};
				\qedge{(w)}{(v2)}{(w1)}{2.5pt}{1pt}{red!70!black}{red!50!white,opacity=0.2};
				\qedge{(w2)}{(v2)}{(v)}{3pt}{1pt}{blue!70!black}{blue!50!white,opacity=0.2};
				\qedge{(v3)}{(w3)}{(v)}{3pt}{1pt}{blue!70!black}{blue!50!white,opacity=0.2};

				\draw (w1) edge [color = red!70!black] (v1);
				\draw (w2) edge [color = blue!70!black] (v2);
				\draw (w1) edge [color = red!70!black] (v2);
				\draw (w3) edge [color = blue!70!black] (v3) ;

				\foreach \i in {w1,w2,w3,v1,v2,v3}
				\fill (\i) circle (1.3pt);
				
				\fill (w) [color = red!80!black] circle (1.3pt);
				\fill (v) [color = blue!80!black] circle (1.3pt);				
				\end{tikzpicture}

		\end{multicols}
		
			\vspace{-1em}
			
		\caption{Forbidden configurations of edges of $\textcolor{red!60!black}{R}$ and $\textcolor{blue!60!black}{B}$.}
		\label{fig:f7}
				
	\end{figure}
	
%	\vspace{-1em}		
	
	Further, repeated applications of \eqref{eq:2deg2} tells us that $|\rr|=|\bb|=5$ is possible only when
	$\rr=\bb=\tikz{\draw (0,0) -- (.2,0) (0,0) -- (.2,.13) (0,0) -- (.2,.26) -- (0,.26) (0,.13) -- (.2,.26);
		\foreach \i in {0,.2} \foreach \j in {0,.13,.26} \fill (\i,\j) circle (.6pt);}$.
	But then $\rr\cap \bb$ contains two disjoint edges, contradicting $\c$-freeness of $\cH$ (see Figure \ref{fig:f7}$(c)$). For the same reason $|\gg|\le 3$. Indeed, otherwise $|\rr|\ge |\bb|\ge |\gg|\ge 4$ and all of these three graphs have the same two vertices of degree larger than one. Thus
	$\rr,\bb,\gg\subseteq\tikz{\draw (0,0) -- (.2,0) (0,0) -- (.2,.13) (0,0) -- (.2,.26) -- (0,.26) (0,.13) -- (.2,.26);
		\foreach \i in {0,.2} \foreach \j in {0,.13,.26} \fill (\i,\j) circle (.6pt);}$ (each of them misses at most one edge)
	and hence the intersection of some two of them contains two disjoint edges, again arriving at a contradiction. Summarizing all these observations so far, we obtain
		\begin{equation}\label{eq:rb9}
			|\rr|+|\bb|\le 9 \qand \deg_{\cF_{12}}(v)=|\gg|\le 3.
		\end{equation}
	Therefore to establish \eqref{eq:310} it remains to show that $\deg_\cH(u) + \deg_\cH(w)\le 10$.
	
	To this end, assume for the sake of a contradiction that $\deg_\cH(u) + \deg_\cH(w)\ge 11$. Then \eqref{eq:rb9} combined with \eqref{eq:f1f2} tells us that
		\[
			\deg_{\cF_1\cup \cF_2}(u) = \deg_{\cF_1\cup \cF_2}(w) = 1
			\qand |\rr|+|\bb|= 9				
		\]
	Without loss of generality we may assume that the edge $f\in \cF_1\cup \cF_2$ with $w\in f$ belongs to $\cF_1$.
	Recalling that $s_1\ge 2$ we infer $|e_1\cap f| =2$. Now, as every edge of $\cF_{12}$ intersects each one of $e_1, e_2, f$, we actually have $\rr \subseteq K^{(2)}_{2,3}[(e_1\cap f)\dcup e_2]$ and thus $|\rr|\le 6$. Therefore, because $|\rr|+|\bb| = 9$ entails $|\rr|\ge 5$, for $\{x,y\}= e_1\cap f$ we have $\deg_{\rr}(x)\ge 2$ and $\deg_{\rr}(y)\ge 2$. Hence \eqref{eq:2deg2} tells us $|\bb|\le 3$ and if $|\rr|=6$, $|\bb|=3$, then $\rr\cap\bb$ contains two disjoint edges, a contradiction (see Figure \ref{fig:f7}$(c)$).
	
\medskip

	Before we move to the proof of \ref{it:s2}-\ref{it:s5} let us show a few simple facts. First note, that for $i=1,2$,
		\begin{equation}\label{eq:degv3}
			\deg_{P_i}(v)\le 3 \quad \textrm{ for all }\quad v \in V_i\setminus S_i.
		\end{equation}
	Indeed, because $v$ is not a 1-cover of $\cF_i$, there exists an edge $f\in \cF_i$ with $v\notin f$. On the other hand, all 2-covers in $P_i$ intersect $f$. Hence $N_{P_i}(v)\subseteq f$ and $\deg_{P_i}(v)\le |f|=3$ follows. Moreover, as for  every edge $h\in \cF_i$ we have $|f\setminus h|\le 2$, one can also deduce that if $v\in h$, then $|N_{P_i}(v)\setminus h|\le 2$. Thus, in view of \eqref{eq:degv3},
		\begin{equation}\label{eq:p7}
			|P_i|\le\begin{cases}
				7, &\textrm{ for } s_i=0,\\
				|V_i|+3, &\textrm{ for } s_i=1,\\
				2|V_i|-2, &\textrm{ for } s_i=2.
			\end{cases}
		\end{equation}

	To see this, take any edge $h\in \cF_i$ and consider neighborhoods in $P_i$ of vertices of $h$. Clearly, as every 2-cover in $P_i$ intersects $h$, we have $P_i\subseteq \left\{p\in {V\choose 2}\colon p\cap h\neq\emptyset \right\}$. For $s_i=0$ observe that if there is at most one 2-cover in $P_i$ entirely contained in $h$, then there are at most six 2-covers in $P_i$ that contain exactly one vertex with $h$.
	Similarly, when $h$ contains two 2-covers (they share a vertex), the number of 2-covers in $P_i$ that contain exactly one vertex with $h$ is at most five; and when $h$ contains three 2-covers, this number is at most three.
	For $s_i=1$ we let $h=\{s,v,w\}$, where $s$ is the unique 1-cover of $\cF_i$. Now, because $\{s,v\},\{s,w\}\in P_i$, $\deg_{P_i}(s) \le |V_i|-1$, $\deg_{P_i}(v)\le 3$, and $\deg_{P_i}(w)\le 3$, we actually have $|P_i|\le (|V_i|-1)+2+2$. For $s_i=2$ similar analysis implies $|P_i|\le (|V_i|-1)+(|V_i|-2)+1$.
	
	Next observe that, as each edge $h\in \ff$ intersects every edge of $\cF_1\cup \cF_2$ and $V[\cF_1]\cap V[\cF_2]=\emptyset$, we have $h=s\cup p$, where $s\in S_i$, $p\in P_j$, $\{i,j\}=\{1,2\}$. Therefore we can split $\ff = \fr \dcup \fl$, where
		\[
			\fr =  \{s\cup p\in \ff \colon s\in S_1, p\in P_2\} \qand
			\fl = \{p\cup s\in \ff \colon p\in P_1, s\in S_2\}.
		\]

	Using the absence of $\c$ and a member of $\p$ in $\cH$, one can prove the following fact. Denote by $B_i\subseteq P_i$, $i=1,2$, the set of 2-covers of $\cF_i$ with at least two neighbors in $\ff$.
	\begin{fact}\label{f:maxdeg}
		For $i=1,2$, $B_i$ is an intersecting family. In particular,
		\begin{eqnarray}\label{eq:maxdeg}
		|\fr|\le |P_2| + (s_1-1)\cdot |B_2|\le |P_2| + (s_1-1)\cdot \max\{3,\Delta(P_2)\},
			\quad \textrm{and}\nonumber \\
		|\fl|\le |P_1| + (s_2-1)\cdot |B_1|\le|P_1| + (s_2-1)\cdot \max\{3,\Delta(P_1)\}.
		\end{eqnarray}
	\end{fact}
	\begin{proof}
		Suppose two 2-covers $p,q\in B_i$ of $\cF_i$ are disjoint and recall that $s_i\le 3$. Then $\cH$ contains either a member of $\p$ (see Figure \ref{fig:fcover}$(a)$) or $\c$ (see Figure \ref{fig:fcover}$(b)$), a contradiction. To see \eqref{eq:maxdeg} recall that the only 2-uniform intersecting families are the triangle and the star, and thus consist of at most $\max\{3,\Delta(P_i)\}$ edges. The inequality $\deg_{\ff}(p)\le s_i$ follows from $N_{\ff}(p)\subseteq S_i$ for every $p\in P_j$, $\{i,j\}=\{1,2\}$.
	\end{proof}
	
\begin{figure}[h!]
		\begin{multicols}{7}
						
							\begin{tikzpicture}[scale = 1]	
							
							\node at (-.1,.4) {\tiny $(a)$};

							\coordinate (v1) at (.3,.4);
							\coordinate (v2) at (.3,0);
							\coordinate (v3) at (.3,-.4);
							
							\coordinate (x) at (1.3,.2);
							\coordinate (y) at (1.1,.5);
							
							\coordinate (a) at (1.3,-.2);
							\coordinate (b) at (1.1,-.5);
							
							\qedge{(y)}{(x)}{(v1)}{2.5pt}{1pt}{black!30!white}{black!10!white,opacity=0.2};
							\qedge{(y)}{(x)}{(v2)}{2.5pt}{1pt}{black!30!white}{black!10!white,opacity=0.2};
							\qedge{(a)}{(b)}{(v2)}{2.5pt}{1pt}{black!30!white}{black!10!white,opacity=0.2};
							\qedge{(a)}{(b)}{(v3)}{2.5pt}{1pt}{black!30!white}{black!10!white,opacity=0.2};
							
							\draw (x) edge [color = green!60!black,thick] (y);
							\draw (a) edge [color = green!60!black,thick] (b);
							
							\foreach \i in {v1,v2,v3,x,y,a,b}
							\fill (\i) circle (1.3pt);		
							\node[green!60!black] at (1.4,.45){\tiny $p$};
							\node[green!60!black] at (1.4,-.45){\tiny $q$};

							\end{tikzpicture}

							\begin{tikzpicture}[scale = 1]				
							\node at (-.1,.4) {\tiny $(b)$};
							
							\coordinate (v1) at (.3,.3);
							\coordinate (v2) at (.3,-.3);
							
							\coordinate (x) at (1.3,.2);
							\coordinate (y) at (1.1,.5);
							
							\coordinate (a) at (1.3,-.2);
							\coordinate (b) at (1.1,-.5);
							
							\qedge{(y)}{(x)}{(v1)}{2.5pt}{1pt}{black!30!white}{black!10!white,opacity=0.2};
							\qedge{(y)}{(x)}{(v2)}{2.5pt}{1pt}{black!30!white}{black!10!white,opacity=0.2};
							\qedge{(a)}{(b)}{(v1)}{2.5pt}{1pt}{black!30!white}{black!10!white,opacity=0.2};
							\qedge{(a)}{(b)}{(v2)}{2.5pt}{1pt}{black!30!white}{black!10!white,opacity=0.2};
							
							\draw (x) edge [color = green!60!black,thick] (y);
							\draw (a) edge [color = green!60!black,thick] (b);
							
							\foreach \i in {v1,v2,x,y,a,b}
							\fill (\i) circle (1.3pt);			
							
							\node[green!60!black] at (1.4,.45){\tiny $p$};
							\node[green!60!black] at (1.4,-.45){\tiny $q$};
							
							\end{tikzpicture}

							\begin{tikzpicture}[scale = 1]		
							\node at (-.8,.5) {\tiny $(c)$};
							
							\foreach \i in {1,...,5}
							\coordinate (v\i) at (18+72*\i:.6cm);
							
							\qedge{(v2)}{(v1)}{(v5)}{3pt}{1pt}{red!70!black}{red!50!white,opacity=0.2};
							\draw (v3) edge [green!70!black, very thick] (v4);

							\foreach \i in {1,...,5}
							\fill (v\i) circle (1.3pt);			
							\end{tikzpicture}
						
							\begin{tikzpicture}[scale = 1]	
							
							\node at (-.2,.4) {\tiny $(d)$};
								\node at (0,-.4) {\tiny $V_1$};
							
							\coordinate (s) at (1.1,0);
							
							\coordinate (v1) at (.4,.4);
							\coordinate (v2) at (.3,0);
							\coordinate (v3) at (.4,-.4);
							
							\coordinate (x) at (1.3,.3);
							\coordinate (y) at (1.1,.5);
							
							\coordinate (a) at (1.3,-.3);
							\coordinate (b) at (1.1,-.5);
							
							\qedge{(y)}{(x)}{(v1)}{2.5pt}{1pt}{black!30!white}{black!10!white,opacity=0.2};
						%	\qedge{(y)}{(x)}{(v2)}{2.5pt}{1pt}{black!30!white}{black!10!white,opacity=0.2};
						%	\qedge{(a)}{(b)}{(v2)}{2.5pt}{1pt}{black!30!white}{black!10!white,opacity=0.2};
							\qedge{(a)}{(b)}{(v3)}{2.5pt}{1pt}{black!30!white}{black!10!white,opacity=0.2};
							
							\qedge{(s)}{(v3)}{(v2)}{2.5pt}{1pt}{black!30!white}{black!10!white,opacity=0.2};
							\qedge{(s)}{(v2)}{(v1)}{2.5pt}{1pt}{black!30!white}{black!10!white,opacity=0.2};
							
							\qedge{(v3)}{(v2)}{(v1)}{6pt}{1pt}{black!30!white}{black!10!white,opacity=0.2};
							
							\draw (x) edge [color = green!60!black,thick] (y);
							\draw (a) edge [color = green!60!black,thick] (b);
							
							\foreach \i in {v1,v2,v3,x,y,a,b}
							\fill (\i) circle (1.3pt);		
							\fill[red!70!black] (s) circle (1.3pt); 
							\node[green!60!black] at (1.4,.5){\tiny $p$};
							\node[green!60!black] at (1.4,-.5){\tiny $q$};
							\node [red!60!black, right] at (s) {\tiny $s$};

							\end{tikzpicture}
							
								\begin{tikzpicture}[scale = 1]	
								
								\node at (-.2,.4) {\tiny $(e)$};
								\node at (0,-.4) {\tiny $V_1$};
								
								\coordinate (s) at (1,.5);
								
								\coordinate (v1) at (.4,.4);
								\coordinate (v2) at (.3,0);
								\coordinate (v3) at (.4,-.4);
								
								\coordinate (x) at (1.3,0);
								\coordinate (y) at (1.2,.2);
								
								\coordinate (a) at (1.3,-.3);
								\coordinate (b) at (1.1,-.5);
								
								\qedge{(y)}{(x)}{(v3)}{2.5pt}{1pt}{black!30!white}{black!10!white,opacity=0.2};
								%	\qedge{(y)}{(x)}{(v2)}{2.5pt}{1pt}{black!30!white}{black!10!white,opacity=0.2};
								%	\qedge{(a)}{(b)}{(v2)}{2.5pt}{1pt}{black!30!white}{black!10!white,opacity=0.2};
								\qedge{(a)}{(b)}{(v3)}{2.5pt}{1pt}{black!30!white}{black!10!white,opacity=0.2};
								
								\qedge{(s)}{(y)}{(v2)}{2.5pt}{1pt}{black!30!white}{black!10!white,opacity=0.2};
								\qedge{(s)}{(v2)}{(v1)}{2.5pt}{1pt}{black!30!white}{black!10!white,opacity=0.2};
								
								\qedge{(v3)}{(v2)}{(v1)}{6pt}{1pt}{black!30!white}{black!10!white,opacity=0.2};
								
								\draw (x) edge [color = green!60!black,thick] (y);
								\draw (a) edge [color = green!60!black,thick] (b);
								
								\foreach \i in {v1,v2,v3,x,y,a,b}
								\fill (\i) circle (1.3pt);		
								\fill[red!70!black] (s) circle (1.3pt); 
								\node[green!60!black] at (1.45,.15){\tiny $p$};
								\node[green!60!black] at (1.4,-.5){\tiny $q$};
								\node [red!60!black, right] at (s) {\tiny $s$};

								\end{tikzpicture}
								
									\begin{tikzpicture}[scale = 1]	
									
									\node at (-.2,.4) {\tiny $(f)$};
									\node at (0,-.4) {\tiny $V_1$};
									
									\coordinate (s) at (1,.5);
									
									\coordinate (v1) at (.4,.4);
									\coordinate (v2) at (.3,0);
									\coordinate (v3) at (.4,-.4);

									\coordinate (c) at (1.1,0);									
									\coordinate (a) at (1.2,-.25);
									\coordinate (b) at (1.1,-.5);
									
									\qedge{(c)}{(a)}{(v3)}{2.5pt}{1pt}{black!30!white}{black!10!white,opacity=0.2};
									%	\qedge{(y)}{(x)}{(v2)}{2.5pt}{1pt}{black!30!white}{black!10!white,opacity=0.2};
									%	\qedge{(a)}{(b)}{(v2)}{2.5pt}{1pt}{black!30!white}{black!10!white,opacity=0.2};
									\qedge{(a)}{(b)}{(v3)}{2.5pt}{1pt}{black!30!white}{black!10!white,opacity=0.2};
									
									\qedge{(s)}{(c)}{(v2)}{2.5pt}{1pt}{black!30!white}{black!10!white,opacity=0.2};
									\qedge{(s)}{(v2)}{(v1)}{2.5pt}{1pt}{black!30!white}{black!10!white,opacity=0.2};
									
									\qedge{(v3)}{(v2)}{(v1)}{6pt}{1pt}{black!30!white}{black!10!white,opacity=0.2};
									
									\draw (a) edge [color = green!60!black,thick] (c);
									\draw (a) edge [color = green!60!black,thick] (b);
									
									\foreach \i in {v1,v2,v3,c,a,b}
									\fill (\i) circle (1.3pt);		
									\fill[red!70!black] (s) circle (1.3pt); 
									\node[green!60!black] at (1.4,0){\tiny $p$};
									\node[green!60!black] at (1.4,-.5){\tiny $q$};
									\node [red!60!black, right] at (s) {\tiny $s$};

									\end{tikzpicture}
									
										\begin{tikzpicture}[scale = 1]	
										
										\node at (-.2,.4) {\tiny $(g)$};
										\node at (0,-.4) {\tiny $V_1$};
										
										\coordinate (s) at (1.1,.3);
										
										\coordinate (v1) at (.4,.4);
										\coordinate (v2) at (.3,0);
										\coordinate (v3) at (.4,-.4);
										
										\coordinate (x) at (.9,.55);
										
										\coordinate (c) at (1.1,0);									
										\coordinate (a) at (1.2,-.25);
										\coordinate (b) at (1.1,-.5);
										
										\qedge{(c)}{(a)}{(v2)}{2.5pt}{1pt}{black!30!white}{black!10!white,opacity=0.2};
										%	\qedge{(y)}{(x)}{(v2)}{2.5pt}{1pt}{black!30!white}{black!10!white,opacity=0.2};
										%	\qedge{(a)}{(b)}{(v2)}{2.5pt}{1pt}{black!30!white}{black!10!white,opacity=0.2};
										\qedge{(a)}{(b)}{(v3)}{2.5pt}{1pt}{black!30!white}{black!10!white,opacity=0.2};
										
										\qedge{(x)}{(s)}{(v1)}{2.5pt}{1pt}{black!30!white}{black!10!white,opacity=0.2};
										\qedge{(s)}{(v2)}{(v1)}{2.5pt}{1pt}{black!30!white}{black!10!white,opacity=0.2};
										
										\qedge{(v3)}{(v2)}{(v1)}{6pt}{1pt}{black!30!white}{black!10!white,opacity=0.2};
										
										\draw (a) edge [color = green!60!black,thick] (c);
										\draw (a) edge [color = green!60!black,thick] (b);
										
										\foreach \i in {v1,v2,v3,c,a,b,x}
										\fill (\i) circle (1.3pt);		
										\fill[red!70!black] (s) circle (1.3pt); 
										\node[green!60!black] at (1.4,0){\tiny $p$};
										\node[green!60!black] at (1.4,-.5){\tiny $q$};
										\node [red!60!black, right] at (s) {\tiny $s$};

										\end{tikzpicture}

\end{multicols}

			\vspace{-1em}
			
		\caption{The illustration of the proofs of Fact \ref{f:maxdeg} and Claim \ref{cl:ss}\ref{it:s2}.}
		\label{fig:fcover}
				
	\end{figure}
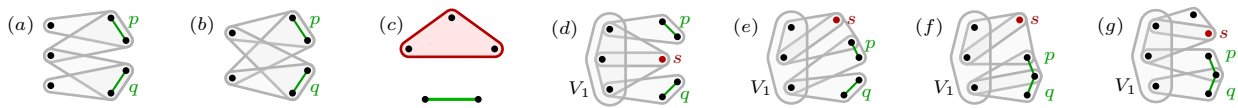
	
%	\vspace{-1em}
	
	For $|V_i|=4$, $i=1,2$, every pair of vertices of $V_i$ is a 2-cover of $\cF_i$ and thereby $P_i=K^{(2)}_4[V_i]$, yielding $|P_i|=6$ and $\Delta(P_i) = 3$. Thus, in this case \eqref{eq:maxdeg} reads as,
		\begin{equation}\label{eq:v4}
			\textrm{If }\quad|V_i|=4\quad \textrm{ then }\quad |\fl|\le 3s_2 + 3 \quad\textrm{or}\quad
			|\fr|\le 3s_1 + 3,\quad\textrm{for $i=1,2$, respecively.}
		\end{equation}

	Now observe that for $|V_i|=5$, $i=1,2$, each 3-edge of $\cF_i$ is disjoint from exactly one pair of vertices of $V_i$ (see Figure \ref{fig:fcover}$(c)$). Therefore, for all distinct $x,y\in V_i$ either $\{x,y\}\in P_i$ or $V_i\setminus \{x,y\} \in \cF_i$, and hence
		\begin{equation}\label{eq:fp}
			\textrm{If }\quad|V_i|=5\quad \textrm{ then }\quad|\cF_i| + |P_i| = {5\choose 2}=10.
		\end{equation}
	Combining this equality with $\Delta(P_2)\le 4$ and $|\fr|\le |P_2| + 4(s_1-1)$ ensured by \eqref{eq:maxdeg}, one gets
		\begin{equation}\label{eq:v5}
			\textrm{If }\quad|V_2|=5\qand s_1 \ge 1\quad\textrm{ then }\quad |\cF_2|+|\fr|\le 4s_1 + 6.
		\end{equation}

\medskip
	
	For the rest of the proof we assume $s_2\le 1$ and if $s_2=1$, denote by $s$ the unique element of $S_2$.
	
\medskip
	\noindent{\it Proof of \ref{it:s2}.}
	We let $s_1=3$ and thereby $|\cF_1|=1$ and $|V_2|=n-3\ge 5$. As $\cH\nsubseteq S_n^{+1}$ and $\cH\nsubseteq SP_n$, there are at least two edges $h,h'\in \fr$ disjoint from $s$. Further, if possible, we choose such $h$, $h'$ so that $p\coloneq h\cap V_2$ and $q\coloneq h'\cap V_2$ are distinct.
	
	First observe, that $p\neq q$. Indeed, otherwise, by our choice of $h$ and $h'$, all edges in $\fr-s$ share the same pair $p\subseteq V_2\setminus \{s\}$.
	This implies that $p\in B_2$ and $|\fr-s|\le 3$.
	By~\eqref{eq:degv3}, $s\notin p$ and $B_2$ is intersecting, the former implies that $|B_2-p|\le 2$.
	Thus, the number of edges in $\fr$ that contain $s$ is at most $(|V_2|-1)+(s_1-1)|B_2-p|\le n$, implying that $|\fr|\le n+3$.
	As every edge of $\cF_2$ intersects both $s$ and $p$, one can estimate
		\[
			|\cF_2|\le 2(|V_2|-3) + 1 = 2n-11.
		\]
	Putting everything together, we obtain for $n\ge 8$,
		\[
			|\cH|=|\cF_1|+|\fl| + |\cF_2| +|\fr|\le 1 + 3 + (2n-11) + (n+3) = 3n-4 < 4n - 11.
		\]

	Now we proceed by induction on $n\ge 8$ and first consider the base case $n=8$.
	For the sake of contradiction suppose $|\cH|\ge 22$. Since by Fact \ref{f:maxdeg} $B_2$ is intersecting, $|B_2|\le 4$ and thus,
		\[
			|\cH|=|\cF_1|+|\cF_{12}^{\rhd}|+|\cF_2|+|\cF_{12}^{\lhd}|
				\overset{\eqref{eq:maxdeg}}{\le}
			1 + 3 + |\cF_2|+|P_2| + (s_1-1)\cdot |B_2|
				\overset{\eqref{eq:fp}}{=} 14 + 2|B_2|\le 22,
		\]
	where we used $|\fl|\le 3$.
	Therefore the equalities go through meaning $|\fl|=3$, $|B_2|=4$, and $\deg_{\ff}(r)=3$ for every $r\in B_2$. This, in turn, entails that the link graph of $s$ in $\cF_{12}$ is a complete bipartite graph $K^{(2)}_{3,4}[V_1\dcup V_2\setminus \{s\}]$ and $B_2$ is a star with center $s$.
	But then, as $p\neq q$, no matter where the edges $h,h'\in \fr-s$ are, $\cH$ contains a minimal 4-path, a contradiction (see Figure~\ref{fig:fcover}$(d)$-$(g)$).

	Next suppose $n\ge 9$ and we shall find a vertex of degree at most four so that we could apply induction and conclude the proof. If there are two edges $f_1, f_2\in \cF_2$ with $f_1\cap f_2 = \{s\}$, set $U\coloneq f_1\cup f_2$. Otherwise, $\cF_2 = K_4^{(3)}-e$ and we define $U\coloneq V[\cF_2]$. Clearly $|U|\le 5$ and thus we can take a vertex $v\in V_2\setminus U$. Now every 2-cover $p\in P_2-s$ of $\cF_2$ is entirely contained in $U$ and therefore the only neighbor in $P_2$ of $v$  is $s$, yielding $\deg_{\fr}(v)=|N_{\fr}(sv)|\le |S_1|=3$. Moreover, because every edge $f\in \cF_2$ contains $s$ and intersects both 2-covers $p,q\in P_2-s$, we have $\deg_{\cF_2}(v) \le 1$. As $\deg_{\fl}(v)=\deg_{\cF_1}(v)=0$, altogether we obtain $\deg_{\cH}(v)\le 4$ and we are done.
	
\medskip
	
	\noindent{\it Proof of \ref{it:s3}.}
	Let $s_1 = 2$ and $s_2=1$, yielding $|V_1|\ge 4$, $|V_2|\ge 4$, $|\cF_1|\le |V_1|-2$, and $|\cF_2|\le {|V_2|-1\choose 2}$.  Moreover, in view of \eqref{eq:p7} one gets $|\fl| \le	|P_1|\cdot s_2 =|P_1| \le 2|V_1|-2$.	
	Now, if there exists an edge $h\in \fr$ disjoint from $S_2$, then because each edge of $\cF_2$ contains $s$ and intersects $h$, we have $|\cF_2|\le 2(|V_2|-3)+1 = 2|V_2|-5$. Further, applying \eqref{eq:p7} combined together with \eqref{eq:maxdeg} yields,
		\[
			|\fr| \le |P_2|+(s_1-1)\cdot \max\{3,\Delta(P_2)\} \le (|V_2|+3) + (|V_2|-1) = 2|V_2| + 2.
		\]
	Summarizing,
		\[
			|\cH|=|\cF_1|+|\fl|+|\cF_2|+|\fr|\le (|V_1|-2) + (2|V_1|-2) + (2|V_2|-5) + (2|V_2|+2)
			= 4n-|V_1|-7 \le 4n-11.
		\]

	Otherwise $s$ is contained in all edges of $\fr$ (so in fact in all edges of $\cF_{12}\cup \cF_2$) and thus $|\fr|\le 2(|V_2|-1)$. Because $\cH\nsubseteq CB_n$, we have $|V_1|\ge 5$, entailing $|V_2| \le n-5$. Then,
		\[
			|\cH|\le (|V_1| -2) + (2|V_1|-2) + {|V_2|-1\choose 2} + 2(|V_2|-1)
			={|V_2|-2\choose 2} + 3n-8\le {n-4\choose 2} +10.
		\]

\medskip	

	Before we proceed observe that for $\{i,j\}=\{1,2\}$ and each $s'\in S_j$, $\cH[V_i\cup \{s'\}]$ is an intersecting family. Indeed, this follows from that $\cF_i= \cH[V_i]$ is intersecting, the pair $(\cF_i,\ff)$ is cross-intersecting, and each edge $h\in \cH[V_i\cup \{s'\}]$ with $s'\in h$ is in $\ff$. Therefore the celebrated Erd\H os--Ko--Rado theorem \cite{EKR } tells us, that for $|V_i|\ge 5$,
		\begin{equation}\label{eq:VsT}
		|\cH[V_i\cup \{s'\}]|\le {|V_i|\choose 2}.
		\end{equation}
	Moreover, if there is an edge $h\in \cH[V_i\cup \{s'\}]$ such that $h\cap S_i=\emptyset$, then for $|V_i|\ge 5$,
		\begin{equation}\label{eq:VsN}
			|\cH[V_i\cup \{s'\}]|\le 3|V_i|-5.
		\end{equation}
	For $|V_i|=5$ the above bound follows from \eqref{eq:VsT}, whereas for $|V_i|\ge 6$ one can use Hilton-Milner theorem (Theorem \ref{th:hm}), as $\cF_i\cup \ff[V_i\cup \{s'\}]$ is a non-trivial intersecting family. This is because only vertices of $S_i$ belong to all edges of $\cF_i$ and $\cF_i\neq \emptyset$, but $h\cap S_i=\emptyset$.

\medskip

	\noindent{\it Proof of \ref{it:s4}.}	
	We let $s_1=s_2=1$, which entails $|V_i|\ge 4$. Observe, that
		\[
			\cH[V_1\cup S_2] = \cF_1 \cup \fl \qqand \cH[V_2\cup S_1] = \cF_2 \cup \fr.
		\]
	Therefore, as clearly for $|V_i|=4$ we have $|\cH[V_i\cup S_j]|\le {5\choose 3}=10$, $\{i,j\}=\{1,2\}$, in view of \eqref{eq:VsT},
		\[
			|\cH| \le \max\left\{10, {|V_1|\choose 2}\right\} + \max\left\{10,{|V_2|\choose 2}\right\} \le {n-4\choose 2}+10,
		\]
	for $n\ge 9$ \footnote{In particular, when $|V_1|, |V_2|\ge 5$, the inequality can be checked using $|V_1||V_2|\ge 4(n-5)\ge \binom n2 - \binom{n-4}2-10$.}, whereas for $n=8$ one gets $|\cH|\le 10+10\le 4\cdot 8 - 11$.
	
\medskip
	\noindent{\it Proof of \ref{it:s5}.}
	Let $s_2=0$ and thereby $\fl = \emptyset$ yielding $\cH = \cF_1\dcup \cF_2\dcup \fr$. Moreover, $|V_2|\ge 4$ and since $\fr=\ff\neq\emptyset$, $s_1\ge 1$, which implies $|\cF_1|\le {|V_1|-1\choose 2}$. If $|V_2|=4$, then $|V_1|=n-4 \ge 4$ entailing $s_1\le 2$, and therefore \eqref{eq:v4} tells us $|\fr| \le  9$. Thus, for $n\ge 8$ we have,
		\[
			|\cH|= |\cF_1|+ |\cF_2| + |\fr|\le {n-5\choose 2} + 4 + 9 \le {n-4\choose 2} + 10.
		\]

	Now, let $|V_2|\ge 5$, and pick any $s' \in S_1$. In view of $S_2=\emptyset$, \eqref{eq:VsN} tells us, $|\cH[V_2\cup \{s'\}]|\le 3|V_2|-5$. Moreover, by the definition of $B_2$,~\eqref{eq:degv3}, \eqref{eq:p7}, and Fact \ref{f:maxdeg} we have $|\fr-s'|\le |P_2|+|B_2|\le 7 + 3=10$. Summarizing,
		\[
			|\cH|\le {|V_1|-1\choose 2} + (3|V_2|-5)+10\le \max\left\{4n-11,{n-4\choose 2} +10\right\},
		\]
	where $|\cH|\le 4n-11$ can be checked for $3\le |V_1|\le 5$, and $|\cH|\le {n-4\choose 2}+10$ for $|V_1|\ge 6$ \footnote{The inequality holds for $|V_2|=5$. For $|V_2|\ge 6$, we have ${n-4\choose 2} - {|V_1|-1\choose 2} = {n-4\choose 2} - {n-|V_2|-1\choose 2}\ge 3(n-6) \ge 3|V_2|$.}.
\end{proof}

%%%%%%%%%%%%%%%%%%%%%%%%%%%%%%%%%%%%%%%%%%%%%%%%%%%%%%%%%%%%%%%%%%%%%%
%																																		%
%																												   						%
%  									nu   =   3 																   						%
%																																		%
%																																		%
%%%%%%%%%%%%%%%%%%%%%%%%%%%%%%%%%%%%%%%%%%%%%%%%%%%%%%%%%%%%%%%%%%%%%%

\subsection{Proof of Lemma \ref{l:M3}}

Let $\cH$ be a connected $\p$-free 3-graph on the set of vertices $V$, $|V|=n\ge 9$, with $\nu(\cH)=3$. Lemma \ref{l:M3} follows from
	\begin{equation}\label{eq:lM3}
		|\cH|\le \binom{n-4}{2}+11,
	\end{equation}
with the equality achieved if and only if $\cH$ is the balloon $B_n$.
To prove this inequality we let
		\[
			\cH = \cF_1\dcup \cF_2\dcup \cF_3\dcup \ff\dcup \fff
		\]
to be a partition guaranteed by Lemma \ref{l:split3}. Set $V_1=V[\cF_1]$, $V_3=V[\cF_3]$ and $V_2=V\setminus (V_1\cup V_3)$, and recall
	\begin{enumerate}[label=\rmlabel]
		\item\label{it:n31} $V[\cF_2]\subset V_2$, $V_1\cap V_3=\emptyset$, and $V[\ff]\cap V_3=\emptyset$,
		\item\label{it:n32} the 3-graphs $\cF_1$, $\cF_2$, and $\cF_3$ are non-empty intersecting families,
		\item\label{it:n33} $\fff\neq\emptyset$,
		\item\label{it:n34} the pairs $(\cF_1\dcup \cF_2\dcup \cF_3\dcup \ff,  \fff)$ and $(\cF_1\dcup \cF_2,\ff)$ are cross-intersecting.
	\end{enumerate}
	
Further, for each $i=1,2,3$ pick an edge $e_i\in \cF_i$ and split the set of edges of $\ff$ into two subsets, $\ff=\ffi\dcup \ffo$, where
		\[
			\ffi=\{f\in \ff: f\subset e_1\cup e_2\} \qand  \ffo=\{f\in \ff: |f\cap e_1|=|f\cap e_2|=1\}.			
		\]
Because every edge of $\ff$ intersects both $e_1$ and $e_2$, we have
		\begin{equation}\label{eq:Hsplit}
			\cH = \cF_1\dcup \cF_2 \dcup \cF_3 \dcup \ffi\dcup \ffo \dcup \fff.
		\end{equation}

The proof of \eqref{eq:lM3} mainly relies on two technical claims enabling us to bound the number of edges in $\ff\cup \fff$. In the first of them we estimate the size of $\fff\cup \ffi$.
		\begin{claim}\label{cl:in}
			$|\fff|+|\ffi|\le 18$. Moreover, if $|\fff|+|\ffi| = 18$ then $\fff\cup \ffi$ is a star.
		\end{claim}
		\begin{proof}
			As every edge $h\in \fff$ intersects each one of $e_1$, $e_2$ and $e_3$, we trivially have $|\fff|\le27$, but this estimate can be improved. Let $G\subseteq K^{(2)}_{3,3}[e_1\cup e_2]$ be an auxiliary bipartite graph with vertex classes $e_1$ and $e_2$, consisting of all pairs $\{u,v\}\in e_1\times e_2$ for which there exists a vertex $w\in e_3$ such that $uvw\in \fff$. It turns out that the number of edges in $\fff$ can exceed $|G|$ only by at most 6,
				\begin{equation}\label{eq:fg}
						|\fff| \le |G| + 6.
				\end{equation}
\def\mr{\textcolor{red!60!black}{\textrm{M}^{\textrm{R}}}}
\def\mb{\textcolor{blue!60!black}{\textrm{M}^{\textrm{B}}}}
\def\mg{\textcolor{green!60!black}{\textrm{M}^{\textrm{G}}}}
\def\mm{\textrm{M}^i}				
			Indeed, clearly any edge of $G$ can be extended to at most 3 edges of $\fff$ (see Figure \ref{fig:f88}$(a)$). However, due to the $\p$-freeness of $\cH$, there can be no two disjoint edges $f_1,f_2\in G$ and three different vertices $w_1,w_2,w_3\in e_3$, such that $f_1w_1, f_1w_2, f_2w_2, f_2w_3$ are all edges in $\fff$, as they would form a minimal 4-path in $\cH$ (see Figure \ref{fig:f88}$(b)$). Similarly, there are no disjoint edges $f_1, f_2, f_3\in G$ and vertices $w_1,w_2\in e_3$ with $f_1w_1, f_2w_1, f_2w_2, f_3w_2\in \fff$ (see Figure \ref{fig:f88}$(c)$). To avoid  such structures, any two disjoint edges in $G$ can be extended, in total, to at most 4 edges of $\fff$, and any three disjoint edges of $G$ can be extended, in total, to at most 5 edges of $\fff$. Therefore, to conclude \eqref{eq:fg} it is enough to observe, that the set of edges of $K^{(2)}_{3,3}$ can be partitioned into three disjoint matchings $M^{(2)}_3$, say $\mr, \mg, \mb$ (see Figure \ref{fig:f88}$(d)$). Now, for each $i\in \{\textcolor{red!60!black}{R},\textcolor{green!60!black}{G},\textcolor{blue!60!black}{B}\}$, $G\cap \mm$ can be extended to at most $|G\cap \mm|+2$ edges of $\fff$.

%\vspace{-1em}
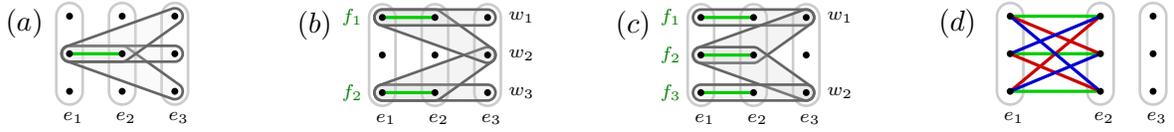
\begin{figure}[h!]
		\begin{multicols}{4}

				\begin{tikzpicture}[scale = 1]

				\coordinate (a) at (.1,1.4);
				\node at (a) {$(a)$};

				\foreach \i in {1,2,3}{
					\foreach \j in {1,2,3}{
						\coordinate (v\i\j) at (.7*\i,.5*\j);
					}}
					
					\node at (.75,.15) {\tiny $e_1$};
					\node at (1.45,.15) {\tiny $e_2$};
					\node at (2.15,.15) {\tiny $e_3$};
					
						\qedge{(v11)}{(v12)}{(v13)}{5pt}{1pt}{white!80!black}{black!10!white,opacity=0};
						\qedge{(v21)}{(v22)}{(v23)}{5pt}{1pt}{white!80!black}{black!10!white,opacity=0};
						\qedge{(v31)}{(v32)}{(v33)}{5pt}{1pt}{white!80!black}{black!10!white,opacity=0};

					\qedge{(v12)}{(v22)}{(v31)}{3pt}{1pt}{white!40!black}{black!10!white,opacity=0.2};
					\qedge{(v12)}{(v22)}{(v32)}{3pt}{1pt}{white!40!black}{black!10!white,opacity=0.2};
					\qedge{(v33)}{(v22)}{(v12)}{3pt}{1pt}{white!40!black}{black!10!white,opacity=0.2};

					\draw [green!80!black, very thick] 
					(v12) edge (v22)
				;
					
					\foreach \i in {1,2,3}{
						\foreach \j in {1,2,3}{
							\fill (v\i\j) circle (1.3pt);}}
					
					\end{tikzpicture}

				\begin{tikzpicture}[scale = 1]	
					\coordinate (b) at (-.2,1.4);			
				\node at (b) {$(b)$};

			\node at (.75,.15) {\tiny $e_1$};
			\node at (1.45,.15) {\tiny $e_2$};
			\node at (2.15,.15) {\tiny $e_3$};
				
				\qedge{(v11)}{(v12)}{(v13)}{5pt}{1pt}{white!80!black}{black!10!white,opacity=0};
				\qedge{(v21)}{(v22)}{(v23)}{5pt}{1pt}{white!80!black}{black!10!white,opacity=0};
				\qedge{(v31)}{(v32)}{(v33)}{5pt}{1pt}{white!80!black}{black!10!white,opacity=0};

					\qedge{(v31)}{(v21)}{(v11)}{3pt}{1pt}{white!40!black}{black!20!white,opacity=0.2};
					\qedge{(v32)}{(v21)}{(v11)}{3pt}{1pt}{white!40!black}{black!20!white,opacity=0.2};
					\qedge{(v13)}{(v23)}{(v32)}{3pt}{1pt}{white!40!black}{black!20!white,opacity=0.2};
					\qedge{(v13)}{(v23)}{(v33)}{3pt}{1pt}{white!40!black}{black!20!white,opacity=0.2};
					
					\draw [green!80!black, very thick] 
					(v11) edge (v21)
					(v13) edge (v23);
					
					\node [green!50!black] at ($(v11)+(-.4,0)$) {\tiny $f_2$};
					\node [green!50!black] at ($(v13)+(-.4,0)$) {\tiny $f_1$};
					
					\node at ($(v33)+(.45,0)$){\tiny $w_1$};
					\node at ($(v32)+(.45,0)$){\tiny $w_2$};
					\node at ($(v31)+(.45,0)$){\tiny $w_3$};
					
					\foreach \i in {1,2,3}{
						\foreach \j in {1,2,3}{
							\fill (v\i\j) circle (1.3pt);}}
					
					\end{tikzpicture}
				
			\begin{tikzpicture}[scale = 1]				
			\node at (b) { $(c)$};
			
				\node at (.75,.15) {\tiny $e_1$};
				\node at (1.45,.15) {\tiny $e_2$};
				\node at (2.15,.15) {\tiny $e_3$};
				
				\qedge{(v11)}{(v12)}{(v13)}{5pt}{1pt}{white!80!black}{black!10!white,opacity=0};
				\qedge{(v21)}{(v22)}{(v23)}{5pt}{1pt}{white!80!black}{black!10!white,opacity=0};
				\qedge{(v31)}{(v32)}{(v33)}{5pt}{1pt}{white!80!black}{black!10!white,opacity=0};

				\qedge{(v33)}{(v23)}{(v13)}{3pt}{1pt}{white!40!black}{black!20!white,opacity=0.2};
				\qedge{(v33)}{(v22)}{(v12)}{3pt}{1pt}{white!40!black}{black!20!white,opacity=0.2};
				\qedge{(v12)}{(v22)}{(v31)}{3pt}{1pt}{white!40!black}{black!20!white,opacity=0.2};
				\qedge{(v11)}{(v21)}{(v31)}{3pt}{1pt}{white!40!black}{black!20!white,opacity=0.2};
				
				\draw [green!80!black, very thick] 
				(v12) edge (v22)
				(v11) edge (v21)
				(v13) edge (v23);
							
								\foreach \i in {1,2,3}{
									\foreach \j in {1,2,3}{
										\fill (v\i\j) circle (1.3pt);}}

					\node [green!50!black] at ($(v11)+(-.4,0)$) {\tiny $f_3$};
					\node [green!50!black] at ($(v13)+(-.4,0)$) {\tiny $f_1$};
					\node [green!50!black] at ($(v12)+(-.4,0)$) {\tiny $f_2$};
					
					\node at ($(v33)+(.45,0)$){\tiny $w_1$};
			%		\node at ($(v32)+(.6,0)$){\tiny $w_2$};
					\node at ($(v31)+(.45,0)$){\tiny $w_2$};			
							
			\end{tikzpicture}
						
			\begin{tikzpicture}[scale = 1]				
				\node at (.5,1.4) {$(d)$};

							\foreach \i in {1,2}
								\foreach \j in {1,2,3}
									\coordinate (v\i\j) at (\i*1.2,.5*\j);

								\foreach \j in {1,2,3}{
									\coordinate (v3\j) at (3.1,.5*\j);	
									\fill (v3\j) circle (1.3pt);	
								}
									
									\node at (1.25,.15) {\tiny $e_1$};
									\node at (2.45,.15) {\tiny $e_2$};
									\node at (3.15,.15) {\tiny $e_3$};
									
									\qedge{(v11)}{(v12)}{(v13)}{5pt}{1pt}{white!80!black}{black!10!white,opacity=0};
									\qedge{(v21)}{(v22)}{(v23)}{5pt}{1pt}{white!80!black}{black!10!white,opacity=0};
									\qedge{(v31)}{(v32)}{(v33)}{5pt}{1pt}{white!80!black}{black!10!white,opacity=0};

					\draw [green!80!black, very thick] 
					(v12) edge (v22)
					(v11) edge (v21)
					(v13) edge (v23);
					
					\draw [red!80!black, very thick] 
					(v12) edge (v21)
					(v11) edge (v23)
					(v13) edge (v22);
					
					\draw [blue!80!black, very thick] 
					(v12) edge (v23)
					(v11) edge (v22)
					(v13) edge (v21);

					\foreach \i in {1,2}
						\foreach \j in {1,2,3}
							\fill (v\i\j) circle (1.3pt);

			\end{tikzpicture}
								
			\end{multicols}
		
			\vspace{-1em}
			
		\caption{Extensions of edges of $G$ and decomposition of $K^{(2)}_{3,3}$ into matchings.}	
		%\caption{The illustrations to the proof of \eqref{eq:fg}.}
		\label{fig:f88}
				
	\end{figure}
	
%	\vspace{-1em}

			Next, let us note that
				\begin{equation}\label{eq:g4f6}
					|G|\ge 4 \quad \textrm{entails}\quad |\ffi|\le 6.
				\end{equation}
			To show this, recall that every edge $f\in \ffi$ intersects each $uvw\in \fff$ and thereby also every $uv\in G$. As there are only five pairwise non-isomorphic subgraphs of $K^{(2)}_{3,3}$ with four edges (all listed in Figure \ref{fig:f99}$(a)$-$(e)$), a simple case analysis enables us to establish \eqref{eq:g4f6}.				

\begin{figure}[h!]
		\begin{multicols}{9}
								
			\begin{tikzpicture}[scale = .7]	
				\coordinate (a) at (.45,1.5);
				\node at (a){\scriptsize  $(a)$};
						
				\foreach \i in {1,2}
					\foreach \j in {1,2,3}
						\coordinate (v\i\j) at (\i,.5*\j);
						
				\qedge{(v11)}{(v12)}{(v13)}{5pt}{.3pt}{black}{black!10!white,opacity=.5};
				\qedge{(v21)}{(v22)}{(v23)}{5pt}{.3pt}{black}{black!10!white,opacity=.5};		

				\draw [green!80!black, very thick] 
				(v13) edge (v23)
				(v13) edge (v22)
				(v12) edge (v23)
				(v12) edge (v22)
				;
				
				\foreach \i in {1,2}{
					\foreach \j in {1,2,3}{
						\fill (v\i\j) circle (1.7pt);}}
				\end{tikzpicture}
								
				\begin{tikzpicture}[scale = .7]		
				\node at (a){\scriptsize $(b)$};	
				
					\qedge{(v11)}{(v12)}{(v13)}{5pt}{.3pt}{black}{black!10!white,opacity=.5};
					\qedge{(v21)}{(v22)}{(v23)}{5pt}{.3pt}{black}{black!10!white,opacity=.5};

				\draw [green!80!black, very thick] 
				(v13) edge (v23)
				(v13) edge (v22)
				(v11) edge (v22)
				(v11) edge (v21)
				;				
				\foreach \i in {1,2}{
					\foreach \j in {1,2,3}{
						\fill (v\i\j) circle (1.7pt);}}				
				\end{tikzpicture}		
						
				\begin{tikzpicture}[scale = .7]			
				\node at (a){\scriptsize $(c)$};		
				
					\qedge{(v11)}{(v12)}{(v13)}{5pt}{.3pt}{black}{black!10!white,opacity=.5};
					\qedge{(v21)}{(v22)}{(v23)}{5pt}{.3pt}{black}{black!10!white,opacity=.5};

				\draw [green!80!black, very thick] 
				(v13) edge (v23)
				(v13) edge (v22)
				(v13) edge (v21)
				(v11) edge (v21)
				;				
				\foreach \i in {1,2}{
					\foreach \j in {1,2,3}{
						\fill (v\i\j) circle (1.7pt);}}				
				\end{tikzpicture}
				
				\begin{tikzpicture}[scale = .7]		
				\node at (a){\scriptsize $(d)$};		
				
					\qedge{(v11)}{(v12)}{(v13)}{5pt}{.3pt}{black}{black!10!white,opacity=.5};
					\qedge{(v21)}{(v22)}{(v23)}{5pt}{.3pt}{black}{black!10!white,opacity=.5};		
								
				\draw [green!80!black, very thick] 
				(v13) edge (v23)
				(v13) edge (v22)
				(v12) edge (v21)
				(v11) edge (v21)
				;				
				\foreach \i in {1,2}{
					\foreach \j in {1,2,3}{
						\fill (v\i\j) circle (1.7pt);}}				
				\end{tikzpicture}
				
				\begin{tikzpicture}[scale = .7]	
				\node at (a){\scriptsize $(e)$};			
				
					\qedge{(v11)}{(v12)}{(v13)}{5pt}{.3pt}{black}{black!10!white,opacity=.5};
					\qedge{(v21)}{(v22)}{(v23)}{5pt}{.3pt}{black}{black!10!white,opacity=.5};		
								
				\draw [green!80!black, very thick] 
				(v13) edge (v23)
				(v13) edge (v22)
				(v12) edge (v22)
				(v11) edge (v21)
				;				
				\foreach \i in {1,2}{
					\foreach \j in {1,2,3}{
						\fill (v\i\j) circle (1.7pt);}}				
				\end{tikzpicture}

				\begin{tikzpicture}[scale = .7]	
				\coordinate (b) at ($(a)+(-.05,0)$);			
				\node at (b){\scriptsize $(f)$};
				
					\qedge{(v11)}{(v12)}{(v13)}{6pt}{.3pt}{black}{black!10!white,opacity=0};
					\qedge{(v21)}{(v22)}{(v23)}{6pt}{.3pt}{black}{black!10!white,opacity=0};

				\foreach \i in {1,2}
				\foreach \j in {1,2,3}
				\coordinate (v\i\j) at (\i,.5*\j);
				
				\qedge{(v13)}{(v23)}{(v12)}{4pt}{.8pt}{black}{black!20!white,opacity=0.2};
				
				\draw [green!80!black, very thick] 
				(v11) edge (v22)
				(v11) edge (v21)
				;
				
				\foreach \i in {1,2}{
					\foreach \j in {1,2,3}{
						\fill (v\i\j) circle (1.7pt);}}
				
				\end{tikzpicture}

				\begin{tikzpicture}[scale = .7]	
				\node at (b){\scriptsize $(g)$};			
				\foreach \i in {1,2}
				\foreach \j in {1,2,3}
				\coordinate (v\i\j) at (\i,.5*\j);
				
					\qedge{(v11)}{(v12)}{(v13)}{6pt}{.3pt}{black}{black!10!white,opacity=0};
					\qedge{(v21)}{(v22)}{(v23)}{6pt}{.3pt}{black}{black!10!white,opacity=0};

				\foreach \i/\j/\k in {13/23/12, 11/12/21} 
				\qedge{(v\i)}{(v\j)}{(v\k)}{4pt}{.8pt}{black}{black!20!white,opacity=0.2};
				
				\draw [green!80!black, very thick] 
				(v11) edge (v21)
				(v11) edge (v22)
				(v13) edge (v22)
				(v13) edge (v23)
				;	
				
				\foreach \i in {1,2}{
					\foreach \j in {1,2,3}{
						\fill (v\i\j) circle (1.7pt);}}
				\end{tikzpicture}

				\begin{tikzpicture}[scale = .7]				
				
				\node at (b){\scriptsize $(h)$};
				
					\qedge{(v11)}{(v12)}{(v13)}{6pt}{.3pt}{black}{black!10!white,opacity=0};
					\qedge{(v21)}{(v22)}{(v23)}{6pt}{.3pt}{black}{black!10!white,opacity=0};

				\foreach \i/\j/\k in {13/22/12, 11/22/21} 
				\qedge{(v\i)}{(v\j)}{(v\k)}{4pt}{.8pt}{black}{black!20!white,opacity=0.2};
				
				\draw [green!80!black, very thick] 
				(v11) edge (v21)
				(v11) edge (v23)
				(v12) edge (v23)
				(v13) edge (v23)
				;
				
				\foreach \i in {1,2}{
					\foreach \j in {1,2,3}{
						\fill (v\i\j) circle (1.7pt);}}
				\end{tikzpicture}
				
				\begin{tikzpicture}[scale = .7]				
				
				\node at (b){\scriptsize $(i)$};
				
					\qedge{(v11)}{(v12)}{(v13)}{6pt}{.3pt}{black}{black!10!white,opacity=0};
					\qedge{(v21)}{(v22)}{(v23)}{6pt}{.3pt}{black}{black!10!white,opacity=0};

				\foreach \i/\j/\k in {13/23/12, 13/23/22,13/22/12, 12/23/22} 
				\qedge{(v\i)}{(v\j)}{(v\k)}{4pt}{.8pt}{black}{black!20!white,opacity=0.2};
				
				\draw [green!80!black, very thick] 
				(v11) edge (v21)
				;
				
				\foreach \i in {1,2}{
					\foreach \j in {1,2,3}{
						\fill (v\i\j) circle (1.7pt);}}
				\end{tikzpicture}
				
		\end{multicols}
		
			\vspace{-1em}
			
		\caption{   All 4-edge subgraphs of $K^{(2)}_{3,3}$ and forbidden edges of $G$ and $F^{in}_{12}$.}
		\label{fig:f99}
				
	\end{figure}
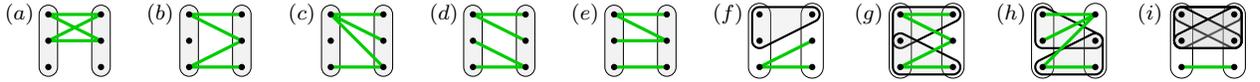

			Finally observe that
				\begin{equation}\label{eq:g5f5}
					\ffi\neq\emptyset  \quad \textrm{entails}\quad |G|\le 7, \qand
					|\ffi|\ge 5 \quad \textrm{yields}\quad |G|\le 5.
				\end{equation}
			Indeed, as $(\ffi,G)$ is cross-intersecting, the existence of any edge in $\ffi$ forbids two pairs from $G$ (see Figure \ref{fig:f99}$(f)$). Moreover, among every 5 edges of $\ffi$ there are two, $f_1, f_2$, sharing at most one vertex. Therefore, as every edge $g\in G$ intersects both $f_1$ and $f_2$, out of all 9 edges of $K^{(2)}_{3,3}$ at least four are forbidden for $G$, yielding $|G|\le 5$ (see Figure \ref{fig:f99}$(g)$,$(h)$).
						
			Now we are ready to finish the proof of Claim \ref{cl:in}. To this end assume
			  	\begin{equation}\label{eq:gf}
					  |\fff|+|\ffi|\ge 18
			  	\end{equation}
			 and note that, in view of \eqref{eq:fg}, this entails $|G| + |\ffi|\ge 12$. Combining this estimate with \eqref{eq:g4f6} and \eqref{eq:g5f5} one can conclude that $|G|\le 3$. Indeed, as $|G|\le 9$ we have $\ffi\neq\emptyset$ and thus $|G|\le 7$. Next, assuming $|G|\ge 4$ we get $|\ffi|\le 6$ and $|G|\le 5$, implying $|G| + |\ffi|\le 11$.
			
			 To exclude $|G|\le 2$ let us recall again that$(\ffi,G)$ is cross-intersecting, and observe that because every edge of $G$ is disjoint from four edges of $\ffi$ (see Figure \ref{fig:f99}$(i)$), $|G|=1$ results $|\ffi|\le 18-4=14$. Similarly, $|G| = 2$ entails $|\ffi|\le 11$. As every edge of $G$ can be extended to at most 3 edges of $\fff$, in both cases $|\fff|+|\ffi|\le 17$. But this, together with $G\neq \emptyset$ guaranteed by \ref{it:n33}, contradicts \eqref{eq:gf}. Thus, $|G|=3$ and thereby $|\ffi|\ge 9$. A quick inspection shows that this is possible only when both $G$ and $\ffi$ are stars with the same center.
		\end{proof}
		
	Our next goal is to bound the number of edges in $\ffo$.
		\begin{claim}\label{cl:out}
			If there exists a vertex $v\in V\setminus (e_1\cup e_2 \cup e_3)$ with $\deg_{\ff}(v)\ge 4$, then $|\cH|\le \binom{n-4}{2}+10$.
		\end{claim}
		\begin{proof}
			We let $v\in V\setminus (e_1\cup e_2 \cup e_3)$ to be a vertex with $\deg_{\ff}(v)\ge 4$. Split the vertex set $V=R\dcup S\dcup V_3$, where
				\[
					R=e_1\cup e_2 \cup \{v\}, \quad
					S = V\setminus (R\cup V_3),
				\]
			and $R\cap V_3=\emptyset$ follows from \ref{it:n31}.
			
			We begin by proving, that every vertex $w\in S$ satisfies
				\begin{equation}\label{eq:deg7}
					\deg_{\cH}(w)\le 7.
				\end{equation}
			Indeed, we let $h\in \fff$ to be an edge guaranteed by \ref{it:n33}, and set $\{x_i\} = h\cap e_i$, $i=1,2,3$. Now \ref{it:n34} tells us, that every edge $f\in \ff$ intersects $h$ and thus contains at least one of the vertices $x_1, x_2$. This entails $\deg_{\ff}(w)\le 5$ (see Figure \ref{fig:f10}$(a)$), and therefore it remains to show that $\deg_{\cF_{1}\cup \cF_2}(w)\le 2$.
			
%\vspace{-1em}
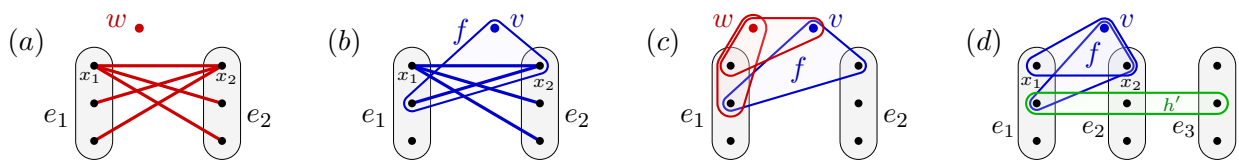
\begin{figure}[h!]
		\begin{multicols}{4}
			
				\begin{tikzpicture}[scale = 1]				
				
					\coordinate (w) at (2.3,2);
					\node [red!70!black] at ($(w)+(-.3,.1)$){$w$};
				
				\coordinate (a) at (.8,1.8);
				\node at (a) {$(a)$};

				\coordinate (v) at (2.5,2);
				
				\foreach \i in {1,2}{
					\foreach \j in {1,2,3}{
						\coordinate (v\i\j) at (1.7*\i,.5*\j);
					}}
					
					\qedge{(v11)}{(v12)}{(v13)}{7pt}{.3pt}{blue!0!black}{black!20!white,opacity=0.2};
					\qedge{(v21)}{(v22)}{(v23)}{7pt}{.3pt}{blue!0!black}{black!20!white,opacity=0.2};
					
					\draw [red!80!black, very thick] 
					(v13) edge (v21)
					(v13) edge (v22)
					(v13) edge (v23)
					(v12) edge (v23)
					(v11) edge (v23)
					;
					
					\foreach \i in {1,2}{
						\foreach \j in {1,2,3}{
							\fill (v\i\j) circle (1.5pt);}}
					
			%		\fill [blue!80!black](v) circle (2pt);
					\fill [red!80!black](w) circle (1.7pt);

		%			\node [blue!70!black] at ($(v)+(.3,.1)$){$v$};
					\node  at ($(v11)+(-.5,.3)$) {$e_1$};
					\node  at ($(v21)+(.5,.3)$){$e_2$};
					\node at ($(v13)+(-.05,-.16)$){\tiny $x_1$};
					\node at ($(v23)+(.06,-.17)$){\tiny $x_2$};

				\end{tikzpicture}
			
				\begin{tikzpicture}[scale = 1]

				\node at (a) { $(b)$};
				
				\coordinate (v) at (2.8,2);

					\qedge{(v11)}{(v12)}{(v13)}{7pt}{.3pt}{blue!0!black}{black!20!white,opacity=0.2};
					\qedge{(v21)}{(v22)}{(v23)}{7pt}{.3pt}{blue!0!black}{black!20!white,opacity=0.2};
					
						\qedge{(v12)}{(v)}{(v23)}{3pt}{.8pt}{blue!80!black}{blue!10!white,opacity=0.2};
				
					\draw [blue!80!black, very thick] 
					(v13) edge (v21)
					(v13) edge (v22)
					(v13) edge (v23)
					(v12) edge (v23)
				;
					
					\foreach \i in {1,2}{
						\foreach \j in {1,2,3}{
							\fill (v\i\j) circle (1.5pt);}}
					
					\fill [blue!80!black](v) circle (1.7pt);

					\node [blue!70!black] at ($(v)+(.3,.1)$){$v$};
				\node  at ($(v11)+(-.5,.3)$) {$e_1$};
					\node  at ($(v21)+(.5,.3)$){$e_2$};
				\node at ($(v13)+(-.05,-.16)$){\tiny $x_1$};
				\node at ($(v23)+(.06,-.25)$){\tiny $x_2$};
					
	%				\node [blue!70!black] at (2.2,.8){ $f$};
						\node [blue!70!black] at (2.35,1.95){$f$};
					
					\end{tikzpicture}

				\begin{tikzpicture}[scale = 1]	
						
				\node at (a) {$(c)$};
				
				\coordinate (w) at (2,2);
				\node [red!70!black] at ($(w)+(-.4,.1)$){$w$};
				
					\qedge{(v11)}{(v12)}{(v13)}{7pt}{.3pt}{blue!0!black}{black!20!white,opacity=0.2};
					\qedge{(v21)}{(v22)}{(v23)}{7pt}{.3pt}{blue!0!black}{black!20!white,opacity=0.2};
					
						\qedge{(v12)}{(v)}{(v23)}{2.8pt}{.8pt}{blue!80!black}{blue!10!white,opacity=0.2};

			%		\draw [blue!80!black, ultra thick] 
				%	(v13) edge (v21)
				%	(v13) edge (v22)
				%	(v13) edge (v23)
				%	(v12) edge (v23)
				%	;
					
						\qedge{(v12)}{(v13)}{(w)}{5pt}{.8pt}{red!80!black}{red!10!white,opacity=0.2};
						\qedge{(v13)}{(w)}{(v)}{3.8pt}{.8pt}{red!80!black}{red!10!white,opacity=0.2};
					
			%		\draw [red!80!black, ultra thick] 
			%		(v13) edge (v)
			%		(v13) edge (v12)
			%		;
					
					\foreach \i in {1,2}{
						\foreach \j in {1,2,3}{
							\fill (v\i\j) circle (1.5pt);}}
					
					\fill [blue!80!black](v) circle (1.7pt);
					\fill [red!80!black](w) circle (1.7pt);
					
						\node [blue!70!black] at ($(v)+(.3,.1)$){$v$};
						\node  at ($(v11)+(-.5,.3)$) {$e_1$};
						\node  at ($(v21)+(.5,.3)$){$e_2$};
						
						\node [blue!70!black] at (2.6,1.45){$f$};
					
					\end{tikzpicture}

			\begin{tikzpicture}[scale = 1]

			\node at (.5,1.8) { $(d)$};
			
				\foreach \i in {1,2,3}{
					\foreach \j in {1,2,3}{
						\coordinate (v\i\j) at (1.2*\i,.5*\j);
					}}

			\coordinate (v) at (2.1,2);

			\qedge{(v11)}{(v12)}{(v13)}{7pt}{.3pt}{black}{black!20!white,opacity=0.2};
			\qedge{(v21)}{(v22)}{(v23)}{7pt}{.3pt}{black}{black!20!white,opacity=0.2};
			\qedge{(v31)}{(v32)}{(v33)}{7pt}{.3pt}{black}{black!20!white,opacity=0.2};
			
			\qedge{(v12)}{(v)}{(v23)}{2.5pt}{.8pt}{blue!80!black}{blue!10!white,opacity=0.2};
			\qedge{(v13)}{(v)}{(v23)}{3.7pt}{.8pt}{blue!80!black}{blue!10!white,opacity=0.2};
			
			\qedge{(v12)}{(v22)}{(v32)}{4pt}{.8pt}{green!70!black}{green!10!white,opacity=0.2};
			
	%		\draw [blue!80!black,very thick] 
	%%		(v13) edge (v22)
		%	(v13) edge (v23)
		%	(v12) edge (v23)
		%	;
			
			\foreach \i in {1,2,3}{
				\foreach \j in {1,2,3}{
					\fill (v\i\j) circle (1.5pt);}}
			
			\fill [blue!80!black](v) circle (1.7pt);

			\node [blue!70!black] at ($(v)+(.3,.1)$){$v$};
			\node  at ($(v11)+(-.44,.15)$) {\small $e_1$};
			\node  at ($(v21)+(-.44,.15)$){\small  $e_2$};
			\node  at ($(v31)+(-.44,.15)$){\small $e_3$};
			\node at ($(v13)+(-.066,-.22)$){\tiny $x_1$};
			\node at ($(v23)+(.06,-.25)$){\tiny $x_2$};
			
			%				\node [blue!70!black] at (2.2,.8){ $f$};
			\node [blue!70!black] at (1.95,1.65){$f$};
			\node [green!50!black] at ($(v22)!.5!(v32)$){\tiny $h'$};
			
			\end{tikzpicture}

			\end{multicols}
		
			\vspace{-1em}
			
		\caption{Possible neighbours of $w\in S$ in $\cF_{12}$ and $\cF_1$. The link graphs of $w$ and $v$ are denoted by red and blue 2-edges, respectively.}	
		%\caption{The illustrations to the proof of \eqref{eq:fg}.}
		\label{fig:f10}
				
	\end{figure}
	
%	\vspace{-1em}	

			For this purpose, recall that in view of \ref{it:n31} every vertex $w\in S$ can have positive degree only in one of the graphs $\cF_1, \cF_2$, say $\cF_1$. Next observe, that there exists an edge $f\in \ff$ disjoint from $\{x_1,w\}$, because only three out of at least four edges of $\ff$ containing $v$ can be incident to $x_1$ (see Figure \ref{fig:f10}$(b)$). Now, repeated application of \ref{it:n34} tells us that every edge $e\in \cF_1$ intersects both $h$ and $f$, and thereby contains $x_1$ and one of two vertices of $f\setminus e_2$. Clearly $w$ is contained in at most two of such edges (see Figure \ref{fig:f10}$(c)$).			
			
			Further we claim that
				\begin{equation}\label{eq:f3}
					|\fff|\le 3,
				\end{equation}
			because every edge $h'\in \fff$ contains both $x_1$ and $x_2$. Indeed, if not, let $h'=x_1'x_2'x_3'$ and say $x_2\neq x_2'$. Then, as in view of \ref{it:n34}, every edge of $\ff$ intersects both $h$ and $h'$, either $N_{\ff}(v)\subseteq \{x_1\}\times e_2$ if $x_1=x'_1$, or $N_{\ff}(v)\subseteq \{x_1'x_2, x_1x_2'\}$ otherwise, contradicting $\deg_{\ff}(v)\ge 4$.
			
			Now we are ready to finish the proof of Claim \ref{cl:out}. To this end denote $|V_3| = t$, and thereby $|S|=n-7-t$, as clearly $|R|=7$. Moreover, we let
				\[
					\cH_S = \{h\in \cH: h\cap S\neq\emptyset\},
				\]
			and observe that \ref{it:n31} entails
				\[
					\cH = \cF_3\dcup \fff \dcup \cH[R] \dcup \cH_S.
				\]
			Next note, that Lemma \ref{l:n7} combined with $\p$-freeness of $\cH$ tells us $|\cH[R]|\le 19$, and \ref{it:n34} yields $x_3\in f$ for each $f\in \cF_3$, causing $|\cF_3|\le {t-1\choose 2}$. Altogether, in view of \eqref{eq:deg7} and \eqref{eq:f3}, for $n\ge 8$,
				\[
					|\cH| =|\cF_3| + |\fff|+|\cH[R]| + |\cH_S| \le {t-1\choose 2} + 3 + 19 + 7(n-7-t) \le {n-4\choose 2} +10,
				\]
			as the left hand side of the last inequality achieves its maximum for either $t=3$ or $t=n-6$.
		\end{proof}
		
	Having established the above claims we proceed with the proof of \eqref{eq:lM3}. To this end, recall that \ref{it:n33} combined with \ref{it:n34} entail, that for each $i=1,2,3$, $\cF_i$ is a star, and thus $|\cF_i|\le {|V_i|-1\choose 2}$. Therefore, in view of \ref{it:n31}, by simple optimization,
		\begin{equation*}\label{eq:Fi}
			|\cF_1| + |\cF_2| + |\cF_3| \le {n-7 \choose 2} + 2,
		\end{equation*}
	with the equality achieved if and only if one of the 3-graphs $\cF_i$, $i=1,2,3$ is a full star on $n-6$ vertices, whereas two remaining 3-graphs each consists of a single edge. Further, we may assume that each vertex $v\in V\setminus (e_1\cup e_2 \cup e_3)$ satisfies $\deg_{\ff}(v)\le 3$, and thereby $|\ffo|\le 3(n-9)$, since otherwise Claim~\ref{cl:out} tells us that $|\cH|\le{n-4\choose 2} + 10$, and \eqref{eq:lM3} follows without the equality. Combining these observations together with \eqref{eq:Hsplit} and Claim~\ref{cl:in} one gets
		\[
			|\cH| =  (|\cF_1|+|\cF_2| +|\cF_3|) +(|\fff|+|\ffi|) +|\ffo|
				\le
			{n-7 \choose 2} + 2 + 18 + 3(n-9) = {n-4\choose 2} + 11,
		\]
	as required. It remains to show that the equality in the above bound is achieved if and only if $\cH$ is a balloon.
		
	Indeed, clearly if $|\cH|={n-4\choose 2} + 11$, then equalities in the above formula go through.
	In particular, if $n=9$, then $\ffo=\emptyset$, $\cF_i$, $i=1,2,3$, is a single edge and $\fff\cup \ffi$ is a star with center in $e_1\cup e_2$. It is easy to see that $\cH=B_9$.
	
	Now assume $n\ge 10$, $|\ffo|=3(n-9)$ and thus $|V[\ff]| = n-3$ yielding, in view of \ref{it:n31}, $|V_3|=3$. Therefore, as $|\cF_1| + |\cF_2| + |\cF_3| = {n-7 \choose 2} + 2$, without loss of generality we may assume that $\cF_1$ is a full star on $n-6\ge 4$ vertices, whereas $|\cF_2|=|\cF_3|=1$.
	Let $c\in e_1$ be the center of $\cF_1$.
	As the pair $(\cF_1, \ff\cup \fff)$ is cross-intersecting and $\fff\cup \ffi$ is a star, the center of $\fff\cup \ffi$ must also be $c$ as for any vertex $v\neq c$, the \emph{full} star $\cF_1$ contains an edge not containing $v$. Finally, since $\fff$ contains all possible 9 edges containing $c$ and the pair $(\fff,\ff)$ is cross-intersecting, every edge of $\ffo$ contains $c$ as well. Altogether $\cF_1\cup \ff\cup \fff$ is a star, whereas $\cF_2$ and $\cF_3$ are single edges, and thereby $\cH=B_n$.

%%%%%%%%%%%%%%%%%%%%%%%%%%%%%%%%%%%%%%%%%%%%%%%%%%%%%%%%%%%%%%%%%%%%%%%%%%%
%																												   									%
%																												  			 						%
%  									chapter 6 														    				   						%
%																																					%
%																																					%
%%%%%%%%%%%%%%%%%%%%%%%%%%%%%%%%%%%%%%%%%%%%%%%%%%%%%%%%%%%%%%%%%%%%%%%%%%%

\section{Proofs of Lemmas \ref{l:c41} and \ref{l:c42}}

 Let $\cH$ be a connected $\p\cup\{\m\}$-free 3-graph on the set of $n$ vertices $V$, $n\ge 8$, such that $\c\subseteq\cH$. Denote by
				\[
					\cc=\{x_1y_1y_2,x_1z_1z_2,x_2y_1y_2,x_2z_1z_2\}
				\]
a copy of $\c$ contained in $\cH$, and set $V[\cc]=Z=\{x_1,x_2,y_1,y_2,z_1,z_2\}$, $W=V\setminus Z$.
			
	Lemmas \ref{l:c41} and \ref{l:c42} are straightforward consequences of the following two lemmas.
			
			\begin{lemma}\label{l:sp}
				If there exist two vertices $u, w\in W$ with degree in $\cH$ at least 5, and moreover either
					\begin{enumerate}[label=\rmlabel]
						\item $|\cH[Z\cup \{u,w\}]|\ge 22$ or
						\item there is a further vertex $v\in W\setminus \{u,w\}$ with $\deg_\cH(v)\ge 5$,
					\end{enumerate}
			then $\cH\subseteq \sp_n$.		
			\end{lemma}
			
			\begin{lemma}\label{l:sk}
				If there exist two vertices $u,w\in W$, such that
					\begin{enumerate}[label = \rmlabel]
						\item $|\cH[Z\cup \{u,w\}]|\ge 22$ and
						\item $\deg_\cH(w) \le 4$,
					\end{enumerate}
				then $\cH\subseteq \sk_n$.
			\end{lemma}
			
		Indeed, assume first that there are two vertices $u,w\in W$, such that $|\cH\cup \{u,w\}|\ge 22$. Then, either $\deg_\cH(u), \deg_\cH(w)\ge 5$ and thus, in view of Lemma \ref{l:sp}, $\cH\subseteq \sp_n$, or the degree in $\cH$ of one of these vertices, say $w$, is at most 4. Then Lemma \ref{l:sk} tells us that $\cH\subseteq \sk_n$.
		
		So let for every pair of vertices $u,w\in W$, $|\cH[Z\cup \{u,w\}]|\le 21$. If there are three vertices $u,w,v\in W$, with the degree in $\cH$ at least 5, then due to Lemma \ref{l:sp}, $\cH\subseteq \sp_n$. Otherwise choose $u,w\in W$ in such a way, that for all $v\in W\setminus \{u,w\}$, $\deg_\cH(v)\le 4$. Then,
			\[
				|\cH| = |\cH[Z\cup \{u,w\}]|+\sum_{v\in W\setminus \{u,w\}}\deg_\cH(v)\le 4n-11.
			\]

		Altogether, either $\cH\subseteq \sp_n$, $\cH\subseteq \sk_n$, or $|\cH|\le 4n-11$. Now, as $\sk_n\nsubseteq \sp_n$, $|\sp_n|=5n-18$, and $|\sk_n|=4n-10$, Lemmas \ref{l:c41} and \ref{l:c42} follows from
			\[
				5n-18\ge 4n-10 > 4n-11,
			\]
		for $n\ge 8$, with the equality only for $n=8$.
				
	\subsection{Preliminaries}
	We begin with a series of technical results which  will be helpful in the proofs of Lemmas \ref{l:sp} and \ref{l:sk}. Throughout, we denote by $u$ and $w$ arbitrary vertices of $W$. The $\p\cup\{\m\}$-freeness of $\cH$ implies that for all  edges $h\in \cH$
		\begin{equation}\label{hZ}
			|h\cap Z|\ge2\;,\qquad h\cap Z\neq\{y_1,y_2\}\;,\qquad h\cap Z\neq\{z_1,z_2\}.
		\end{equation}
	Let us partition $\cH$ into four edge-disjoint sub-3-graphs,
		\[
			\cH=\cH_Z\cup \cH^0\cup \cH^1\cup \cH^2,
		\]
	where $\cH_Z=\cH[Z]$ and, for $i=0,1,2$,
		\[
			\cH^i=\{h\in \cH\setminus \cH_Z: |h\cap \{x_1,x_2\}|=i\}.
		\]
	The first inequality in (\ref{hZ}) implies that the link graph $L_{\cH}(w)$ of every vertex $w\in W$ is entirely contained in $\binom Z2$. Moreover, the above partition of $\cH$ induces a corresponding partition of $L_{\cH}(w)$,
			\[
				L_{\cH}(w)=H^0(w)\cup H^1(w)\cup H^2(w),
			\]
	where $H^i(w)=L_{\cH^i}(w)=\{e\in L_\cH(w): |e\cap \{x_1,x_2\}|=i\}$. Observe, that
			\begin{equation}\label{eq:h02}
				|H^0(w)|\le 4\;, \quad |H^1(w)|\le 8\;,\textrm{ and } \quad |H^2(w)|\le 1,
			\end{equation}
	where the first inequality holds, because in view of \eqref{hZ}, $\{y_1,y_2\}, \{z_1,z_2\}\notin L_{\cH}(w)$, and thus $H^0(w)$ is a subgraph of the 4-cycle $y_1z_1y_2z_2$.
						
			\medskip
			
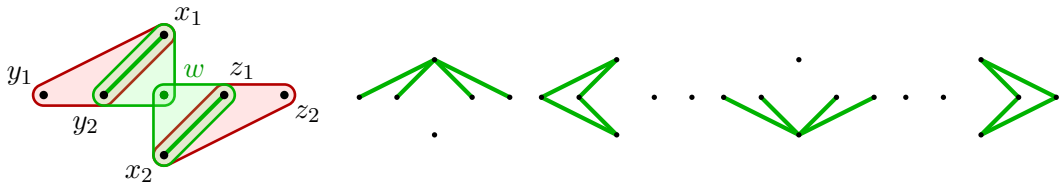
\begin{figure}[h!]
				\begin{tabular}{c c c c c}
					\begin{tikzpicture}[scale=.8]
					\coordinate (x2) at (0,0);
					\coordinate (x1) at (0,2);	
					\coordinate (y1) at (-2,1);
					\coordinate (y2) at (-1,1);
					\coordinate (z1) at (1,1);
					\coordinate (z2) at (2,1);
					\coordinate (w) at (0,1);
					\qedge{(x2)}{(z1)}{(z2)}{5pt}{1pt}{red!70!black}{red!50!white,opacity=0.2};
					\qedge{(x1)}{(y2)}{(y1)}{5pt}{1pt}{red!70!black}{red!50!white,opacity=0.2};
					\qedge{(w)}{(y2)}{(x1)}{5pt}{1pt}{green!70!black}{green!50!white,opacity=0.2};
					\qedge{(w)}{(z1)}{(x2)}{5pt}{1pt}{green!70!black}{green!50!white,opacity=0.2};
					\draw (x1) edge [color=green!70!black, ultra thick] (y2);
					\draw (x2) edge [color=green!70!black, ultra thick] (z1);
					\foreach \i in {x1,x2,y1,y2,z1,z2,w}
					\fill  (\i) circle (2pt);
					\fill [green!60!black] (w) circle (2pt);
					\node at (x1) [above right] {$x_1$};
					\node at (-1.3,.5) {$y_2$};
					\node at (x2) [below left] {$x_2$};
					\node at (1.3,1.4) {$z_1$};
					\node at (y1) [above left ] {$y_1$};
					\node  at (z2)[below right]{$z_2$};
					\node at (.5,1.4)  [color=green!60!black] {$w$};
					\end{tikzpicture}
					
					&\begin{tikzpicture}[scale=.5]
					\coordinate (x2) at (0,0);
					\coordinate (x1) at (0,2);
					
					\coordinate (y1) at (-2,1);
					\coordinate (y2) at (-1,1);
					
					\coordinate (z1) at (1,1);
					\coordinate (z2) at (2,1);
					
					\coordinate (w) at (0,1);

					\draw [white] (0,-1.5) rectangle (0,-1.5);
					\draw (x1) edge [color=green!70!black, ultra thick] (y1);
					\draw (x1) edge [color=green!70!black, ultra thick] (y2);
					\draw (x1) edge [color=green!70!black, ultra thick] (z1);
					\draw (x1) edge [color=green!70!black, ultra thick] (z2);
					
					\foreach \i in {x1,x2,y1,y2,z1,z2}
					\fill  (\i) circle (2pt);
					\end{tikzpicture}
					
					&\begin{tikzpicture}[scale=.5]
					
					\draw [white] (0,-1.5) rectangle (0,-1.5);
					\coordinate (x2) at (0,0);
					\coordinate (x1) at (0,2);
					
					\coordinate (y1) at (-2,1);
					\coordinate (y2) at (-1,1);
					
					\coordinate (z1) at (1,1);
					\coordinate (z2) at (2,1);
					
					\coordinate (w) at (0,1);
					
					\draw (x1) edge [color=green!70!black, ultra thick] (y1);
					\draw (x1) edge [color=green!70!black, ultra thick] (y2);
					
					\draw (x2) edge [color=green!70!black, ultra thick] (y1);
					\draw (x2) edge [color=green!70!black, ultra thick] (y2);
					
					\foreach \i in {x1,x2,y1,y2,z1,z2}
					\fill  (\i) circle (2pt);
					
					\end{tikzpicture}

					&\begin{tikzpicture}[scale=.5]
					\draw [white] (0,-1.5) rectangle (0,-1.5);
					\coordinate (x2) at (0,0);
					\coordinate (x1) at (0,2);
					
					\coordinate (y1) at (-2,1);
					\coordinate (y2) at (-1,1);
					
					\coordinate (z1) at (1,1);
					\coordinate (z2) at (2,1);
					
					\coordinate (w) at (0,1);
					
					\draw (x2) edge [color=green!70!black, ultra thick] (y1);
					\draw (x2) edge [color=green!70!black, ultra thick] (y2);
					\draw (x2) edge [color=green!70!black, ultra thick] (z1);
					\draw (x2) edge [color=green!70!black, ultra thick] (z2);
					
					\foreach \i in {x1,x2,y1,y2,z1,z2}
					\fill  (\i) circle (2pt);
					\end{tikzpicture}
					
					&\begin{tikzpicture}[scale=.5]
					\draw [white] (0,-1.5) rectangle (0,-1.5);
					\coordinate (x2) at (0,0);
					\coordinate (x1) at (0,2);
					
					\coordinate (y1) at (-2,1);
					\coordinate (y2) at (-1,1);
					
					\coordinate (z1) at (1,1);
					\coordinate (z2) at (2,1);
					
					\coordinate (w) at (0,1);
					
					\draw (x1) edge [color=green!70!black, ultra thick] (z1);
					\draw (x1) edge [color=green!70!black, ultra thick] (z2);
					
					\draw (x2) edge [color=green!70!black, ultra thick] (z1);
					\draw (x2) edge [color=green!70!black, ultra thick] (z2);
					
					\foreach \i in {x1,x2,y1,y2,z1,z2}
					\fill  (\i) circle (2pt);
					
					\end{tikzpicture}
				\end{tabular}
				
				\caption{A minimal 4-path $y_1y_2x_1wx_2z_1z_2$ in $\cH$ and all possible edges of link graphs $H^1(w)$.}
				\label{fig:F0}
				\vspace{-1em}
			\end{figure}

	Our first result describes the structure of $H^1(w)$ and, as a consequence, halves the upper bound on $|H^1(w)|$.
			\begin{fact}\label{f0}
				For every $w\in W$, $H^1(w)$ is either a star (with the center at $x_1$ or $x_2$) or a subgraph of one of the 4-cycles: $C_y=x_1y_1x_2y_2$ or  $C_z=x_1z_1x_2z_2$. In particular, $|H^1(w)|\le 4$.
			\end{fact}
			\begin{proof}
			If there were  two disjoint edges in $H^1(w)$, one contained in $C_y$ and the other  in $C_z$, say $\{x_1,y_2\}$ and $\{x_2,z_1\}$, then $y_1y_2x_1wx_2z_1z_2$ would form a minimal 4-path  in $\cH$, a contradiction (see Figure \ref{fig:F0}). So, either all edges of $H^1(w)$ are contained in one of the cycles, $C_y$ or $C_z$, or they form a star.
			\end{proof}
			
		It is convenient to break the 3-graph $\cH_Z$  into three further sub-3-graphs,
			\[
				\cH_Z=\cH_Z^0\cup \cH_Z^1\cup \cH_Z^2,\quad\textrm{where}\quad
				 \cH^i_Z= \{h\in \cH_Z: |h\cap \{x_1,x_2\}|=i\}, \quad i=0,1,2.
			\]
		Note that $C\subseteq \cH_Z^1$, $|\cH_Z^0|\le \binom43=4$, $|\cH_Z^1|\le2\binom42= 12$, and $|\cH_Z^2|\le\binom 41= 4$.
			
			The next result lists several  basic observations on the above defined subgraphs, all stemming from the $\p$-freeness of $\cH$. 		
			
			\begin{fact}\label{f00}
				Let $h,h'\in C$ be disjoint and let, for some vertex $w\in W$, an edge $e\in H^1(w)$ be contained in $h$. Then there is no edge $f\in \cH-w$ disjoint from $e$ and intersecting both $h$ and $h'$. Consequently, for any two distinct vertices $u,w\in W$, the following properties hold,
				\begin{enumerate}[label=\rmlabel]
					\item \label{it:f01} if $e\in H^1(u)$ and $e'\in H^1(w)$ are disjoint, then there exist  disjoint  $h,h'\in C$ such that $e\subset h$ and $e'\subset h'$;
					\item\label{it:f02} the pair of 2-graphs $(H^0(u),H^1(w))$ is cross-intersecting;
					\item\label{it:f025} if $e\in H^1(w)$ and $f\in \cH^0_Z \cup \cH^1_Z$, $f\not\in C$, then $e\cap f\neq\emptyset$;
					\item\label{it:f03} if $H^1(w)\neq \emptyset$ then $|\cH_Z^1|\le 10$;
					\item\label{it:f04} if $|H^1(u)\cup H^1(w)| \ge 2$ then $|\cH_Z^1|\le 9$ and $|\cH_Z|\le 15$;
					\item \label{it:f05} if $|H^1(u)\cup H^1(w)|\ge 3$ then $|\cH_Z^1|\le 8$ and $|\cH_Z|\le 13$;
					\item \label{it:f06} if $|H^1(w)|= 4$ then $|\cH_Z|\le 12$;
					\item \label{it:f066} if $|H^1(u)\cup H^1(w)|\ge 7$ then $|\cH_Z|\le 8$;
					\item \label{it:f07} if $H^1(w)$ is a star with four edges and the center $x_1$ or $x_2$, then $\cH\subseteq \sp_n$. %(see Figure \ref{fig:F00} (a)--(d)).
				\end{enumerate}
			\end{fact}
			
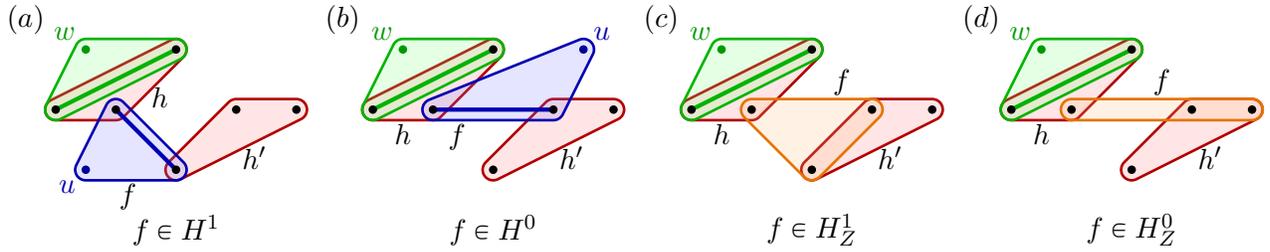
\begin{figure}[h!]
				\begin{multicols}{4}
					\begin{tikzpicture}[scale=.8]
					\coordinate (x2) at (0,0);
					\coordinate (x1) at (0,2);	
					\coordinate (y1) at (-2,1);
					\coordinate (y2) at (-1,1);
					\coordinate (z1) at (1,1);
					\coordinate (z2) at (2,1);
					\coordinate (w) at (-1.5,0);
					\coordinate (u) at (-1.5,2);
					\qedge{(x2)}{(z1)}{(z2)}{5pt}{1pt}{red!70!black}{red!50!white,opacity=0.2};
					\qedge{(x1)}{(y2)}{(y1)}{5pt}{1pt}{red!70!black}{red!50!white,opacity=0.2};
					\qedge{(w)}{(y2)}{(x2)}{5pt}{1pt}{blue!70!black}{blue!50!white,opacity=0.2};
					\qedge{(u)}{(x1)}{(y1)}{5pt}{1pt}{green!70!black}{green!50!white,opacity=0.2};
					\draw (x2) edge [color=blue!70!black, ultra thick] (y2);
					\draw (x1) edge [color=green!70!black, ultra thick] (y1);
					\foreach \i in {x1,x2,y1,y2,z1,z2}
					\fill  (\i) circle (2pt);
					\fill [blue!60!black] (w) circle (2pt);
					\fill [green!60!black] (u) circle (2pt);
					\node at (u) [color=green!60!black, above left] {$w$};
					\node at (1.3,0.2) {$h'$};
					\node at (-0.27,1.2) {$h$};
					\node at (-.8,-0.45) {$f$};
					\node at (w) [color=blue!70!black,below left] {$u$};
					\node at (0,-1){$f\in \cH^1$};
					\node at (-2.5,2.5) {$(a)$};
					\end{tikzpicture}\\
									
					\begin{tikzpicture}[scale=.8]
					\node at (1.3,0.2) {$h'$};
					\node at (-1.5,0.6) {$h$};
					\node at (-0.6,0.55) {$f$};
					\coordinate (x2) at (0,0);
					\coordinate (x1) at (0,2);	
					\coordinate (y1) at (-2,1);
					\coordinate (y2) at (-1,1);
					\coordinate (z1) at (1,1);
					\coordinate (z2) at (2,1);
					\coordinate (w) at (1.5,2);
					\coordinate (u) at (-1.5,2);
					\qedge{(x2)}{(z1)}{(z2)}{5pt}{1pt}{red!70!black}{red!50!white,opacity=0.2};
					\qedge{(x1)}{(y2)}{(y1)}{5pt}{1pt}{red!70!black}{red!50!white,opacity=0.2};
					\qedge{(w)}{(z1)}{(y2)}{5pt}{1pt}{blue!70!black}{blue!50!white,opacity=0.2};
					\qedge{(u)}{(x1)}{(y1)}{5pt}{1pt}{green!70!black}{green!50!white,opacity=0.2};
					\draw (y2) edge [color=blue!70!black, ultra thick] (z1);
					\draw (x1) edge [color=green!70!black, ultra thick] (y1);
					\foreach \i in {x1,x2,y1,y2,z1,z2}
					\fill  (\i) circle (2pt);
					\fill [blue!60!black] (w) circle (2pt);
					\fill [green!60!black] (u) circle (2pt);
					\node at (u) [color=green!60!black,above left] {$w$};
					\node at (w) [color=blue!70!black, above right] {$u$};
					\node at (0,-1){$f\in \cH^0$};
					\node at (-2.5,2.5) {$(b)$};
					\end{tikzpicture}\\
					
					\begin{tikzpicture}[scale=.8]
					\node at (1.3,0.2) {$h'$};
					\node at (-1.5,0.6) {$h$};
					\node at (0.5,1.45) {$f$};
					\coordinate (x2) at (0,0);
					\coordinate (x1) at (0,2);	
					\coordinate (y1) at (-2,1);
					\coordinate (y2) at (-1,1);
					\coordinate (z1) at (1,1);
					\coordinate (z2) at (2,1);
					\coordinate (w) at (-1.5,0);
					\coordinate (u) at (-1.5,2);
					\qedge{(x2)}{(z1)}{(z2)}{5pt}{1pt}{red!70!black}{red!50!white,opacity=0.2};
					\qedge{(x1)}{(y2)}{(y1)}{5pt}{1pt}{red!70!black}{red!50!white,opacity=0.2};
					\qedge{(z1)}{(x2)}{(y2)}{5pt}{1pt}{orange!90!black}{orange!50!white,opacity=0.2};
					\qedge{(u)}{(x1)}{(y1)}{5pt}{1pt}{green!70!black}{green!50!white,opacity=0.2};
					\draw (x1) edge [color=green!70!black, ultra thick] (y1);
					\foreach \i in {x1,x2,y1,y2,z1,z2}
					\fill  (\i) circle (2pt);
					\fill [green!60!black] (u) circle (2pt);
					\node at (u) [color=green!60!black, above left] {$w$};
					\node at (0,-1){$f\in \cH_Z^1$};
					\node at (-2.5,2.5) {$(c)$};
					\end{tikzpicture}\\
					
					\begin{tikzpicture}[scale=.8]
					\node at (1.3,0.2) {$h'$};
					\node at (-1.5,0.6) {$h$};
					\node at (0.5,1.45) {$f$};
					\coordinate (x2) at (0,0);
					\coordinate (x1) at (0,2);	
					\coordinate (y1) at (-2,1);
					\coordinate (y2) at (-1,1);
					\coordinate (z1) at (1,1);
					\coordinate (z2) at (2,1);
					\coordinate (w) at (-1.5,0);
					\coordinate (u) at (-1.5,2);
					\qedge{(x2)}{(z1)}{(z2)}{5pt}{1pt}{red!70!black}{red!50!white,opacity=0.2};
					\qedge{(x1)}{(y2)}{(y1)}{5pt}{1pt}{red!70!black}{red!50!white,opacity=0.2};
					\qedge{(y2)}{(z1)}{(z2)}{5pt}{1pt}{orange!90!black}{orange!50!white,opacity=0.2};
					\qedge{(u)}{(x1)}{(y1)}{5pt}{1pt}{green!70!black}{green!50!white,opacity=0.2};
					\draw (x1) edge [color=green!70!black, ultra thick] (y1);
					\foreach \i in {x1,x2,y1,y2,z1,z2}
					\fill  (\i) circle (2pt);
					\fill [green!60!black] (u) circle (2pt);
					\node at (u) [color=green!60!black, above left] {$w$};
					\node at (0,-1){$f\in \cH^0_Z$};
					\node at (-2.5,2.5) {$(d)$};
					\end{tikzpicture}
				\end{multicols}
				\caption{Illustration to the proof of Fact \ref{f00}.}
				\label{fig:F00}
				\vspace{-1em}
			\end{figure}

		\begin{proof}
			Suppose that $h,h'\in C$, $w\in W$, $e\in H^1(w)$, and $f\in \cH-w$ are such that $h\cap h'=\emptyset$, $e\subset h$, and  $f\cap h=h\setminus e$ and $f\cap h'\neq\emptyset$. Then, regardless of the location of $f$, the 3-edges $we, h, f$, and $h'$ form a minimal 4-path in $\cH$ (see Figure \ref{fig:F00}), contradicting the $\p$-freeness of $\cH$. So the main statement is proved, and consequently, \ref{it:f01}-\ref{it:f025} follow. Indeed, if \ref{it:f01}, \ref{it:f02}, or \ref{it:f025} were not true, then we would be looking at the forbidden configurations in Figure \ref{fig:F00}$(a)$, $(b)$ or $(c),(d)$, respectively.
			
			In turn, \ref{it:f025} implies \ref{it:f03}-\ref{it:f066}. Indeed, \ref{it:f03} follows from the bound $|\cH_Z^1|\le12$ as, in view of \ref{it:f025}, $H^1(w)\neq\emptyset$ excludes two edges from $\cH_Z^1$. Similarly, in \ref{it:f04}, considering five different cases with respect to the location of the two edges of $H^1(u)\cup H^1(w)$, we may exclude (by applying \ref{it:f025}) at least 3 edges of $\cH_Z^1$ and at least 5 edges of $\cH_Z^0\cup \cH_Z^1$. By the same token, in \ref{it:f05}, we exclude at least 4 edges of $\cH_Z^1$ and at least 7 edges of $\cH_Z^0\cup \cH_Z^1$. For the proof of \ref{it:f06}, recall that Fact \ref{f0} tells us that $H^1(w)$ is either a 4-arm star or one of the cycles $C_y$ or $C_z$. In both cases, via \ref{it:f025}, it wipes out at least 4 edges of $\cH^1_Z$ and  at least 8 edges of $\cH^0_Z\cup \cH^1_Z$. We leave case \ref{it:f066} for the Reader.
			
			Finally, to prove \ref{it:f07}, assume, without loss of generality, that $H^1(w)$ is a 4-edge star with the center $x_1$. Now observe that by \ref{it:f01}-\ref{it:f025} every edge of $\cH$, except for $x_2y_1y_2$ and $x_2z_1z_2$ (which form $P=P^{(3)}_2$ disjoint from $\{x_1\}$), contains both $x_1$ and a member of $V[P]=\{y_1,y_2,x_2,z_1,z_2\}$, entailing $\cH\subseteq \sp_n$. Indeed, \ref{it:f02} yields that $H^0(u)=\emptyset$ and $H^1(u)\subseteq H^1(w)$ holds by~\ref{it:f01}. Moreover \ref{it:f025} tells us that $H^0_Z=\emptyset$ and whenever $f\in H^1_Z\setminus C$, $x_1\in f$.			
		\end{proof}
			
			\begin{cor}\label{f2}
				For all distinct $u,w\in W$,
				\begin{enumerate}[label=\rmlabel]
					\item\label{it:f22} $H^1(u)\neq\emptyset \Rightarrow |H^0(w)|\le 2$;
					\item\label{it:f21} $H^0(u)\neq\emptyset \Rightarrow |H^1(w)|\le 2$.
				\end{enumerate}
			\end{cor}
			 \begin{proof}
			 	Observe that for any edge $e\in H^1(u)$, there exist at most two edges in $H^0(w)$ which intersect $e$. Similarly, by Fact \ref{f0}, for any edge $e\in H^0(u)$ there are at most two edges in $H^1(w)$ sharing a vertex with $e$. Consequently, by Fact \ref{f00}\ref{it:f02}, both assertions follow.
			 \end{proof}
			
			\begin{cor}\label{minmax}
				If there is a vertex $u\in W$ with $\deg_\cH(u)\ge6$, then the degree of every vertex $w\in W\setminus\{u\}$ is at most 3. Moreover, if additionally $\deg_\cH(w)=3$, then $|H^2(w)|=1$, that is, $wx_1x_2\in \cH$.
			\end{cor}
			 \begin{proof}
			 	Let $\deg_\cH(u)\ge6$ and let $w\in W\setminus \{u\}$. By \eqref{eq:h02} and Fact \ref{f0}, both sets, $H^0(u)$ and $H^1(u)$, must be nonempty  and at least one of them of size at least three, say $|H^1(u)|\ge3$. But then,  by Corollary \ref{f2}\ref{it:f21}, $|H^1(w)|\le2$ and $H^0(w)=\emptyset$. Hence,
				 	\[
					 	\deg_\cH(w)=|L_{\cH}(w)|=|H^0(w)|+|H^1(w)|+|H^2(w)|\le 0+2+1=3
				 	\]
				and if $\deg_\cH(w)=3$, then $|H^2(w)|=1$.
			 \end{proof}
			
			\begin{fact}\label{f8}
				If $ux_1x_2,wx_1x_2\in \cH^2$ and  $e\in H^0(w)$, then there is no  $f\in \cH[V\setminus \{x_1,x_2,u,w\}]$ with $f\cap e\neq \emptyset$. It follows that $\cH^0_Z=\emptyset$. Moreover, if $|H^0(w)|\ge 2$, then for every $v\in W\setminus \{u,w\}$, we have $H^0(v)=\emptyset$.
			\end{fact}
			\begin{proof}
				To prove the first statement, it is enough to observe that whenever $ux_1x_2,x_1x_2w\in \cH$, $e\in H^0(w)$, and $f\in \cH[V\setminus \{x_1,x_2,u,w\}]$, $f\cap e\neq \emptyset$, then edges $ux_1x_2$, $x_1x_2w$, $we$, and $f$ form a minimal 4-path in $\cH$. As $e$ uses two of the four vertices of $V[\cH^0_Z]$, there is no room for an $f\in \cH^0_Z$ with $f\cap e=\emptyset$, and so $\cH_Z^0=\emptyset$. Furthermore, if $|H^0(w)|\ge 2$ and  $e'\in H^0(v)$ then $f= e' v \in \cH[V\setminus \{x_1,x_2,u,w\}]$ and there exists $e\in H^0(w)$ such that $f\cap e=e'\cap e\neq\emptyset$, a contradiction.
			\end{proof}
						
\begin{figure}[h!]
				\centering
				\begin{multicols}{4}
					\begin{tikzpicture}[scale=.7]
				
					\coordinate (x2) at (0,0);
				\coordinate (x1) at (0,2);
					
					\coordinate (y1) at (-2,1);
					\coordinate (y2) at (-1,1);
					
					\coordinate (z1) at (1,1);
					\coordinate (z2) at (2,1);
					
				\coordinate (u) at (0,1.5);
					\coordinate (w) at (.5,1);

					\qedge{(y1)}{(y2)}{(z1)}{6pt}{1pt}{red!80!black}{red!50!white,opacity=0.2};		
					\qqedge{(x2)}{($(y1)+(-.3,0)$)}{(x1)}{5pt}{1pt}{red!80!black}{red!50!white,opacity=0.2}{16pt};
					\draw (y2) edge [color=green!70!black, ultra thick] (z1);
					\draw (x1) edge [color=blue!70!black, ultra thick] (x2);
					
					\foreach \i in {x1,x2,y1,y2,z1,z2}
					\fill  (\i) circle (2pt);
					\fill (u)  [blue!70!black] circle (3pt);
					\fill [green!70!black] (w) circle (3pt);

					\node at ($(w) + (0,-.4)$)[green!50!black] {$w$};
					\node at ($(u)+(.3,0)$) [blue!50!black]{$u$};
					\node at (x2) [below right]{$x_2$};
					\node at (x1)[above right] {$x_1$};
					
					\node at ($(y1)+(-.45,-.45)$) {$y_1$};
				
					\node at ($(z1)+(.4,-.4)$){$z_1$};
					\node at ($(z2)+(.3,-.3)$) {$z_2$};

				\node at (0,-1) {{\tiny $\cbb{f_u},\{x_1,x_2,y_1\},\{y_1,y_2,z_1\}, \cgg{f_w}$}};
					
					\end{tikzpicture}\\
					
					\begin{tikzpicture}[scale=.7]
							
				\coordinate (x2) at (0,0);
					\coordinate (x1) at (0,2);					
					\coordinate (y1) at (-2,.7);
					\coordinate (y2) at (-1,1.3);
					
					\coordinate (w) at (-.3,1.1);
					
					\coordinate (z1) at (1,.7);
					\coordinate (z2) at (2,.7);
					
					\qedge{(y1)}{(z1)}{(z2)}{5pt}{1pt}{red!80!black}{red!50!white,opacity=0.2};		
					\qqedge{(x2)}{($(y2)+(-.1,0)$)}{(x1)}{5pt}{1pt}{red!80!black}{red!50!white,opacity=0.2}{9pt};
					\draw (y2) edge [color=green!70!black, ultra thick] (z1);
					\draw (x1) edge [color=blue!70!black, ultra thick] (x2);
					
					\foreach \i in {x1,x2,y1,y2,z1,z2}
					\fill  (\i) circle (2pt);
				
					\node at ($(u)+(.3,0)$) [blue!50!black]{$u$};
				
					\node at (x2)[below right] {$x_2$};
					\node at (x1)[above right] {$x_1$};
					
					\node at (y1)[below left] {$y_1$};
					\node at (y2)[above left] {$y_2$};
					
					\node at ($(z1)+(.2,.4)$) {$z_1$};
					\node at (z2)[below right] {$z_2$};
					
					\fill (u)  [blue!70!black] circle (3pt);
					\fill [green!70!black] (w) circle (3pt);
					
					\node at (0,-1) {{\tiny $\cbb{f_u},\{x_1,x_2,y_2\},\cgg{f_w},\{y_1,z_1,z_2\}$}};
					\node at ($(w)+(0,.2)$) [green!50!black]{{\tiny $w$}};
					
					\end{tikzpicture}\\
					
					\begin{tikzpicture}[scale=.7]	
					\coordinate (y2) at (-1,1);
					
					\coordinate (z1) at (1,1.3);
					\coordinate (w) at (.4,1.22);
					
					\qedge{(y1)}{(y2)}{(z2)}{5pt}{1pt}{red!80!black}{red!50!white,opacity=0.2};	
					\qqedge{(x1)}{($(z1)+(.1,0)$)}{(x2)}{5pt}{1pt}{red!80!black}{red!50!white,opacity=0.2}{9pt};	
					\draw (y2) edge [color=green!70!black, ultra thick] (z1);
					\draw (x1) edge [color=blue!70!black, ultra thick] (x2);
					
					\foreach \i in {x1,x2,y1,y2,z1,z2}
					\fill  (\i) circle (2pt);
					
					\fill (u)  [blue!70!black] circle (3pt);
					\fill [green!70!black] (w) circle (3pt);

					\node at ($(u)+(-.3,0)$) [blue!50!black]{$u$};
					
					\node at (x2)[below left] {$x_2$};
					\node at (x1)[above left] {$x_1$};
					
					\node at (y1)[below left] {$y_1$};
					\node at (y2)[above left] {$y_2$};
					
					\node at (z1)[above right] {$z_1$};
					\node at (z2)[below right] {$z_2$};
					
					\node at ($(w)+(-.1,.17)$) [green!50!black]{{\tiny $w$}};
					
					\node at (0,-1) {{\tiny $\cbb{f_u},\{x_1,x_2,z_1\},\cgg{f_w},\{y_1,y_2,z_2\}$}};

					\end{tikzpicture}\\
					
					\begin{tikzpicture}[scale=.7]	
					\coordinate (x2) at (0,0);
					\coordinate (x1) at (0,2);
					
					\coordinate (y1) at (-2,1);
					\coordinate (y2) at (-1,1);
					
					\coordinate (z1) at (1,1);
					\coordinate (z2) at (2,1);
					
					\coordinate (u) at (0,1.5);
					\coordinate (w) at (-.5,1);
					
					\qedge{(z2)}{(y2)}{(z1)}{6pt}{1pt}{red!80!black}{red!50!white,opacity=0.2};		
					\qqedge{(x1)}{($(z2)+(.3,0)$)}{(x2)}{5pt}{1pt}{red!80!black}{red!50!white,opacity=0.2}{16pt};
					
					\draw (y2) edge [color=green!70!black, ultra thick] (z1);
					\draw (x1) edge [color=blue!70!black, ultra thick] (x2);
					
					\foreach \i in {x1,x2,y1,y2,z1,z2}
				\fill  (\i) circle (2pt);
					\fill (u)  [blue!70!black] circle (3pt);
					\fill [green!70!black] (w) circle (3pt);
					
					\node at ($(w) + (0,-.4)$)[green!50!black] {$w$};
					\node at ($(u)+(-.3,0)$) [blue!50!black]{$u$};
					
					\node at (x2) [below left]{$x_2$};
					\node at (x1)[above left] {$x_1$};
					
					\node at ($(y1)+(-.25,-.25)$) {$y_1$};
					\node at ($(y2)+(-.35,-.35)$) {$y_2$};

					\node at ($(z2)+(.45,-.45)$) {$z_2$};
					
					\node at (0,-1) {{\tiny $\cbb{f_u},\{x_1,x_2,z_2\},\{y_2,z_1,z_2\},\cgg{f_w}$}};
					
					\end{tikzpicture}
					
				\end{multicols}
				
				\caption{A minimal 4-paths with edges \cbb{$f_u$ (blue)} and \cgg{$f_w$ (green)}.}
				\label{fig:f02}
				\vspace{-1em}
			\end{figure}
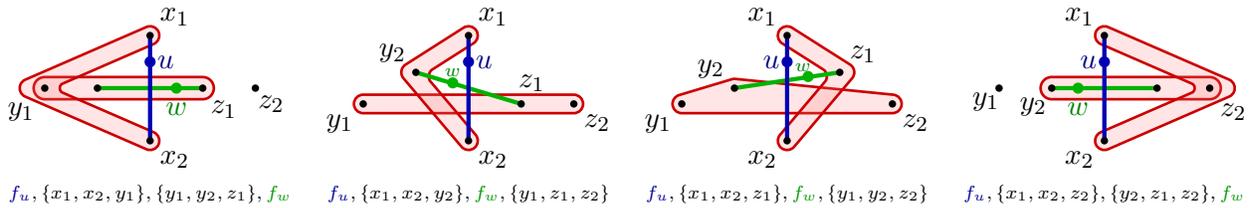

			\begin{fact}\label{f7}
				If $ux_1x_2\in \cH^2$ and $H^0(w)\neq\emptyset$, then $|\cH^0_Z|+|\cH^2_Z|\le 4$.
			\end{fact}
			\begin{proof}
				Let $f_u=ux_1x_2\in \cH^2$ and $f_w\in \cH^0$, $w\in f_w$. Without loss of generality we may assume that $f_w = wy_2z_1$. Observe that $\cH^0_Z\cup \cH_Z^2$ can be partitioned into four  pairs of edges,
					\[
						\{x_1x_2y_1,y_1y_2z_1\},\quad\{x_1x_2y_2,y_1z_1z_2\},\quad
						\{x_1x_2z_1,y_1y_2z_2\},\quad\{x_1x_2z_2,y_2z_1z_2\},
					\]
				such that each of them, together with edges $f_u$ and $f_w$, forms a minimal 4-path in $\cH$ (see Figure \ref{fig:f02}). Consequently, from each of these pairs only one edge may belong to $\cH_Z$.
			\end{proof}
			
			\begin{fact}\label{f457}
				If $ux_1x_2\in \cH^2$ and $|H^0(w)|\ge 2$, then $|\cH^1_Z|\le 8$, $|\cH^2_Z|\le 2$ and $|\cH_Z|\le 12$.
			\end{fact}
			\begin{proof}
				Let $f_u=ux_1x_2\in \cH^2$ and $e,e'\in H^0(w)$. Regardless of whether $e\cap e'=\emptyset$ or not, every $f\in \cH^1_Z$ intersects at least one of  $e$ or $e'$.  Suppose that there is $f\in \cH_Z$, disjoint from exactly one of the edges $e$ and $e'$, say $e$. Then $f_u$, $f$, $e'w$, and $we$ form a minimal 4-path in $\cH$, a contradiction.  Since there are exactly two edges of $\cH^1_Z$ disjoint from $e$ and two other edges of $\cH^1_Z$ disjoint from $e'$, we have $|\cH^1_Z|\le12-4= 8$. In view of Fact \ref{f7}, this implies that
					\[
						|\cH_Z| = |\cH_Z^1| + |\cH_Z^0| + |\cH_Z^2| \le 8 + 4 = 12.
					\]

				Similarly, there exists in $\cH^2_Z$ at least one edge  intersecting $e$ and disjoint from $e'$ and at least one edge intersecting $e'$ and disjoint from $e$, implying $|\cH^2_Z|\le 4-2= 2$.
			\end{proof}
						
			\begin{fact}\label{f1}
				If $|H^0(u)|\ge 3$ and $|H^0(w)|\ge 3$, then $|\cH_Z|\le 13$.
			\end{fact}
			\begin{proof}
				Observe that, if $e\in H^0(w)$ and $e'\in H^0(u)\cap H^0(w)$ are two disjoint edges, then there is no $f\in \cH_Z$ with $f\cap e'=\emptyset$, because otherwise edges $f$, $ew$, $we'$ and $e'u$ would form a minimal 4-path in $\cH$. As there are exactly four triples in ${Z \choose 3}$ disjoint from $e'$, the presence of such $e,e'$ eliminates 4 edges from $\cH_Z$ (see Figure~\ref{fig:f8}$(a),(b)$).
				
				Further note, that since $|H^0(u)|\ge 3$, $|H^0(w)|\ge 3$, and both $H^0(u)$, $H^0(w)$ are subgraphs of the cycle $y_1z_1y_2z_2$, there are at least two edges $e,e'\in H^0(u)\cap H^0(w)$. If $e\cap e'=\emptyset$, then every triple in ${Z\choose 3}$ disjoint from $e'$ intersects $e$ and vice versa (see Figure~\ref{fig:f8}$(c)$). Therefore, $|\cH_Z|\le 20-2\cdot 4=12$, better than needed.
				
				Otherwise $e$ and $e'$ share a vertex, and there are two further edges $\hat e, \hat e'\in H^0(u)\cup H^0(w)$, such that $e\cap \hat e=\emptyset$ and $e'\cap \hat e'=\emptyset$ (see Figure \ref{fig:f8}$(d)$). Hence, we can apply the above elimination scheme to these two pairs. As there is exactly one triple in ${Z\choose 3}$ disjoint from both $e$ and $e'$, by sieve principle, we eliminate from $\cH_Z$ exactly $4+4-1=7$ edges, leading to the required bound $|\cH_Z|\le 20-7=13$.
			\end{proof} 	
				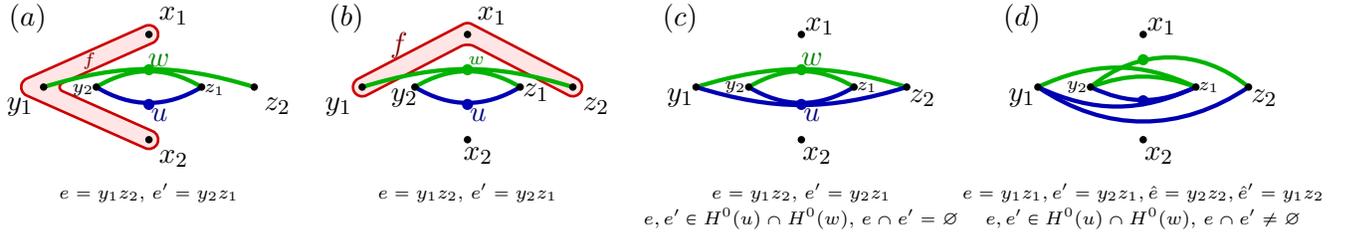
\begin{figure}[h!]
					\centering
					\begin{multicols}{4}
						\begin{tikzpicture}[scale=.7]
						
						\coordinate (x2) at (0,0);
						\coordinate (x1) at (0,2);
						
						\coordinate (y1) at (-2,1);
						\coordinate (y2) at (-1,1);						
						\coordinate (z1) at (1,1);
						\coordinate (z2) at (2,1);
						
						\coordinate (u) at (0,.67);
						\coordinate (w) at (0,1.33);
						
						\qqedge{(x2)}{($(y1)+(-0.25,0)$)}{(x1)}{5pt}{1pt}{red!80!black}{red!50!white,opacity=0.2}{15pt};
						\draw (y1) edge [bend left =17, color=green!70!black, ultra thick] (z2);
						\draw (y2) edge [bend left =31, color=green!70!black, ultra thick] (z1);
						\draw (y2) edge [bend right = 31, color=blue!70!black, ultra thick] (z1);
						
						\foreach \i in {x1,x2,y1,y2,z1,z2}
						\fill  (\i) circle (2pt);
						
						\fill (u) [blue!70!black] circle (3pt);
						\fill (w)[green!70!black] circle (3pt);
						
						\node at (x2)[below right] {$x_2$};
						\node at (x1)[above right] {$x_1$};
						
						\node at (y1)[below left] {$y_1$};
						\node at (-1.25,.95) {{\tiny $y_2$}};
						
						\node at (1.25,.95) {{\tiny $z_1$}};
						\node at (z2)[below right] {$z_2$};
						
						\node at ($(u)+(.2,-.2)$) [blue!50!black] {$u$};
						\node at ($(w)+(.19,.19)$) [green!50!black] {$w$};
						
						\node at (-1.13,1.5) [red!50!black]{{\tiny $f$}};
						
						\node at (0,-1){\tiny $e=y_1z_2$, $e'=y_2z_1$};
						
						\node at (-2.3,2.3){$(a)$};

						\end{tikzpicture}\\
						
						\begin{tikzpicture}[scale=.7]
						\qqedge{(y1)}{($(x1)+(0,.04)$)}{(z2)}{5pt}{1pt}{red!80!black}{red!50!white,opacity=0.2}{5pt};
						\draw (y1) edge [bend left =17, color=green!70!black, ultra thick] (z2);
						\draw (y2) edge [bend left =31, color=green!70!black, ultra thick] (z1);
						\draw (y2) edge [bend right = 31, color=blue!70!black, ultra thick] (z1);
						
						\foreach \i in {x1,x2,y1,y2,z1,z2}
						\fill  (\i) circle (2pt);
						
						\fill (u) [blue!70!black] circle (3pt);
						\fill (w)[green!70!black] circle (3pt);
						
						\node at (0.2,-0.3) {$x_2$};
						\node at (x1)[above right] {$x_1$};
						
						\node at (y1)[below left] {$y_1$};
						\node at (-1.2,.8) {$y_2$};
						
						\node at (1.3,.9) {$z_1$};
						\node at (z2)[below right] {$z_2$};
						
						\node at ($(u)+(.2,-.2)$) [blue!50!black] {$u$};
						\node at ($(w)+(.17,.17)$) [green!50!black] {{\tiny $w$}};
						
						\node at (-1.3,1.8) [red!50!black]{{$f$}};
						
						\node at (0,-1){\tiny $e=y_1z_2$, $e'=y_2z_1$};
						\node at (-2.3,2.3){$(b)$};
						\end{tikzpicture}\\
						
						\begin{tikzpicture}[scale=.7]

						\draw (y1) edge [bend left =17, color=green!70!black, ultra thick] (z2);
						\draw (y2) edge [bend left =31, color=green!70!black, ultra thick] (z1);
						\draw (y2) edge [bend right = 31, color=blue!70!black, ultra thick] (z1);
						\draw (y1) edge [bend right =17, color=blue!70!black, ultra thick] (z2);
						
						\foreach \i in {x1,x2,y1,y2,z1,z2}
						\fill  (\i) circle (2pt);
						
						\fill (u) [blue!70!black] circle (3pt);
						\fill (w)[green!70!black] circle (3pt);
						
						\node at (0.3,-0.3) {$x_2$};
						\node at ($(x1)+(.35,.15)$) {$x_1$};
						
						\node at ($(y1)+(-.3,-.2)$) {$y_1$};
						\node at (-1.25,1) {{\tiny $y_2$}};
						
						\node at (1.25,1) {{\tiny $z_1$}};
						\node at ($(z2)+(.3,-.2)$) {$z_2$};
						
						\node at ($(u)+(.2,-.2)$) [blue!50!black] {$u$};
						\node at ($(w)+(.2,.2)$) [green!50!black] {$w$};
						
						\node at (0,-1){\tiny $e=y_1z_2$, $e'=y_2z_1$} ;
						
						\node at (0,-1.5) {\tiny $e,e'\in H^0(u)\cap H^0(w)$, $e\cap e'=\emptyset$};
						\node at (-2.3,2.3){$(c)$};
						
						\end{tikzpicture}\\
						
						\begin{tikzpicture}[scale=.7]
						
						\coordinate (u) at (0,.74);
						\coordinate (w) at (0,1.52);
						
						\draw (y2) edge [bend left =40, color=green!70!black, ultra thick] (z2);
						\draw (y2) edge [bend left =20, color=green!70!black, ultra thick] (z1);
						\draw (y1) edge [bend left =25, color=green!70!black, ultra thick] (z1);
						\draw (y2) edge [bend right = 25, color=blue!70!black, ultra thick] (z1);
						\draw (y1) edge [bend right =25, color=blue!70!black, ultra thick] (z1);
						\draw (y1) edge [bend right =34, color=blue!70!black, ultra thick] (z2);
						
						\foreach \i in {x1,x2,y1,y2,z1,z2}
						\fill  (\i) circle (2pt);
						
						\fill (u) [blue!70!black] circle (3pt);					
						\fill (w)[green!70!black] circle (3pt);
						
						\node at (0.3,-0.3) {$x_2$};
						\node at ($(x1)+(.35,.15)$) {$x_1$};
						
						\node at ($(y1)+(-.3,-.2)$) {$y_1$};
						\node at (-1.25,1) {{\tiny $y_2$}};
						
						\node at (1.25,1) {{\tiny $z_1$}};
						\node at ($(z2)+(.3,-.2)$) {$z_2$};
						
						\node at (0,-1){\tiny $e=y_1z_1, e'=y_2z_1, \hat e = y_2z_2, \hat e'=y_1z_2$} ;
						
						\node at (0,-1.5) {\tiny $e,e'\in H^0(u)\cap H^0(w)$, $e\cap e'\neq\emptyset$};
						
						\node at (-2.3,2.3){$(d)$};						
						\end{tikzpicture}
						
					\end{multicols}
					
					\caption{$(a), (b)$: a minimal 4-path $\{f,ew,we',e'u\}$; $(c), (d)$: $e,e'\in H^0(u)\cap H^0(w)$.}
					\label{fig:f8}
					\vspace{-1em}
				\end{figure}

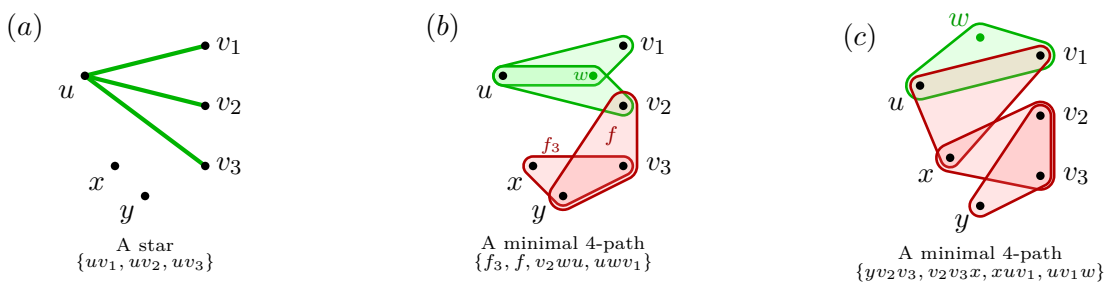
\begin{figure}[h!]
				\centering
				\begin{multicols}{3}
					\begin{tikzpicture}[scale=.8]
						
					\coordinate (u) at (0,2.5);
					
					\coordinate (v1) at (2,3);		
					\coordinate (v2) at (2,2);
					\coordinate (v3) at (2,1);
					
					\coordinate (x) at (.5,1);
					\coordinate (y) at (1,.5);
					
					\coordinate (w) at (1.5,2.5);

					\draw (u) edge [color=green!70!black, ultra thick] (v1);
					\draw (u) edge [color=green!70!black, ultra thick] (v2);
					\draw (u) edge [color=green!70!black, ultra thick] (v3);
					
					\foreach \i in {u,v1,v2,v3,x,y}
					\fill  (\i) circle (2pt);
							\node at (u)[below left] {$u$};
					\node at ($(x)+(-.3,-.3)$) {$x$};		
					\node at ($(y)+(-.3,-.3)$) {$y$};
					
					\node at ($(v1)+(.4,0)$){$v_1$};
					\node at ($(v2)+(.4,0)$){$v_2$};
					\node at ($(v3)+(.4,0)$){$v_3$};	
						
					\node at (1,-.3) {\tiny A star};
					\node at (1,-.6){\tiny $\{uv_1, uv_2, uv_3\}$};
					
					\node at (-1,3.3){$(a)$};

					\end{tikzpicture}\\
					
					\begin{tikzpicture}[scale=.8]
						\qedge{(w)}{(u)}{(v1)}{4.5pt}{1pt}{green!70!black}{green!50!white,opacity=0.2};
					\qedge{(u)}{(w)}{(v2)}{4.5pt}{1pt}{green!70!black}{green!50!white,opacity=0.2};
						\qedge{(y)}{(x)}{(v3)}{4.5pt}{1pt}{red!70!black}{red!50!white,opacity=0.2};
					\qedge{(y)}{(v2)}{(v3)}{6.5pt}{1pt}{red!70!black}{red!50!white,opacity=0.2};
						\foreach \i in {u,v1,v2,v3,x,y}
					\fill  (\i) circle (2pt);
					
					\fill (w) [green!70!black] circle (2pt);
					
					\node at (u)[below left] {$u$};
					\node at ($(x)+(-.3,-.3)$) {$x$};		
					\node at ($(y)+(-.4,-.3)$) {$y$};
					
					\node at ($(v1)+(.5,0)$){$v_1$};
					\node at ($(v2)+(.6,0)$){$v_2$};
					\node at ($(v3)+(.6,0)$){$v_3$};	
					
					\node at ($(x)+(.3,.35)$) [red!60!black] {\tiny $f_3$};
					
					\node at ($(v3)+(-.2,.5)$) [red!60!black] {{\scriptsize $f$}};
					
					\node at ($(w)+(-.2,0)$)[green!60!black]{\tiny $w$};	
					
					\node at (1,-.3) {\tiny A minimal 4-path};
					\node at (1,-.6){\tiny $\{f_3,f, v_2wu,uwv_1\}$};
					
					\node at (-1,3.3){$(b)$};
				
					\end{tikzpicture}\\
					\begin{tikzpicture}[scale=.8]

					\coordinate (w) at (1,3.3);
					\coordinate (x) at (.5,1.3);
					
					\qedge{(u)}{(w)}{(v1)}{6.5pt}{1pt}{green!70!black}{green!50!white,opacity=0.2};
										\qedge{(x)}{(u)}{(v1)}{4.5pt}{1pt}{red!70!black}{red!50!white,opacity=0.2};
					\qedge{(x)}{(v2)}{(v3)}{6.5pt}{1pt}{red!70!black}{red!50!white,opacity=0.2};
					\qedge{(y)}{(v2)}{(v3)}{5pt}{1pt}{red!70!black}{red!50!white,opacity=0.2};
						\foreach \i in {u,v1,v2,v3,x,y}
					\fill  (\i) circle (2pt);
					
					\fill (w) [green!70!black] circle (2pt);
					
					\node at ($(u)+(-.4,-.3)$) {$u$};
					\node at ($(x)+(-.4,-.3)$) {$x$};		
					\node at ($(y)+(-.3,-.3)$) {$y$};
					
					\node at ($(v1)+(.6,0)$){$v_1$};
				\node at ($(v2)+(.6,0)$){$v_2$};
					\node at ($(v3)+(.6,0)$){$v_3$};	
					\node at ($(w)+(-.35,.3)$) [green!60!black] {$w$};
					
					\node at (1,-.3) {\tiny A minimal 4-path};
					\node at (1,-.6){\tiny $\{yv_2v_3, v_2v_3x, xuv_1,uv_1w\}$};
					
					\node at (-1,3.3){$(c)$};

					\end{tikzpicture}
					
				\end{multicols}
				\caption{$(a)$: a star in $H(w)$; $(b)$, $(c)$: minimal 4-paths in $\cH$.}
				\label{fig:f9}
				\vspace{-1em}
			\end{figure}

			\begin{fact}\label{f9}
%				If $\Delta(L_{\cH}(w))\ge 3$, then $|\cH_Z|\le 14$.
				If $S_4^{(2)}\subseteq L_{\cH}(w) $, then $|\cH_Z|\le 14$.
			\end{fact}
		\begin{proof}
		Recall that $L_{\cH}(w)\subseteq \binom Z2$.
				Let $\{u v_1, uv_2,uv_3\}$ be in $L_{\cH}(w)$ and set $Z\setminus \{u,v_1,v_2, v_3\}=\{x,y\}$ (see Figure~\ref{fig:f9}$(a)$). If for some $i\in [3]$, $f_i=xyv_i\in \cH_Z$, then  none of the six triples $f\subset Z$, such that $u\notin f$ and $|f\cap f_i|=2$, can belong to $\cH_Z$, since otherwise the edges $f_i,f, uwv_j, uwv_k$, $\{i,j,k\}=\{1,2,3\}$, would form a minimal 4-path, contradicting the $\p$-freeness of $\cH$. Thus, $|\cH_Z|\le 20-6=14$ (see Figure \ref{fig:f9}$(b)$).
			
			Therefore, assume now that $f_1, f_2, f_3\notin \cH_Z$. Again by the $\p$-freeness of $\cH$, from each of the three disjoint sets of triples,
				\[
					\{uv_1x,v_2v_3x, v_2v_3y\},\quad \{uv_2x,v_1v_3x, v_1v_3y\},\quad\{uv_3x, v_1v_2x, v_1v_2y\},
				\]
			at most two triples may belong to $\cH_Z$ and, consequently, $|\cH_Z|\le 20 - 3 - 3=14$ (see Figure \ref{fig:f9}$(c)$).
		\end{proof}
			
		\begin{fact}\label{f3}
			If $L_{\cH}(w)=S_5^{(2)}$ is a star with the center in $\{y_1, y_2, z_1, z_2\}$ and $|\cH_Z|\ge 14$, then $H\subseteq \sk_n$.
		\end{fact}
		\begin{proof}
			Without loss of generality we may assume that $y_1$ is the center of the star $L_{\cH}(w)$. Thus, $L_{\cH}(w)=\{y_1v\colon v\in A\}$, where $A=\{x_1,x_2,z_1,z_2\}$. Let us denote by $K_A$ the complete 3-graph on  $A$. We will prove that
				\[
						\cH\subseteq K_A \cup S(y_1,A) = SK_n,
				\]
			which boils down to showing that for each edge $f\in \cH$  with  $f\nsubseteq A$ we have  $f\cap A\neq\emptyset$ and $y_1\in f$.
						
			Recall that each $f\in \cH$ satisfies $|f\cap Z|\ge 2$ and $f\cap Z \neq \{y_1,y_2\}$. Therefore, for all $f\in \cH$, we have $f\cap A \neq\emptyset$. Consequently, we only need to show that if $f\nsubseteq A$, then $y_1\in f$.
			
			Let us begin with $f\in \cH_Z$. By Fact \ref{f00}\ref{it:f025},  for all  $f\in (H^0_Z\cup \cH^1_Z)\setminus \cc$ we have $f\cap\{x_1,y_1\}\neq\emptyset$ and $f\cap\{x_2,y_1\}\neq\emptyset$, and thus, $y_1\in f$. So, we are done with $\cH_Z$, except that we still need to rule out the presence of the edge $x_1x_2y_2$ in $\cH$.
			
			The above established fact that $y_1\in f$ for all $f\in(\cH^0_Z\cup \cH^1_Z)\setminus \cc$ implies that $|\cH^0_Z|\le 3$ and $|\cH^1_Z|\le 8$, and, in turn, $|\cH^2_Z|=|\cH_Z|-|\cH^0_Z|-|\cH_Z^1|\ge 14-3-8= 3$. But triples $x_1x_2y_2$, $x_1x_2z_1$, $wy_1z_1$, and $wy_1z_2$ form a minimal 4-path, and the same is true with $x_1x_2z_1$ replaced by $x_1x_2z_2$ and the last two edges reversed. Thus, in order to satisfy $|\cH^2_Z|\ge 3$, we must have $x_1x_2y_2\notin \cH_Z^2$, while $x_1x_2z_1, x_1x_2z_2\in \cH_Z^2$.
			
			Turning to the edges of $\cH\setminus\cH_Z$, recall that all edges of $\cH$ containing $w$ contain also $y_1$. Next, fix an arbitrary vertex $u\in V\setminus (Z\cup \{w\})$, and observe, that $x_1x_2z_1\in \cH_Z^2$ entails $ux_1x_2\notin \cH$, and thus $H^2(u)=\emptyset$, because otherwise $\cH$ would contain a minimal 4-path consisting of edges $ux_1x_2$, $x_1x_2z_1$, $wy_1z_1$, and $wy_1z_2$.
			Finally, note that, by Fact \ref{f00}\ref{it:f02}, all edges of $H^0(u)$ intersect $\{x_1,y_1\}$ and $\{x_2,y_1\}$, while all  edges of $H^1(u)$ intersect $\{y_1,z_1\}$ and $\{y_1,z_2\}$. This implies that for all $e\in L_{\cH}(u)=H^0(u)\cup H^1(u)$, the condition $y_1\in e$ holds. In summary, for all $f\in \cH\setminus\cH_Z$, we have $y_1\in f$, which ends the proof.
			\end{proof}

		\subsection{Proofs of Lemmas \ref{l:sp} and \ref{l:sk}.}
			
			\begin{proof}[Proof of Lemma \ref{l:sp}]
				Assume, for the sake of a contradiction, that the assumptions of Lemma \ref{l:sp} are satisfied, but $\cH\nsubseteq \sp_n$. Let $u,w\in W$ be two vertices with degree in $\cH$ at least $5$. In view of Corollary~\ref{minmax}, we actually have
				\[
				\deg_\cH(u)=\deg_\cH(w) = 5.
				\]
				Then Corollary \ref{f2} combined with \eqref{eq:h02} and Fact \ref{f0} tells us that this is possible only if $|H^2(u)|=|H^2(w)|=1$ and one of the following is true:
				\begin{enumerate}[label=\rmlabel]
					\item \label{it:l11} $|H^0(u)|=|H^0(w)|=4$;
					\item \label{it:l12} $|H^1(u)|=|H^1(w)|=4$;
					\item \label{it:l13} $|H^0(u)|=|H^1(u)|=|H^0(w)|=|H^1(w)|=2$.
				\end{enumerate}			
				Case \ref{it:l11} is impossible -- otherwise, the vertices $y_1z_2uy_2z_1wx_1x_2$ would form a minimal 4-path in $\cH$.%
				
				If we are in case \ref{it:l12}, then because $\cH\nsubseteq \sp_n$, Fact \ref{f0} together with Fact \ref{f00}\ref{it:f07}, ensures that both $H^1(u)$ and $H^1(w)$ are 4-cycles, either $C_y=x_1y_1x_2y_2$ or $C_z=x_1z_1x_2z_2$. Now, Fact \ref{f00}\ref{it:f01} entails, that exactly one of them, say $H^1(u)$, equals $C_y$, whereas the other one $H^1(w)=C_z$. But then $|H^1(u)\cup H^1(w)|\ge 7$, and thus, in view of Fact \ref{f00}\ref{it:f066}, $|\cH_Z|\le 8$ yielding
				\[
				|\cH[Z\cup \{u,w\}]| = |\cH_Z| + \deg_\cH(u) + \deg_\cH(w)\le 8+5+5=18<22.
				\]

				Therefore there exists a vertex $v\in W\setminus \{u,w\}$, with $\deg_\cH(v)\ge 5$. Another application of Corollary~\ref{f2}\ref{it:f21} with $v$ in place of $u$ says, that $H^0(v)=\emptyset$, again by Facts \ref{f0} and \ref{f00}\ref{it:f07}, either $H^1(v)=C_y$ or $H^1(v)=C_z$. But, because already $H^1(u)=C_y$ and $H^1(w)=C_z$, in view of Fact \ref{f00}\ref{it:f01} this is impossible, namely, we arrive at a contradiction.
				
				Finally, in  case \ref{it:l13}, one can observe that, by Fact \ref{f457}, $|\cH^1_Z|\le 8$ and $|\cH^2_Z|\le 2$, while, by Fact \ref{f8}, $\cH_Z^0=\emptyset$ and for every $v\in W\setminus \{u,w\}$, $H^0(v)=\emptyset$.
				Altogether, we get  $|\cH_Z|\le 10$ and, consequently, $|\cH[Z\cup \{u,w\}]| \le 20$. Hence there is a vertex $v\in W\setminus \{u,w\}$ with $\deg_\cH(v)\ge 5$. Now, Corollary~\ref{f2}\ref{it:f21} says $|H^1(v)|\le 2$, and so $\deg_\cH(v)=|H^0(v)|+|H^1(v)|+|H^2(v)|\le 3$, yielding a contradiction with $\deg_\cH(v)\ge 5$.
			\end{proof}

		\begin{proof}[Proof of Lemma \ref{l:sk}]
			Let $u,w\in W$ be two vertices, such that
				\begin{enumerate}[label=\rmlabel]
					\item\label{it:2a} $|\cH[Z\cup \{u,w\}]|\ge 22$ and
					\item\label{it:2b} $\deg_\cH(w)\le 4$.
				\end{enumerate}
			We will show that $\cH\subseteq SK_n$, which will end the proof. Set $\hat \cH=\cH[Z\cup\{u,w\}]$.
			
			Because $L_\cH(u)\subseteq \binom Z2$, the connectivity of $\cH$ implies $\cH[Z\cup \{u\}]\neq K^{(3)}_6 \cup K_1$, and thereby, in view of Lemma \ref{l:n7}, $|\cH[Z\cup \{u\}]|\le 19$. Consequently, $\deg_\cH(w)\ge 3$, yielding that at least one of the graphs, $H^0(w)$ or $H^1(w)$, is not empty. Hence, by \eqref{eq:h02}, Fact \ref{f0} and Corollary \ref{f2}, $\deg_\cH(u)\le 4+2+1=7$. Similarly, $\deg_\cH(u)\ge 3$.
			
			Suppose that $\deg_\cH(u)\ge 6$. Then, in view of the bound $\deg_\cH(w)\ge 3$, Corollary \ref{minmax} tells us that $\deg_\cH(w) = 3$ and so $H^2(w)\neq\emptyset$. In addition, as $|H^0(u)|+|H^1(u)|\ge 5$, either $|H^0(u)|\ge 2$ or $|H^1(u)|= 4$, implying, together with Facts \ref{f457} (with $u$ and $w$ swapped) and \ref{f00}\ref{it:f06}, that $|\cH_Z|\le 12$. Therefore \ref{it:2a} entails, that $\deg_\cH(u)=7$ which, in turn, results $|H^0(u)|\ge 2$ and $H^2(u)\neq\emptyset$. But then Facts \ref{f457} and \ref{f8} yield that $|\cH_Z^1|\le 8$, $|\cH_Z^2|\le 2$, and $\cH^0_Z=\emptyset$.  Consequently,
				\[
					|\hat \cH|=|\cH^0_Z|+|\cH^1_Z|+|\cH^2_Z|+\deg_\cH(u)+\deg_\cH(w)\le 0 + 8 + 2 + 7 + 3 =20,
				\]
			contradicting \ref{it:2a}.
			
			Hence, from now on, we assume that $\deg_\cH(u)\le 5$. Then, in view of \ref{it:2a} and \ref{it:2b}, it follows that $|\cH_Z|\ge 22-5-4= 13$, implying, via Fact \ref{f00}\ref{it:f06}, that both $|H^1(u)|\le 3$ and $|H^1(w)|\le 3$. We split the proof into three cases according to the emptiness  of  $H^2(u)$ and $H^2(w)$. In particular we will show that if at least one of these graphs is not empty, then $|\hat \cH|\le 21$, contradicting \ref{it:2a}.
			
			\medskip
			
				\noindent{\it Case 1. }{\bf $\mathbf{H^2(u)\neq \emptyset}$ and $\mathbf{H^2(w)\neq \emptyset}$}.  If, in addition, $H^0(u)=H^0(w)=\emptyset$, then either $|H^1(u)\cup H^1(w)|=2$ and so, by Fact \ref{f00}\ref{it:f04},
				\[
					|\hat \cH| = |\cH_Z| +\deg_\cH(u)+\deg_\cH(w)\le 15 + 3 + 3=21,
				\]
				or  $|H^1(u)\cup H^1(w)|\ge 3$. Then, in view of Fact \ref{f00}\ref{it:f05}, $|\cH_Z|\le 13$ and, again,  $|\hat \cH| \le 13 + 4 + 4=21$.
				
				Therefore we may assume, that $H^0(u)\cup H^0(w)\neq\emptyset$ yielding, together with Fact \ref{f8}, $\cH^0_Z=\emptyset$. Moreover,  since $|\cH_Z|\ge 13$, Fact \ref{f457} tells us that both $|H^0(u)|\le 1$ and $|H^0(w)|\le 1$. Finally, by Fact~\ref{f00}\ref{it:f03}-\ref{it:f05}, either $|H^1(u)\cup H^1(w)|=1$ and thus $|\cH^1_Z|\le 10$, $|H^1(u)\cup H^1(w)|=2$, entailing $|\cH^1_Z|\le 9$, or $|H^1(u)\cup H^1(w)|\ge 3$ and then $|H^1_Z|\le 8$.	
				That is, $|\cH_Z^1|+|H^1(u)|+|H^1(w)|\le 14$ and the equality holds only if $|H^1(u)|=|H^1(w)|=3$. Altogether, in all of these cases, as $|\cH_Z^2|\le 4$, and $|H^0(u)|+|H^2(u)|+|H^0(w)|+|H^2(w)|\le 4$,
					\[
						|\hat \cH| = |\cH^0_Z| + |\cH^2_Z| + (|\cH_Z^1|+|H^1(u)|+|H^1(w)|) + |H^0(u)|+|H^2(u)|+|H^0(w)|+|H^2(w)|  \le 21,
					\]
				unless $|H^1(u)|=|H^1(w)|=3$, in which case $|\hat\cH|\le 21$ by using $\deg_\cH(u)+\deg_\cH(w)\le 9$ and $|H^1_Z|\le 8$.

				\noindent{\it Case 2. }{\bf $\mathbf{H^2(u)\neq\emptyset}$ and $\mathbf{H^2(w)=\emptyset}$} (the proof of the case $H^2(u)=\emptyset$ and $H^2(w)\neq\emptyset$ is similar). Recall, that $|H^1(w)|\le 3$ and $|\cH_Z|\ge 13$ which implies, together with Fact \ref{f457}, that $|H^0(w)|\le 1$. Therefore, $\deg_\cH(w)\le 3$, because otherwise, $|H^1(w)|=3$ and $|H^0(w)|=1$. But then, by Facts \ref{f00}\ref{it:f05} and \ref{f7}, $|\cH^1_Z|\le 8$ and $|\cH_Z^0|+|\cH_Z^2|\le 4$, contradicting $|\cH_Z|\ge 13$. Hence, $\deg_\cH(u)+\deg_\cH(w)\le 8$  from which we infer that $|\cH_Z|\ge 14$ and, consequently, by Fact \ref{f00}\ref{it:f05}, $|H^1(w)|\le 2$. Thus, $|H^0(w)|=1$ and $|H^1(w)|=2$. But then, again by Facts \ref{f00}\ref{it:f04} and \ref{f7}, $|\cH_Z| \le 9 + 4 < 14$, a contradiction.				
				\medskip
				
				\noindent {\it Case 3. }$\mathbf{H^2(u)=H^2(w)=\emptyset}$. First observe, that \eqref{eq:h02}, Fact \ref{f0}, and Corollary \ref{f2} tell us ${\deg_\cH(u),\deg_\cH(w)\le 4}$ and, consequently, \ref{it:2a} yields $|\cH_Z|\ge 14$. Thus, by Fact \ref{f00}\ref{it:f05},
					\[
						|H^1(u)\cup H^1(w)|\le 2.
					\]
				Note also that both
					\[	
						|H^0(u)|\le 2\quad\mbox{ and }\quad|H^0(w)|\le 2,
					\]
				 because otherwise Corollary \ref{f2} and $\deg_\cH(u), \deg_\cH(w)\ge 3$ entail, that $|H^0(u)|\ge 3$ and $|H^0(w)|\ge 3$. This, however, together with Fact \ref{f1} implies $|\cH_Z|\le 13$, a contradiction.
								
				 Now, in view of Fact \ref{f3}, to finish the proof it is enough to show that at least one of the graphs $L_\cH(u)$ or $L_\cH(w)$, is a star $S_5^{(2)}$ with the center in $\{y_1, y_2, z_1, z_2\}$. To this end observe, that if $|H^1(u)\cup H^1(w)|= 2$ and either $|H^0(u)|=|H^0(w)|=1$ or $|H^1(u)|=|H^1(w)|=1$, then $\deg_\cH(u)=\deg_\cH(w)=3$ and, in view of Fact~\ref{f00}\ref{it:f04}, $|\cH_Z|\le15$, yielding $|\hat \cH|\le 15+3+3= 21$, a contradiction with \ref{it:2a}.
				
				 Otherwise there exists an edge $e\in H^1(u)\cap H^1(w)$, and at least one of the graphs, $H^0(u)$ or $H^0(w)$, say $H^0(w)$, has two edges. We let $\{v\}=e\cap \{y_1,y_2,z_1, z_2\}$ and note, that due to Fact \ref{f00}\ref{it:f02}, every edge of $H^0(u)\cup H^0(w)$ contains $v$. Therefore $S^{(2)}_4\subseteq L_\cH(w)$ entailing, together with Fact \ref{f9}, $|\cH_Z|\le 14$, and thereby $\deg_\cH(u)=\deg_\cH(w)=4$. In particular, $H^0(w)$ has two edges both containing $v$. Finally, a repeated application of Fact \ref{f00}\ref{it:f02} reviles, that $L_\cH(u)$ is a star $S^{(2)}_5$ with the center $v$, as required.
			\end{proof}

%%%%%%%%%%%%%%%%%%%%%%%%%%%%%%%%%%%%%%%%%%%%%%%%%%%%%%%%%%%%%%%%%%%%%%
%																																		%
%																																		%
%											Ramsey numbers																%
%																																		%
%																																		%
%%%%%%%%%%%%%%%%%%%%%%%%%%%%%%%%%%%%%%%%%%%%%%%%%%%%%%%%%%%%%%%%%%%%%%

\section{Ramsey numbers}\label{RamNum}

\subsection{Shorter paths} Before we turn to proving Theorem \ref{RN}, let us briefly discuss Ramsey numbers for 3-uniform minimal paths of shorter length.
Observe that the family $\cP_2$ consists of two 3-graphs, each being a pair of overlapping edges, either in one (a bow) or two vertices (a kite). Therefore, $\cP_2$-free 3-graphs are necessarily matchings, that is, consist of disjoint edges only. Consequently, $\ex_3(n;\cP_2)=\lfloor n/3\rfloor$, and, by \eqref{RamTur}, 
\[
R(\cP_2;r)=\min\left\{n: \frac{\binom n3}{\lfloor n/3\rfloor}>r\right\},
\]
or, asymptotically, $R(\cP_2;r)\sim\sqrt{2r}$, as $r\to\infty$. For small $r$, in particular, $R(\cP_2;2)=R(\cP_2;3)=4$, while $R(\cP_2;4)=5$.
In \cite{Axe}, the two 3-graphs belonging to $\cP_2$ were considered separately. It was shown there that $R(bow;r)\sim\sqrt{6r}$, while $R(kite;r)\in\{r+1,r+2,r+3\}$ depending on the divisibility of $r$ by 6. It is, perhaps, interesting to see the drop from $\sqrt{6r}$ to $\sqrt{2r}$ when the bow is accompanied by the kite.

The family $\cP_3$ also consists of two 3-graphs, among them the linear path $P_3$. For the latter, an easy lower bound by a construction of Gyarfas and Raeisi \cite{GR} says that   $R(P_3;r)\ge r+6$.
It was proved in a series of papers (\cite{Jac}, \cite{JPR2}, \cite{PR}, and \cite{Pol}) that, indeed, $R(P_3;r)= r+6$ for $r\le10$. The trivial upper bound, $R(P_3;r)\le 3r$, stemming from \eqref{RamTur} was improved down to $R(P_3;r)<1.98r$ in \cite{LP}.

Turning to minimal paths of length 3, there is a similar lower bound $R(\cP_3;r)\ge r+5$. Using the known value of $\ex_3(7;\cP_3)=15$ determined in \cite{MV}, it follows by \eqref{RamTur} that indeed $R(\cP_3;2)= 7$. With a bit more effort, observing that  a  connected $\cP_3$-free 3-graph must be intersecting and using the Hilton-Milner Theorem \ref{th:hm}, one can also show that $R(\cP_3;r)= r+5$ for $r\le7$. The range of $r$, for which $R(\cP_3;r)= r+5$ is certainly  wider, but to prove it one would need more sophisticated tools, like the third order Tur\'an number $\ex^{(3)}_3(n;\cP_3)$.
 %This number  is still not known for $r=4$; we only have
%$9\le R(\cP_3;4)\le R(\cP_3;4)\le R(P_3;4)=10$.

\subsection{Proof of Theorem \ref{RN}}
Let us start with a general lower bound on $R(\p;r)$ based on the slightly modified construction given by Gy\'{a}rf\'{a}s and Raeisi in \cite{GR}. We let
	\[
		s_r=\max\left\{s\in \mathbb{Z}: \sum_{k=6}^s\binom k2\le r-1\right\}\quad\mbox{and}\quad
		t_r=\max\left\{t\in\mathbb{Z}: \binom t3\le r\right\}.
	\]
\begin{prop}\label{lower}
	For all $ r\ge1$,
		\[
			R(\p;r) \ge r+\max\{s_r,t_r\}+1 \ge r +  \sqrt[3]{6r}+ 1.
		\]
\end{prop}

Note, that
	\[
	{\tiny s_r = \begin{cases}
	5, &\textrm{ for } 1\le r\le 15, \cr
	6, &\textrm{ for } 16\le r\le 36, \cr
	7, &\textrm{ for } 37\le r\le 64, \cr
	8, &\textrm{ for } 65\le r\le 100, \cr
	9, &\textrm{ for } 101\le r\le 145, \cr
	\cdots
	\end{cases},\qquad
	t_r = \begin{cases}
	3, &\textrm{ for } 1\le r\le 3, \cr
	4, &\textrm{ for } 4\le r\le 9, \cr
	5, &\textrm{ for } 10\le r\le 19, \cr
	6, &\textrm{ for } 20\le r\le 34, \cr
	7, &\textrm{ for } 35\le r\le 55, \cr
	\cdots
	\end{cases}, 		\quad\textrm{\normalsize and thus }\quad
	R(\p;r) \ge\begin{cases}
	r+6, &\textrm{ for } r\ge 1, \cr
	r+7, &\textrm{ for } r\ge 16, \cr
	r+8, &\textrm{ for } r\ge 35, \cr
	r+9, &\textrm{ for } r\ge 56, \cr
	r+10, &\textrm{ for } r\ge 84, \cr
	\cdots
	\end{cases}}
	\]
In particular, for $r\ge20$ we have $t_r\ge s_r$.

\begin{proof}
	Set $m=\max\{s_r,t_r\}$ and let $V(K^{(3)}_{r+m})=\{1,2,\dots,r+m\}$. If $m=s_r$, for $i=1,\dots, r-1$, color every edge of $K^{(3)}_{r+m}$ whose minimum vertex is $i$ by color $i$. In addition, apply different colors from $\{1,\dots,r-1\}$ to all edges with minimum vertex in the set $\{r,r+1,\dots, r+m-6\}$. Note that there are exactly $\sum_{k=6}^m\binom k2 \le r-1$ such edges. Moreover, the edges of color $i$ form a starplus, so no monochromatic copy of a  minimal 4-path  has been created in any of the first $r-1$ colors. The remaining uncolored edges form a complete 3-graph $K^{(3)}_6$ on the last 6 vertices $r+m-5,\dots,r+m$ and we color them by color $r$. As a minimal 4-path has at least 7 vertices, there is no member of $\p$ in color $r$ as well.

	If $m=t_r$, the construction is even simpler. For $i=1,\dots, r$,  color every edge of $K_{r+m}$ whose minimum vertex is $i$ by color $i$. In addition, apply different colors from $\{1,\dots,r\}$ to
	all $\binom m3\le r$ edges spanned on the vertices $r+1,\dots,r+m$. Again, each color is a starplus, so no monochromatic copy of a minimal 4-path  has been created.
\end{proof}

%In the rest of this section we will show that for $r\le 4$ the lower bound given in Proposition \ref{lower} is best possible, which will complete the proof of Theorem \ref{RN}.

\begin{proof}[Proof of Theorem \ref{RN}]
	In view of Proposition \ref{lower}, we only need to show the upper bound on $R(\p;r)$. For $r=1$ there is nothing to prove so let us begin with $r=2$ and $n=8$. For this purpose observe that in every 2-coloring of $K^{(3)}_8$ at least one color takes at least ${8\choose 3}/2=28>22=\ex(8;\p)$ edges, and so, due to Theorem~\ref{P4}, contains a member of $\p$. Moreover, the same averaging argument entails that this is true for every 8-vertex 3-graph with at least 45 edges.

	Now, let $r=3$ and $n=9$. With an eye on the case $r=4$, we are going to prove, for $r=3$, a slightly stronger result. An $r$-coloring which does not yield a monochromatic member of $\p$ is referred to as \emph{proper}. Let $\cH_9$ be a 9-vertex 3-graph with at least $\binom 93-2=82$ edges and let a proper 3-coloring of $\cH_9$ be given. Then, there is a color with at least $\lceil 82/3\rceil = 28>27=\ex^{(2)}(9;\p)$ edges and thus, since the coloring is proper,  by Theorems \ref{P4} and \ref{e2}, that color must be a subset of $S^{+1}_9$. After removing the center of that star as well as the unique edge not containing it, we obtain a proper 2-coloring of an 8-vertex 3-graph with at least $\binom 83-3=53$ edges, which, as it is shown above, contains a monochromatic member of $\p$, a contradiction.

	Finally, consider the case $r=4$ and $n=10$. To this end let a proper 4-coloring of all $\binom{10}3=120$ edges of $K^{(3)}_{10}$ be given. If there is a color which is a subset of either $S_{10}^{+1}$ or $SP_{10}$, then we remove its center together with at most two additional edges. As a result, we obtain a proper 3-coloring of a 9-vertex 3-graph with at least $\binom 93-2=82$ edges, which, as shown above, contains a monochromatic copy of a member of $\p$, a contradiction. Otherwise, in view of Theorems \ref{P4}, \ref{e2}, and \ref{e3}, each of the four colors has exactly 30 edges and is isomorphic to $SK_{10}$. But this is impossible, because in $K^{(3)}_{10}$ every vertex has degree ${9\choose 2}=36$, whereas in $SK_{10}$ each vertex has its degree in $\{4,11,26\}$. Clearly, $36$ can not be obtained as a sum of four numbers from $\{4,11,26\}$ and we are done.
\end{proof}

\section{Open problems}

It would be interesting, though tedious, to calculate higher order Tur\'an numbers for $\p$, that is, $\ex_3^{(s)}(n;\p)$, $s\ge4$, and, using them, to pin down Ramsey numbers $R(\p,r)$ for  $5\le r\le r_0$, for some $r_0\ge5$.

 Another  challenging project would be to determine \emph{for all $n$} the Tur\'an number $\ex_3(n;\C_4)$, where, recall $\C_4$ is the family of all minimal 3-uniform cycles with four edges. Kostochka, Mubayi, and Verstraete showed in \cite{KMV} that for \emph{large $n$}
\[
\ex_3(n;\C^3_4)=\binom{n-1}2+\left\lfloor\frac{n-1}3\right\rfloor.
\]
Gunderson, Polcyn, and Ruci\'nski in \cite{GPR} confirmed this formula  for $n\le7$.

Tur\'an numbers for longer minimal paths and cycles seem to be currently out of reach if one desires the exact values  \emph{for all $n$}.

\end{document}